\documentclass[12pt]{article}
\usepackage{amsmath}
\usepackage{amssymb}
\newcommand{\integ}[1]{\ensuremath{\int_{0}^{T}\!\!\int_{\Omega}{#1}
    \,dx\,dt}} 
\newcommand{\dif}[0]{\ensuremath{\mathrm{d}}}
\newtheorem{thm}{Theorem}
\newtheorem{lemma}[thm]{Lemma}
\newtheorem{cor}[thm]{Corollary}
\newtheorem{prop}[thm]{Proposition}
\newtheorem{definition}[thm]{Definition}
\DeclareMathOperator{\divergence}{div}
\DeclareMathOperator*{\essinf}{ess\,liminf}
\DeclareMathOperator{\esssup}{ess\,sup}

\newcommand{\Rn}{\ensuremath{\mathbb{R}^{n}}}
\begin{document}
\title{\textsc {Regularity of Supersolutions}}
\author{Peter Lindqvist}
\date{{\small Norwegian University of Science and Technology}}
\maketitle

\bigskip

{\footnotesize \textsf{These notes were written up after my lectures at
      Cetraro in June 2009 during the summer course ``Regularity estimates
      for nonlinear elliptic and parabolic problems'',  organized by
      C.I.M.E. in Italy.  Chapter 7 was written in 2011. New results were incorporated in a revision 2015. The topic is the sub-
      and supersolutions; they are like the stepchildren of Regularity
      Theory, since the proper solutions usually get most of the
      attention. Not now! My objective is the supersolutions of the
      $p$-Laplace Equation. The notes are a\, t\,o\,r\,s\,o: vital
      parts are missing. The fascinating story about the $p$-Laplace
      equation and its  solutions is not told here, the text
      being focused on supersolutions. 
      Generalizations to other
      equations are excluded. ''Less is more.''}}

\section{Introduction}
\label{introduction}
 
The regularity theory for  \emph{solutions} of certain parabolic
differential equations of the type 
\begin{equation} \label{Aeqn}
\frac{\partial u}{\partial t} =
 \divergence\textbf{A}(x,t,u,\!\nabla u)
\end{equation} 
is a well developed topic, but when it comes to (semicontinuous)
\emph{supersolutions} and \emph{subsolutions} a lot remains to be
done.  Supersolutions are often auxiliary tools as in the
celebrated Perron method, for example, but they are also interesting
in their own right. They appear as solutions to obstacle problems and
variational inequalities. 

As  a mnemonic rule
 $$\frac{\partial v}{\partial t} \geq \divergence\textbf{A}(x,t,v,\!\nabla
v)$$
for smooth supersolutions. We shall restrict our exposition mainly to the basic equations
$$\frac{\partial v}{\partial t} \geq \divergence (|\nabla v|^{p-2}\nabla v)\qquad \text{and}\qquad
\divergence (|\nabla u|^{p-2}\nabla u) \leq 0,$$
although the methods have a wider applicability. To avoid the splitting in different cases we usually keep
$p > 2.$

Our supersolutions are required to be
 lower semicontinuous but are not
assumed to be differentiable in any sense: part of the theory is to
prove that they have Sobolev derivatives. If one instead studies weak
supersolutions that by definition belong to a Sobolev space, then one
has the task to prove that they are semicontinuos. Unfortunately, the
weak supersolutions do not form a good closed class under monotone
convergence.  For \emph{bounded}
functions the definitions yield the same class of supersolutions.

 The modern theory of  viscosity solutions, created by Lions,
Crandall, Evans, Ishii, Jensen, and others, relies on the appropriately
defined \emph{viscosity supersolutions}, which are merely lower
semicontinuous functions by their definition. For second order
equations, these are often the same functions as those supersolutions
that are encountered in potential theory. The link enables one to study
the regularity properties also of the viscosity supersolutions. This is the case for the  so-called Evolutionary p-Laplace equation:
$$\frac{\partial v}{\partial t} = \divergence (|\nabla v|^{p-2}\nabla v).$$
We will restrict our exposition to this equation and we only treat the \textsf{slow diffusion case} $\mathbf{p>2}.$

To sum up, we shall deal with three different definitions of supersolutions:
\begin{itemize}
\item \textbf{Weak supersolutions}. They belong to the natural Sobolev space and satisfy the equation in weak form with test functions under the integral sign.
\item \textbf{Viscosity supersolutions}. The differential inequality is valid at points of contact for test functions touching from below.
\item \textbf{$p$-supercaloric functions}. They are defined as in Potential Theory via the comparison principle. 
\end{itemize}
The definitions are given later. As we will see in Chapter 7, the viscosity supersolutions and the $p$-supercaloric functions are the same functions. Therefore we often use the term ''viscosity supersolution'' as a label also for $p$-supercaloric functions.

As an example of what we have in mind, consider \textsf{the Laplace equation}
$\Delta u = 0$ and recall that a superharmonic
function is a lower semicontinuous function satisfying a comparison
principle with respect to the harmonic functions. An analogous
definition comes from the super meanvalue property. General
superharmonic functions are not differentiable in the classical
sense. Nonetheless, the following holds.
\begin{prop}
\label{Riesz}
Suppose that $v$ is a superharmonic function defined in $\mathbb{R}^{n}$. Then
the Sobolev derivative $\nabla v$ exists and
 $$ \int_{B_{R}} |\nabla v|^{q} \, dx  < \infty$$
whenever $0 < q < \frac{n}{n-1}$. Moreover,
$$\int_{\Rn} \langle \nabla v,\nabla \eta \rangle \,dx  \geq 0$$
for $\eta \geq 0,\,\eta \in C_{0}^{\infty}(\Rn)$.
\end{prop}

The fundamental solution $v(x) = |x|^{2-n}$ ($= -\log(|x|)$, when $n=2$)
is a superharmonic function showing that the summability exponent $q$ is
sharp. We seize the opportunity to mention that the superharmonic
functions are exactly the same as the viscosity supersolutions of the
Laplace equation. In other words, \emph{a viscosity supersolution has a
gradient in Sobolev's sense}. As an example, the Newtonian potential
$$v(x) = \sum \frac{c_{j}}{|x-q_{j}|^{n-2}},$$
where the rational points $q_{j}$ are numbered and the $c_{j}$'s are
convergence factors, is a superharmonic function, illustrating that
functions in the Sobolev space can be infinite in a dense set.
 ---The proof of the proposition follows
from Riesz's representation theorem, a classical result according to
which we have a
harmonic function plus a Newtonian potential. This was about the Laplace
equation.

A similar theorem holds for the viscosity supersolutions (= the
$p$-super\-har\-mo\-nic functions) of the
so-called \textsf{p-Laplace equation}
$$\nabla\!\cdot \!\left(|\nabla v|^{p-2}\nabla v \right) = 0$$
but now 0 $< q < \frac{n(p-1)}{n-1}$ in the counterpart to Proposition \ref{Riesz}. (Strictly speaking, we obtain a
proper 
Sobolev space only for $p > 2 - \frac{1}{n}$, because $q < 1$
otherwise.)
 The principle of
superposition is not valid and, in particular, Riesz's representation
theorem is no longer available. The original proof in [L1] was based on the
obstacle problem in the calculus of variations and on the so-called
weak Harnack inequality. At
present, the simplest proof seems to rely upon an approximation with
so-called infimal convolutions
$$v_{\varepsilon}(x) = \inf_{y} \left\{v(y) +
    \frac{|x-y|^{2}}{2\,\varepsilon}\right\}\:,\quad \varepsilon > 0.$$
At each point $v_{\varepsilon}(x) \nearrow v(x)$. They are viscosity
supersolutions, if the original $v$ is. Moreover, they are (locally)
Lipschitz continuous and hence differentiable a.e. Therefore the
approximants $v_{\varepsilon}$ satisfy expedient \emph{a priori} estimates, which, to
some extent, can be passed over to the original function $v$ itself.

Another kind of results is related to the \emph{pointwise behaviour}. The
viscosity supersolutions are pointwise defined. At \emph{each} point we have
$$v(x) = \essinf_{y \rightarrow x}v(y)$$
where \emph{essential limes inferior} means that sets of measure zero
are neglected in the calculation of the lower limit. In the linear
case $p = 2$ the result seems to be due to Brelot, cf. [B]. So much about the
p-Laplace equation for now. The theory extends to a wider class of
elliptic equations of the type
$$\divergence\textbf{A}(x,u,\!\nabla u)= 0.$$

For parabolic equations like
 $$\frac{\partial u}{\partial t} =
\sum_{i,j}\frac{\partial}{\partial x_{i}} \biggl(\Bigl|\sum_{k,m}
      a_{k,m}\frac{\partial u}{\partial x_{k}}\frac{\partial
        u}{\partial x_{m}}\Bigr|^{\frac{p-2}{2}}\!
    a_{i,j}\frac{\partial u}{\partial x_{j}}\biggr),$$
where the matrix $(a_{i,j})$ satisfies the
ellipticity condition $$ \sum a_{i,j}\xi_{i}\xi_{j} \geq \gamma
|\xi|^{2},$$ the situation is rather similar, although technically much
more demanding. Now the use of infimal convolutions as approximants
offers considerable simplification, at least in comparison with the
original proofs in [KL1]. We will study 
the \textsf{Evolutionary p-Laplace Equation}
\begin{equation}\frac{\partial u}{\partial t} = \divergence(|\nabla
  u|^{p-2}\nabla u) 
\end{equation}
where $u = u(x,t)$, restricting ourselves to \textsf{the slow diffusion
case $p > 2$.}

We shall encounter two different classes of supersolutions, depending on whether they belong to $L^{p-1}_{loc}(\Omega_T)$ or not. Depending on this, seemingly little distinction, the classes are widely apart. Those that belong to $L^{p-1}_{loc}(\Omega_T)$ have been much studied, since they have many good properties and satisfy a differential equation where the right-hand side is a Riesz measure. The others are less known. 

 The celebrated Barenblatt solution\footnote{``Einen
  wahren wissenschaftlichen Werth erkenne ich --- auf dem Felde der
  \emph{Mathematik} --- nur in concreten mathematischen Wahrheiten,
  oder sch\"{a}rfer ausgedr\"{u}ckt, 'nur in mathematischen
  Formeln'. Diese allein sind, wie die Geschichte der Mathematik
  zeigt, das Unverg\"{a}ngliche. Die verschiedenen Theorien f\"{u}r
  die Grundlagen der Mathematik (so die von Lagrange) sind von der
  Zeit weggeweht, aber die Lagrangesche Resolvente ist geblieben!''
  \textsf{Kronecker} 1884}
\begin{equation}
\label{Bar}
\mathfrak{B}_{p}(x,t) = \left\{\begin{array}{ll}
                              t^{-n/\lambda}\left[C-\frac{p-2}{p}\lambda^{1/(1-p)}\left(\frac{|x|}{t^{1/\lambda}}\right)^{\frac{p}{p-1}}\right]_{+}^{\frac{p-1}{p-2}} \:&\mbox{if $ t > 0,$}\\
0\ &\mbox{if $ t \leq 0,$}
\end{array}
\right.
\end{equation}
where $\lambda = n(p-2)+ p,$ is the leading example of a viscosity supersolution (=
p-supercaloric function). It plays the r\^{o}le of a fundamental solution, although the Principle of Superposition is naturally lost. It has a compact support in the
$x-$variable for each fixed instance $t$. Disturbances propagate with
finite speed and an interface (moving boundary) appears. Notice that
$$\int_{0}^{T} \! \! \int_{|x| < 1}|\nabla \mathfrak{B}_{p}(x,t)|^{p}\,dx \,dt  =  \infty$$
due to the singularity at the origin. Thus $ \mathfrak{B}_{p}$ fails to be a weak
supersolution in a domain containing the origin.\footnote{In my
  opinion, a definition of ``supersolutions'' that excludes the
  fundamental solution cannot be regarded as entirely
  satisfactory.}

A supersolution of a totally different kind is provided by the example
\begin{equation}
\label{separable1}
\mathfrak{M}(x,t) =
\begin{cases}
 t^{-\frac{1}{p-2}} \mathfrak{U}(x) \quad \text{if}\quad t > 0,\\
0,\qquad\text{if}\quad t \leq 0,
\end{cases}
\end{equation}
where the function $\mathfrak{U} > 0$ is the solution to an auxiliary elliptic equation. Here $\mathbf{p>2}$ (such a separable solution does not exist for $p = 2$).
 Our
 main theorem is:
\begin{thm}
\label{Mainpar}
Let $p > 2$ and suppose that $v = v(x,t)$ is a viscosity supersolution
in the domain $\Omega_{T}$ in $\Rn \times \mathbb{R}$. There are two disjoint cases:

\medskip

{\textsf Class} $\mathfrak{B}:\quad \mathbf{v \in L^{p-2}_{loc}(\Omega_T)}$. Then
$$v \in L^{q}_{loc}(\Omega_{T})\:\: \mbox{whenever}\:\: q < p-1+\frac{p}{n}$$
and the Sobolev derivative
$$\nabla v = \left(\frac{\partial v}{\partial x_{1}}, \cdots
  ,\frac{\partial v}{\partial x_{n}} \right)$$
exists and
$$\nabla v \in L^{q}_{loc}(\Omega_{T})\:\: \mathit{whenever}\:\: q <
p-1+\frac{1}{n+1}.$$
The summability exponents are sharp. Moreover,
$$\underset{\Omega_{T}}{\int\!\!\int} \left(-v\eta_{t} + \langle |\nabla v|^{p-2}\nabla
  v,\!\nabla \eta \rangle\,\right) dx\,dt \geq 0$$
for all $\eta \geq 0,\,\eta \in C_{0}^{\infty}(\Omega_{T})$.

\medskip

{\textsf Class} $\mathfrak{M}:\quad \mathbf{v \not \in L^{p-2}_{loc}(\Omega_T)}$. Then there exists a time $t_0,\,0<t_0<T,$ such that

$$\liminf_{\substack{(y,t)\to(x,t_0)\\t>t_0}} \bigl(v(y,t)(t-t_0)^{\frac{1}{p-2}}\bigr) > 0 \qquad\text{for all}\quad x\in \Omega.$$ In particular, $$\lim_{\substack{(y,t)\to(x,t_0)\\t>t_0}} v(y,t) \equiv +\infty\qquad\text{for all}\quad x\in \Omega.$$

\end{thm}

\bigskip

The occurrence of the void gap
$$[\,p-2,\,p-1+\tfrac{p}{n})$$
is arresting: either the function belongs to $L^{p-1+\frac{p}{n}-0}_{loc}$ or does not even belong to $L^{p-2}_{loc}$.
The two classes are deliberately labelled after their most representative members. The Barenblatt solution (\ref{Bar}) belongs to class $\mathfrak{B}$ and it shows that the exponents are sharp. Notice
that the time derivative is not included in the
statement. Actually, the time derivative need not be a function, as the
example $v(x,t) = 0$, when $t \leq 0$, and $v(x,t) = 1$, when $t >
0$ shows. Dirac's delta appears! It
is worth our while to emphasize that the gradient $\nabla v$ is not
present in the definitions of the viscosity supersolutions and the
p-supercaloric functions.

The separable solution (\ref{separable1}) belongs to class $\mathfrak{M}.$ Notice that at time $t=0$ it blows totally up, \textsf{its infinities fill the whole space} $\Rn.$ A convenient criterion to guarantee that a viscosity supersolution $v:\Omega \times (t_1,t_2)\to (-\infty,\infty]$ is \emph{not} of class  $\mathfrak{M}$ is that on the lateral boundary
$$\limsup_{(x,\tau) \to (\xi,t)}v(x,t) < \infty\qquad\text{for all} \quad (\xi,t)\in \partial \Omega \times (t_1,t_2).$$

An important feature is that the viscosity supersolutions are defined
at each point, not just almost everywhere in their domain. When it
comes to the pointwise behaviour,
one may  even exclude all future times so that only the instances $\tau < t$
are used for the calculation of $v(x,t)$, as in  the next theorem. (It
is also, of course, valid without restriction to the past times.)
\begin{thm}
\label{Ess}
Let $p \geq 2$. A viscosity supersolution of the Evolutionary
$p$-Laplace Equation satisfies
$$v(x,t) = \essinf_{\underset{\tau < t}{(y,\tau) \rightarrow (x,t)}} v(y,\tau)$$
at each interior point (x,t).
\end{thm}

In the calculation of \emph{essential limes inferior}, sets of
$(n+1)$-dimensional Lebesgue measure zero are neglected. We mention an
immediate consequence, which does not seem to be easily obtained by
other methods.
\begin{cor}
Two viscosity supersolutions that coincide almost everywhere are equal
at \emph{each} point.
\end{cor}

A general comment about the method employed in these notes is appropriate.
We do not know about proofs for viscosity supersolutions that would
totally stay within that framework. It must be emphasized that the
proofs are carried out for those supersolutions that are defined as
one does in Potential Theory, namely through comparison principles,
and then the results are valid even for the viscosity supersolutions,
just because, incidentally, they are the same functions. The
identification\footnote{A newer proof based on [JJ] of this
    fundamental identification is  included in these notes in Chapter
    7, replacing
  the reference [JLM].} of these two classes of
``supersolutions'' is not a
quite obvious fact. This limits the applicability of the method.

In passing, we also treat the measure data equation 
$$\frac{\partial v}{\partial t} - \nabla \cdot (|\nabla v|^{p-2}\nabla
v) = \mu$$
where the right-hand side is a Radon measure. It follows quite easily
from Theorem 2 that \textsf{each viscosity supersolution of class $\mathfrak{B}$ induces a
  measure} and
is a solution to the measure data equation. (The reversed problem,
which starts with a given measure $\mu$ instead of a given function
$v$, is a much investigated topic, cf. [BD].)

Some other equations that are susceptible of this kind of analysis are the
\textsf{ Porous Medium Equation}\footnote{The Porous Medium
  Equation
 is not well suited
for the viscosity theory (it is not ``proper''), although the comparison
principle works well. It is not $\nabla
v$ but $\nabla(|v|^{m-1}v)$ that is guaranteed to exist.}
$$\frac{\partial u}{\partial t} = \Delta (|u|^{m-1}u)$$ and
$$\frac{\partial (|u|^{p-2}u)}{\partial t} = 
\nabla\!\cdot \!\left(|\nabla u|^{p-2}\nabla u \right)  ,$$
but it does not seem to be known which equations of the form
$$\frac{\partial u}{\partial t} = F(x,t,u,\nabla u,D^{2}u)$$
enjoy the property of having their viscosity supersolutions in some
local Sobolev $x$-space. I hope that this could be a fruitful
research topic for the younger readers. ---I thank \textsc{T. Kuusi}
and \textsc{M. Parviainen} for a careful reading of the 
manuscript. The first version of these notes has appeared in [L3].

\section{The Stationary Equation}
\label{stationary}

For reasons of exposition\footnote{Chapter 3 is pretty independent of
  the present chapter.}, we begin with the stationary equation
\begin{equation}
\Delta_{p}u \equiv \divergence(|\nabla u|^{p-2}\nabla u) = 0,
\end{equation}
which offers some simplifications not present in the time dependent
situation. In principle, here we keep $p \geq 2$, although the theory
often allows that $1 <p  < 2$ at least with minor
changes. Moreover, the cases $p > n$, $p = n$, and $p < n$ often
require separate proofs. We sometimes skip the borderline case $p=n$.

 The fundamental solution
$$c|x|^{\frac{p-n}{p-1}}$$
does not belong to the Sobolev space $W^{1,p}_{loc}(\Rn).$ The problem is the origin. It is good to keep this in mind, when learning the definition below.

\begin{definition}
We say that $u \in W^{1,p}_{loc}(\Omega)$ is a \emph{weak solution}
in $\Omega$, if
$$\int_{\Omega} \langle |\nabla u|^{p-2} \nabla u,\nabla \eta \rangle
\,dx  = 0$$
for all $\eta \in C_{0}^{\infty}(\Omega)$. If, in addition, $u$ is
continuous, it is called  a $p$\emph{-harmonic function}. 

We say that $u \in W^{1,p}_{loc}(\Omega)$ is a \emph{weak supersolution}
in $\Omega$, if the integral is $\geq 0$ for all nonnegative $\eta \in C_{0}^{\infty}(\Omega)$. It is a \emph{weak subsolution} if the integral is $\leq 0.$
\end{definition}
The terminology suggests that ''super $\geq$ sub''. 

\begin{lemma} [Comparison Principle]
\label{ellcomp}
Let $u$ be a weak subsolution and $v$ a weak supersolution, and $u,v \in W^{1,p}(\Omega).$ If
$$\liminf_{x \to \xi}v(x) \geq \limsup_{x \to \xi}u(x)\quad\text{when}\quad \xi \in \partial \Omega$$
then $v \geq u$ almost everywhere in $\Omega.$
\end{lemma}
\emph{Proof:}
To see this, choose the test function $\eta = (u+\varepsilon -v)_+$ in the equations for $v$ and $u$ and subtract these:
\begin{equation*}
\int_{\{\varepsilon +u>v\}}\langle|\nabla v|^{p-2}\nabla v-|\nabla u|^{p-2}\nabla u,\nabla v-\nabla u\rangle \,dxdy\, \leq\, 0.
\end{equation*}
See Lemma 9. The integrand is strictly positive when $\nabla v \not = \nabla u,$
since $$\langle|b|^{p-2}b-|a|^{p-2}a,b-a\rangle \geq 2^{2-p}|b-a|^p,\quad p \geq2$$ holds for vectors. The result follows.\qquad $\Box$

\bigskip

 The weak solutions can, in
accordance with the elliptic regularity theory, be made continuous
after a redefinition in a set of Lebesgue measure zero. The H\"{o}lder
continuity estimate
\begin{equation}
\label{ellholder}
|u(x) -u(y)| \leq L\,|x - y|^{\alpha}
\end{equation}
holds when $x,y \in B(x_{0},r), \, B(x_{0},2r)\subset \subset \Omega$;
here $\alpha$ depends on $n$ and $p$ while $L$ also depends on the
norm $\|u\|_{p,B(x_{0},2r)}.$ We omit the proof. 
The continuous
weak solutions are called \emph{$p$-harmonic functions}\footnote{Thus the
2-harmonic functions are the familiar harmonic functions encountered
in Potential Theory.}.
In fact, even the gradient is continuous. One has $u \in
C_{loc}^{1,\alpha}(\Omega)$, where $\alpha = \alpha(n,p)$. This deep
result of N. Ural'tseva will not be needed here. According to [T1]
positive solutions obey the Harnack inequality. 
\begin{lemma}[Harnack's Inequality]
If the $p$-harmonic function $u$ is nonnegative in the ball
 $B_{2r} = B(x_{0},2r)$, then
$$\max_{\overline{B_{r}}}u \leq C_{n,p}\min_{\overline{B_{r}}}u\, .$$
\end{lemma}

\medskip
The $p$-Laplace equation is the \emph{Euler-Lagrange equation} of a
variational integral.
Let us recall the \emph{Dirichlet problem} in a bounded domain
$\Omega$. Let $f \in C(\overline{\Omega}) \cap W^{1,p}(\Omega)$
represent the boundary values. Then there exists a unique function $u$
in
$C(\Omega) \cap W^{1,p}(\Omega)$ such that $ u-f \in
W^{1,p}_{0}(\Omega)$ and
$$\int_{\Omega} |\nabla u|^{p}\,dx \leq \int_{\Omega}
|\nabla(u+\eta)|^{p}\,dx $$
for all $\eta \in C_{0}^{\infty}(\Omega)$. The minimizer is
$p$-harmonic. If the boundary $\partial \Omega$ is regular enough, the
boundary values are attained in the classical sense:
$$ \lim_{x \rightarrow \xi}u(x) = f(\xi),\:\:\: \xi \in \Omega.$$

When it comes to the super- and subsolutions, several definitions are
currently being used. We need the following ones:

\begin{description}
 \item{(1) {\small \textbf{weak supersolutions}}} (test functions under the integral sign)
 \item{(2) {\small \textbf{$p$-superharmonic functions}}} (defined via a comparison
   principle)
 \item{(3) {\small \textbf{viscosity supersolutions}}} (test functions evaluated at points
   of contact)
\end{description}

\textsf{The $p$-superharmonic functions and the viscosity supersolutions are
exactly the same functions}, see Chapter 7, [JJ], or [JLM]. They are not assumed to have any
derivatives. In contrast, the weak supersolutions are by their
definition required to belong to the Sobolev space
$W^{1,p}_{loc}(\Omega)$ and therefore their Caccioppoli estimates are at
our disposal. As we will see, locally \emph{bounded} $p$-superparabolic
functions (= viscosity supersolutions) are, indeed, weak
supersolutions, having Sobolev derivatives as they should. To this one
may add that the weak supersolutions are p-superharmonic functions,
 provided that the issue of semicontinuity be properly handled.
\begin{definition}\label{ellsup}
We say that a function $v: \Omega \rightarrow (-\infty,\infty]$ is 
\emph{p-superharmonic} in $\Omega$, if 
 \begin{description}
  \item{(i)} \ \  $v$ is finite in a dense subset
  \item{(ii)}\  $v$ is lower semicontinuous
  \item{(iii)} in each subdomain $D \subset \subset \Omega$\quad $v$ obeys the
    comparison principle:\\ if $h \in C(\overline{D})$ is $p$-harmonic in
    $D$, then the implication
    $$v|_{\partial D} \geq h|_{\partial D} \ \ \Rightarrow  \ \ v \geq h$$
    is valid. 
\end{description}
\end{definition}

\emph{Remarks}. For $p = 2$ this is the classical definition of
superharmonic functions due to  F. Riesz. It is
sufficient\footnote{This is not quite that simple in the parabolic
  case.} to assume
that $v \not \equiv \infty$ instead of (i). The fundamental solution
$|x|^{(p-n)/(p-1)}$, is not a weak supersolution in $\Rn$, merely
because it fails to belong to the right Sobolev space, but it is
 $p$-superharmonic.

\bigskip

 Some examples are the following functions.

\begin{eqnarray*}
(n-p)|x|^{- \frac{n-p}{p-1}} \;\;\; (n \not = p),\ \ \  
V(x) = \int\! \frac{\varrho(y)\,dy}{|x-y|^{n-2}} \;\;\; (p = 2, n \geq
3),\\
V(x) = \sum \frac{c_{j}}{|x-q_{j}|^{(n-p)/(p-1)}}\;\;\; (2 < p < n),\\
V(x) = \int\! \frac{\varrho(y)\,dy}{|x-y|^{(n-p)/(p-1)}} \;\;\; (2 < p <
n),\ \
v(x) = \min\{v_{1},v_{2}, \cdots ,v_{m}\}.
\end{eqnarray*}
The first example is \emph{the fundamental solution}, which fails to belong
to the ``natural'' Sobolev space $W_{loc}^{1,p}(\Rn)$.\footnote{Therefore it is
not a weak supersolution, but it is a viscosity supersolution and a
$p$-superharmonic function.}  The second is
the Newtonian potential. In the third example the $c_{j}$'s
are positive convergence factors and the  $q_{j}$'s are the rational
points; the superposition of fundamental solutions is credited to
M. Crandall and J. Zhang, cf. [CZ]. The last example says that one may take
the pointwise minimum of a finite number of $p$-superharmonic
functions, which is an essential ingredient in the celebrated Perron method, cf. [GLM]. 

The next definition is from the theory of viscosity solutions. One
defines them as being both viscosity super- and subsolutions, since it
is not practical to do it in one stroke.

\begin{definition}
Let $p \geq 2$. A function $v: \Omega \rightarrow (-\infty,\infty]$ is
called a \emph{viscosity supersolution}, if
 \begin{description}
  \item{(i)}\ \  $v$ is finite in a dense subset
  \item{(ii)}\  $v$ is lower semicontinuous
  \item{(iii)} whenever $x_{0} \in \Omega$ and $\phi \in
    C^{2}(\Omega)$ are such that $v(x_{0}) = \phi(x_{0})$ and $v(x) >
    \phi(x)$ when $x \not = x_{0}$, we have
$$\divergence\! \left(|\nabla \phi(x_{0})|^{p-2}\nabla \phi(x_{0}) \right)
\leq 0.$$
\end{description}
\end{definition}

\emph{Remarks}. The differential operator is evaluated only at the
point of contact. The singular case $1 < p < 2$ requires a
modification \footnote{There is no requirement when
  $\nabla \phi$ is 0, see [JLM].}, if it so happens that $\nabla
\phi(x_{0}) = 0$. Notice that each point has its own family of test
functions. If there is no test function touching from below at
$x_{0}$, then there is no requirement: the point passes for
free. Please, notice that nothing is said about the gradient $\nabla
v$, it is $\nabla \phi(x_{0})$ that appears.

\begin{thm}\label{uusi}
A $p$-superharmonic function is a viscosity supersolution.
\end{thm}
\emph{Proof:} Let $v$ be a $p$-superharmonic function in the domain $\Omega.$
In order to prove that $v$ is a viscosity supersolution, we use an indirect proof.
 Our \textsf{antithesis} is that there exist a
point $x_0 \in \Omega$ and a test function $\phi$
touching $v$ from below at $x_0$ and  satisfying the inequality
 $\Delta_p \phi(x_0) > 0$. Replacing
the test function with
$$\phi(x) - |x-x_0|^4$$
 we may further assume that the strict inequality $\phi(x)
< v(x)$ is valid when $x \not = x_0.$ The subtracted fourth power does not affect
$\Delta_p \phi (x_{0}).$ By
continuity we can assure that the strict inequality  $\Delta_p \phi(x) >
0$ holds in a small neighbourhood $U$ of the point $x_0$. 
  Now $\phi$ is
 $p$-\emph{sub}harmonic in
$U$ and by adding a small positive constant $m$, say $$ 2m = \max_{\,\partial
  U}\bigl(v(x) - \phi(x)\bigr)$$ we arrive at the following
situation. The function $\phi + m$ is $p$-subharmonic in $U$, which contains $x_0,$ and it is $\leq v$ on its boundary $\partial U.$
 By the comparison principle
 $\phi(x) + m  \leq v(x)$ in $U.$
This is a contradiction at the point $x = x_0.$ This proves the claim.\qquad $\Box$

\bigskip

The functions in the next lemma, the continuous weak supersolutions,
form a more tractable subclass, when it
comes to \emph{a priori} estimates, since
they are differentiable in Sobolev's sense.

\begin{lemma}\label{Crelle}
Let $v \in C(\Omega) \cap
 W^{1,p}(\Omega)$. Then the following
conditions are equi\-va\-lent:
\begin{description}
\item{(i)\ \ } $\int_{\Omega}|\nabla v|^{p} \, dx \leq
    \int_{\Omega}|\nabla (v+ \eta)|^{p} \, dx$ \quad when \quad $\eta \geq 0,\,
    \eta \in C_{0}^{\infty}(\Omega),$
\item{(ii) \ }$\int_{\Omega} \langle |\nabla v|^{p-2} \nabla v,\nabla
  \eta
 \rangle
\,dx  \geq 0$  \quad when \quad $\eta \geq 0,\, \eta \in C_{0}^{\infty}(\Omega),$
\item{(iii)}\  $v$ is $p$-superharmonic.
\end{description}
They imply that $v$ is a viscosity supersolution.
\end{lemma}

\emph{Proof}: The equivalence of (i) and (ii) is plain. So is the
necessity of (iii), stating that the comparison principle must hold. The crucial part is the sufficiency of (iii), which will
be established by the help of an \emph{obstacle problem} in the calculus of
variations. The function $v$ itself will act as an obstacle for the
admissible functions in the minimization of the $p$-energy
$\int_{D}|\nabla v|^{p} \, dx$ and it also induces the boundary values
in the subdomain $D$. If $D$ is a regular subdomain of $\Omega$, then
there exists a unique minimizer, say $w_{v}$, in the class
$$\mathcal{F}_{v} = \{w \in C(\overline{D}) \cap W^{1,p}(D)|\  w \geq v,\, w =
v\  \mathit{on}\  \partial D \}.$$
The crucial part is the continuity of $w_{v}$, cf. [MZ]. The solution of
the obstacle problem automatically  has
the property (i), and hence also (ii). We claim that $w_{v} = v$ in
$D$, from which the desired conclusion thus follows. The minimizer is a
$p$-harmonic function in the open set $\{w_{v} > v\}$ where the obstacle
does not hinder. On the boundary of this set $w_{v} = v$. Hence
the comparison principle, which $v$ is assumed to obey, can be
applied. It follows that $w_{v} \leq v$ in the same set. To avoid a
contradiction it must be the empty set. The conclusion is that  $w_{v}
= v$ in $D$, as desired. One can now deduce that (iii) is
sufficient. $\Box$

\medskip
A function, whether continuous or not, belonging to $W^{1,p}(\Omega)$
and satisfying (ii) in the previous lemma is called a \emph{weak
  supersolution}. 
For completeness we record below that weak supersolutions are
semicontinuous ``by nature''.

\begin{prop}
A weak supersolution $v \in W^{1,p}(\Omega)$ is lower semicontinuous
(after redefinition in a set of measure zero). We can define
$$ v(x) = \essinf_{y \rightarrow x} v(y)$$
pointwise. This representative is a $p$-superharmonic function.
\end{prop}

\emph{Proof:} The case $p > n$ is clear, since then the Sobolev space
contains only continuos functions (Morrey's inequality). In the range
$p < n$ we claim that 
$$ v(x) = \essinf_{y \rightarrow x} v(y)$$
at a.e. $x \in \Omega.$ The proof follows from this, because the
right-hand side is always lower semicontinuous. We omit two demanding
steps. First, it is required to establish that $v$ is locally
bounded from below. (See Theorem \ref{semico} below.) This is standard regularity theory. Second, for non-negative
functions  we use ''the weak
Harnack estimate''\footnote{Harnack's inequality can
  be replaced by the more elementary estimate
$$ \mathrm{ess}\sup_{B_{r}}\left(v(x_{0})-v(x)\right)_{+} \leq
  \frac{C}{|B_{2r}|}    \int_{B(x_{0},2r)}(v(x_{0})-v(x))_{+}\,dx$$
as a starting point for the proof. It follows immediately that also
$$ \mathrm{ess}\sup_{B_{r}}\left(v(x_{0})-v(x)\right) \leq
  \frac{C}{|B_{2r}|}    \int_{B(x_{0},2r)}|v(x_{0})-v(x)|\,dx.$$ If
  $x_{0}$ is a Lebesgue point, the integral approaches zero as
  $r\rightarrow0$ and it follows that
$$\essinf_{x \rightarrow x_{0}} v(x) \geq v(x_{0}).$$
The opposite inequality holds for ``arbitrary'' functions at their
Lebesgue points. (See the end of Chapter 4.)}
\begin{equation}
\label{avehar}
\left(\frac{1}{|B_{2r}|}\int_{B_{2r}}v^{q}\,dx \right)^{\frac{1}{q}} \leq
C\,\mathrm{ess}\inf_{B_{r}}v,
\end{equation}
when $q < n(p-1)/(n-p),\, C = C(n,p,q).$ This comes from the
celebrated Moser iteration, cf. [T1]. Taking\footnote{If $p\geq
  2n/(n+1)$ does not hold, we need a larger $q$.} $q=1$ and using
the non-negative function $v(x) - m(2r),$ where
$$ m(r) = \mathrm{ess}\inf_{B_{r}} v,$$
we have
\begin{gather*}
0 \leq \frac{1}{|B_{2r}|} \int_{B_{2r}}v\,dx\, -\, m(2r)\\ =
\frac{1}{|B_{2r}|}
  \int_{B_{2r}}(v(x)-m(2r))\,dx
\leq
C(m(r)-m(2r)).
\end{gather*}
Since $m(r)$ is monotone, $m(r)-m(2r) \rightarrow 0$ as $r \rightarrow
0$.
It follows that
$$\essinf_{y \rightarrow x_{0}} v(y) =  \lim_{r \rightarrow
  0}m(2r) =  \lim_{r \rightarrow
  0} \frac{1}{|B_{2r}|}    \int_{B(x_{0},2r)}v(x)\,dx$$
at \emph{each} point $x_{0}$. Lebesgue's differentiation theorem states that
the limit of the average on the right-hand side coincides with
$v(x_{0})$ at almost every point  $x_{0}$.

\bigskip

\begin{lemma}[Caccioppoli] Let $v \in C(\Omega) \cap W^{1,p}(\Omega)$ be a
$p$-superharmonic function. Then
$$\int_{\Omega} \zeta ^{p} |\nabla v|^{p}\,dx \leq p^{p}(\underset{\zeta
  \not = 0}{ \rm{osc}}\, v)^{p}\int_{\Omega}|\nabla \zeta|^{p}\,dx$$
holds for non-negative  $\zeta \in C_{0}^{\infty}(\Omega)$. If $v
\geq 0$, then
$$\int_{\Omega} \zeta ^{p}v^{-1- \alpha} |\nabla v|^{p}\,dx \leq
\left(\frac{p}{\alpha}\right)^{p}\int_{\Omega}v^{p-1- \alpha}|\nabla
  \zeta|^{p}\,dx$$
when $\alpha > 0$.
\end{lemma}

\emph{Proof}: To prove the first estimate, fix $\zeta$ and let $L =
\sup v$ taken over the set where $\zeta \not = 0$. Use the test
function
$$\eta = (L-v(x))\zeta(x)^{p}$$
in Lemma 9(ii) and arrange the terms.

To prove the second estimate, first replace $v(x)$ by $v(x) +
\varepsilon$, if needed, and use
$$\eta = v^{-\alpha}\zeta^{p}.$$
The rest is clear. $\Box$

\medskip
The special case $\alpha = p-1$ is appealing, since the right-hand
member of the inequality
\begin{equation}
\label{Gran}
\int_{\Omega} \zeta ^{p} |\nabla \log v|^{p}\,dx \leq 
\Bigl(\frac{p}{p-1} \Bigr)^{p}\int_{\Omega}|\nabla \zeta|^{p}\,dx
\end{equation}
is independent of the non-negative function $v$ itself.

We aim at approximating $v$ with functions for which Lemma \ref{Crelle}  is
valid. To this end, let $v$ be lower semicontinuous and bounded in
$\Omega$:
 $$0\leq v(x) \leq L.$$
Define
\begin{equation}
\label{ellinfimal}
v_{\varepsilon}(x) = \inf_{y \in \Omega} \left\{v(y) +
    \frac{|x-y|^{2}}{2\,\varepsilon}\right\}\:,\quad \varepsilon > 0.
\end{equation} 
Then 
\begin{itemize}
\item $v_{\varepsilon}(x) \nearrow v(x)$ as $\varepsilon \rightarrow
  0+$
\item $v_{\varepsilon}(x) - |x|^{2}/2\varepsilon$ is locally concave
  in $\Omega$
\item $v_{\varepsilon}$ is locally Lipschitz continuous in $\Omega$
\item The Sobolev gradient $\nabla v_{\varepsilon}$ exists and belongs
  to $L^{\infty}_{loc}(\Omega)$
\item The second Alexandrov  derivatives $D^2 v_{\varepsilon}$ exist. \footnotesize See Section 7.\normalsize
\end{itemize} 
The next to last assertion follows from Rademacher's theorem about Lipschitz
functions, cf. [EG]. Thus these ``infimal convolutions'' are rather regular. A
most interesting property for a bounded viscosity supersolution is
the following:

\begin{prop} \label{sitten} If $v$ is a bounded viscosity supersolution in $\Omega,$ the approximant  $v_{\varepsilon}$ is a viscosity
  supersolution in the open subset of $\Omega$ where
$${\rm dist}(x,\partial \Omega) > \sqrt{2L\varepsilon}.$$ Similarly, if $v$ is a $p$-superharmonic function, so is $v_{\varepsilon}.$ 
\end{prop}

\emph{Proof:} First, notice that for $x$ as required above, the infimum
is attained at some point $y = x^{\star}$ comprised in $\Omega$. The
possibility that $x^{\star}$ escapes to the boundary of $\Omega$ is prohibited
by the inequalities
\begin{gather*}
\frac{|x-x^{\star}|^{2}}{2 \varepsilon} \leq  \frac{|x-x^{\star}|^{2}}{2
  \varepsilon} + v(x^{\star}) = v_{\varepsilon}(x) \leq v(x) \leq L,\\
|x-x^{\star}| \leq \sqrt{2L\varepsilon} < {\rm dist}(x,\partial
\Omega).
\end{gather*}
This explains why the domain shrinks a little. Now we give two proofs.

\textsf{Viscosity proof}: Fix a point $x_{0}$  so that also $x^{\star}_{0}
\in \Omega$. Assume that the test function $\varphi$ touches
$v_{\varepsilon}$ from below at  $x_{0}$. Using 
\begin{eqnarray*}
\varphi(x_{0}) & =  v_{\varepsilon}(x_{0}) & =   \frac{|x_{0}-x^{\star}_{0}|^{2}}{2
  \varepsilon} + v(x^{\star}_{0}) \\
\varphi(x) & \leq   v_{\varepsilon}(x)  & \leq   \frac{|x-y|^{2}}{2
  \varepsilon} + v(y)
\end{eqnarray*}
we can verify that the function
$$\psi(x) = \varphi(x+x_{0}-x^{\star}_{0}) -\frac{|x_{0}-x^{\star}_{0}|^{2}}{2
  \varepsilon} $$
touches the original function $v$ from below at the point
$x^{\star}_{0}$. Since $x^{\star}_{0}$ is an interior point, the
inequality
$$\divergence \left(|\nabla \psi(x^{\star}_{0})|^{p-2}\nabla
  \psi(x^{\star}_{0})\right) \leq 0$$
holds by assumption. Because
$$\nabla\psi(x^{\star}_{0}) = \nabla\varphi(x_{0}),
\ \ D^{2}\psi(x^{\star}_{0}) = D^{2}\varphi(x_{0}),$$
we also have that
$$\divergence \left(|\nabla \varphi(x_{0})|^{p-2}\nabla
  \varphi(x_{0})\right) \leq 0$$
at the original point $x_{0}$, where $\varphi$ was touching
$v_{\varepsilon}$. Thus $v_{\varepsilon}$ fulfills the requirement in
the definition.

\textsf{Proof by Comparison}: Now we assume that $v$ is $p$-superharmonic in $\Omega$ and show that $v_{\varepsilon}$ is  $p$-superharmonic in  
$$\Omega_{\varepsilon} = \left\{x \in \Omega|\ \text{dist}(x,\partial \Omega) >
\sqrt{2L\varepsilon} \right\}.$$ We have to verify the comparison
principle for $v_{\varepsilon}$. To this end, let $D \subset \subset
\Omega_{\varepsilon}$
be a subdomain and suppose that $h \in C(\overline{D})$ is a
$p$-harmonic function so that $v_{\varepsilon}(x) \geq h(x)$ on the
boundary $\partial D$ or, in other words,
$$ \frac{|x-y|^{2}}{2 \varepsilon} + v(y) \geq h(x)\ \mathrm{when}\ x
\in \partial D,\  y \in \Omega.$$ 
Thus, writing $y = x + z$, we have
$$w(x) \equiv v(x+z) + \frac{|z|^{2}}{2 \varepsilon} \geq h(x),\ \ x \in
  \partial D$$
whenever $z$ is a small fixed vector. But also $w = w(x)$ is a
$p$-superharmonic function in $\Omega_{\varepsilon}$. By the
comparison principle $w(x) \geq h(x)$ in $D$. Given any point $x_{0}$
in $D$, we may choose $z = x_{0}^{\star} - x_{0}$. This yields
$v_{\varepsilon}(x_{0}) \geq h(x_{0})$. Since  $x_{0}$ was arbitrary,
we have verified that
$$ v_{\varepsilon}(x) \geq h(x),\ \ \mathrm{when}\  x \in D.$$
This concludes the proof.
\medskip
We record the following result.

\begin{cor}
\label{mant}
If $v$ is a bounded $p$-superharmonic function,
the approximant $v_{\varepsilon}$ is a weak supersolution in
$\Omega_{\varepsilon}$, i.e.
\begin{equation}
\int_{\Omega_{\varepsilon}} \langle |\nabla v_{\varepsilon}|^{p-2} \nabla
v_{\varepsilon}, \nabla \eta \rangle
\,dx  \geq 0\end{equation} when
 $\eta \geq 0$, $\eta \in C_{0}^{\infty}(\Omega_{\varepsilon}).$
\end{cor}

\emph{Proof}: This is a combination of the Proposition and Lemma \ref{Crelle}.\qquad  $\Box$

\medskip
The Caccioppoli estimate for $v_{\varepsilon}$ reads

$$\int_{\Omega} \zeta ^{p} |\nabla v_{\varepsilon}|^{p}\,dx \leq (p
L)^{p}\int_{\Omega}|\nabla \zeta|^{p}\,dx,$$
when $\varepsilon$ is so small that the support of $\zeta$ is in
$\Omega_{\varepsilon}$. By a compactness argument (a subsequence of)
$\nabla v_{\varepsilon}$ is locally weakly convergent in
$L^{p}(\Omega)$. We conclude that $\nabla v$ exists in Sobolev's sense
and that
$$\nabla v_{\varepsilon}\ \ \rightharpoonup\ \ \nabla v\ \
\mathrm{weakly\ in}\ L^{p}_{loc}(\Omega).$$
By the weak lower semicontinuity of the integral also
$$\int_{\Omega} \zeta ^{p} |\nabla v|^{p}\,dx \leq (p
L)^{p}\int_{\Omega}|\nabla \zeta|^{p}\,dx.$$
We have proved the first part of the next theorem.

\begin{thm}
\label{boundell}
Suppose that $v$ is a \emph{bounded} $p$-superharmonic function in
$\Omega$. Then the Sobolev gradient $\nabla v$ exists and $v \in
W^{1,p}_{loc}(\Omega)$. Moreover,
\begin{equation}
\int_{\Omega} \langle |\nabla v|^{p-2} \nabla v,\nabla
  \eta
 \rangle
\,dx  \geq 0\end{equation} 
for all $\eta \geq 0, \eta \in C_{0}^{\infty}(\Omega)$.
\end{thm}

\emph{Proof}: To conclude the proof, we show that the convergence
$\nabla v_{\varepsilon} \rightarrow \nabla v$ is strong in
$L^{p}_{loc}(\Omega)$, so that we may pass to the limit under the
integral sign in (8). To this end, fix a function $\theta \in
C^{\infty}_{0}(\Omega),\ 0 \leq \theta \leq 1$ and use the test
function $\eta = (v- v_{\varepsilon}) \theta$ in the equation for 
$v_{\varepsilon}$. Then
$$\int_{\Omega} \langle |\nabla v|^{p-2}\nabla v -|\nabla
v_{\varepsilon}|^{p-2}
\nabla v_{\varepsilon},\nabla((v-v_{\varepsilon})\theta) \rangle\,dx$$
$$\leq \int_{\Omega} \langle |\nabla v|^{p-2}\nabla v,
\nabla((v-v_{\varepsilon})\theta) \rangle\, dx \longrightarrow 0,$$
where the last integral approaches zero because of the weak
convergence. The first integral splits into the sum
\begin{eqnarray*}
\int_{\Omega}\theta \langle |\nabla v|^{p-2}\nabla v -|\nabla
v_{\varepsilon}|^{p-2}\nabla v_{\varepsilon},\nabla(v-v_{\varepsilon})
\rangle\,dx\\
+ \int_{\Omega}(v - v_{\varepsilon}) \langle |\nabla v|^{p-2}\nabla v -|\nabla
v_{\varepsilon}|^{p-2}\nabla v_{\varepsilon},\nabla \theta \rangle\, dx.
\end{eqnarray*}
 The last integral approaches zero because its absolute value is
 majorized by
$$\Bigl( \int_{D} (v-v_{\varepsilon})^{p}\,dx \Bigr)^{1/p}
\biggl[\Bigl(\int_{D}|\nabla v|^{p}\,dx \Bigr)^{(p-1)/p}+
  \Bigl(\int_{D}|\nabla v_{\varepsilon}|^{p}\,dx \Bigr)^{(p-1)/p}
\biggr] \| \nabla \theta \|,    $$
where $D$ contains the support of $\theta$ and $\|
v-v_{\varepsilon}\|_{p}$ approaches zero. Thus we have
established that $$ \int_{\Omega}\theta \langle |\nabla v|^{p-2}\nabla
v -|\nabla v_{\varepsilon}|^{p-2}\nabla v_{\varepsilon},\nabla
(v-v_{\varepsilon}) \rangle\,dx$$
approaches zero. Now the strong convergence of the gradients follows
from the vector inequality
\begin{equation}
\label{Vect}
2^{2-p}|b - a|^{p} \leq \langle |b|^{p-2}b - |a|^{p-2}a,b -a \rangle
\end{equation}
valid for $p > 2$. $ \Box$

\medskip
 It also follows that the Caccioppoli estimates in
Lemma \ref{caccioell} are valid for locally bounded $p$-superharmonic functions.
The case when $v$ is unbounded can be reached via the truncations
$$v_{k} = \min\{v(x),k\},\ \ \ k= 1,2,3,\dots,$$ because Theorem \ref{boundell}
holds for these locally bounded functions.
Aiming at a local result, we may  just by adding a constant 
assume that $v \geq 0$
in $\Omega$. The situation with $v = 0$ on the boundary $\partial
\Omega$ offers expedient simplifications. We shall describe an
iteration procedure,  under this extra assumption. See [KM].

\begin{lemma}
\label{Kilpmal}
Assume that $v \geq 0$ and that $v_{k} \in W_{0}^{1,p}(\Omega)$ when
$k = 1,2,\dots$ Then
$$\int_{\Omega}|\nabla v_{k}|^{p}\,dx \leq  k \int_{\Omega}|\nabla
v_{1}|^{p}\,dx $$
and, in the case $1 < p < n$
$$\int_{\Omega} v^{\alpha}\,dx \leq C_{\alpha} \left(1+\int_{\Omega}|\nabla
v_{1}|^{p}\,dx \right)^{\frac{n}{n-p}}$$
whenever $\alpha < \frac{n(p-1)}{n-p}$.
\end{lemma}

\emph{Proof}: Let $j$ be a large index and use the test functions
$$\eta_{k} = (v_{k} - v_{k-1})-(v_{k+1} - v_{k}),\ \ \
k=1,2,\cdots,j-1$$
in the equation for $v_{j}$, i.e.
$$\int_{\Omega} \langle|\nabla v_{j}|^{p-2}\nabla v_{j}, \nabla
\eta_{k}\rangle \,dx \geq 0 . $$
Indeed, $\eta_{k} \geq 0$. We obtain 
\begin{align*}
A_{k+1}& =
\int_{\Omega} \langle|\nabla v_{j}|^{p-2}\nabla v_{j}, \nabla
v_{k+1} - \nabla v_{k} \rangle \,dx \\
&\leq  \int_{\Omega} \langle|\nabla v_{j}|^{p-2}\nabla v_{j}, \nabla
v_{k} - \nabla v_{k-1} \rangle \,dx = A_{k}.
\end{align*}
Thus
$$A_{k+1} \leq A_{1} = \int_{\Omega}|\nabla v_{1}|^{p}\,dx$$
and hence
$$A_{1} +A_{2 } + \cdots + A_{j}\, \leq \,j A_{1}.$$
The ``telescoping'' sum becomes
$$\int_{\Omega}|\nabla v_{j}|^{p}\,dx
 \leq j \int_{\Omega}|\nabla v_{1}|^{p}\,dx.$$
This was the first claim.

If $1 < p < n$, it follows from Tshebyshev's and Sobolev's
 inequalities that
\begin{equation*}
j|j \leq v \leq 2j|^{\frac{1}{p*}}  \leq
  \left(\int_{\Omega}v_{2j}^{p*}\,dx \right)^{\frac{1}{p*}} 
\leq S\left(\int_{\Omega} |\nabla v_{2j}|^{p}\,dx \right)^{\frac{1}{p}}
\leq S(2j)^{\frac{1}{p}} A_{1}^{\frac{1}{p}},
\end{equation*}
where $p* = np/(n-p)$. We arrive at the estimate

$$|j \leq v \leq 2j| \leq C j^{- \frac{n(p-1)}{n-p}}
A_{1}^{\frac{n}{n-p}}$$
for the measure of the level sets. To conclude the proof we write

$$ \int_{\Omega}v^{\alpha}\, dx  =  \int_{v \leq 1} v^{\alpha} \, dx + 
 \sum_{j=1}^{\infty} \int_{2^{j-1} < v \leq 2^{j}} v^{\alpha} \, dx.$$
Since

$$ \int_{2^{j-1} < v \leq 2^{j}} v^{\alpha} \, dx \leq C\, 2^{j \alpha}
2^{-(j-1)\frac{n(p-1)}{n-p}} A_{1}^{\frac{n}{n-p}},$$
the series converges when $\alpha$ is as prescribed. $\Box$

\medskip
It remains to abandon the restriction about zero boundary values and
to estimate
$$\int_{\Omega}|\nabla v_{1}|^{p} \, dx.$$
The \textsf{reduction to zero boundary values} is done locally in a
ball $B_{2r} \subset \subset \Omega$. Suppose first that $v \in
C(\overline{B_{2r}}) \cap W^{1,p}(B_{2r}),\  v \geq 0$, and define
$$ w = \left\{ \begin{array}{ll}
v \ \mbox{in $\overline{B_{r}}$}\\
h \ \mbox{in $B_{2r} \backslash \overline{B_{r}}$}
\end{array}
\right.$$
where $h$ is the $p$-harmonic function in the annulus having outer
boundary values zero and inner boundary values $v$. Now $h \leq
v$. The so defined $w$ is $p$-superharmonic in $B_{2r}$, which follows
by comparison. It is quite essential that the original $v$ was defined
in a domain larger than $B_{r}$! We also have 
$$\int_{B_{2r}} |\nabla w|^{p}\, dx \leq C r^{n-p}
(\max_{B_{2r}}w)^{p}$$
after some estimation.\footnote{It is important to include the whole
  $B_{2r}$. Of course, the Caccioppoli estimate (Lemma 11) will do
  over any smaller ball $B_{\varrho}$, $\varrho < 2r$. To get the
  missing estimate, say over the boundary annulus $B_{2r} \setminus
  B_{3r/2}$, the test function $\eta = \zeta h$ works in Definition 5,
  where $\zeta = 1$ in the annulus and $= 0$ on $\partial B_{r}$. The
  zero boundary values of the weak solution $h$ were essential.} 

Finally, if $v \in W^{1,p}(B_{2r})$ is semicontinuous and bounded
 (but not necessarily continuous), then we
first modify the approximants $v_{\varepsilon}$ defined as in (7) and
obtain $p$-superharmonic functions $w_{\varepsilon}$. Since
$0 \leq w_{\varepsilon} \leq v_{\varepsilon} \leq v$, the previous estimate
becomes
$$\int_{B_{2r}} |\nabla w_{\varepsilon}|^{p}\, dx \leq C r^{n-p}
(\max_{B_{2r}}v)^{p}$$
and, by the weak lower semicontinuity of the integral, we can pass to
the limit as $\varepsilon$ approaches zero. We end up with a
$p$-superharmonic function $w \in W^{1,p}_{0}(B_{2r})$ such that $w = v$
in $B_{r}$ and, in particular,
$$\int_{B_{r}}|\nabla v|^{p} \, dx \leq  \int_{B_{2r}}|\nabla w|^{p} \,
dx \leq C r^{n-p} (\max_{B_{2r}} v)^{p}.$$
This is the desired modified function. Now, repeat the procedure 
with every function $\min\{v(x),k\}$ in sight.
 We obtain
$$\int_{B_{r}}|\nabla v_{1}|^{p} \,dx \leq \int_{B_{2r}}|\nabla w_{1}|^{p} \,dx \leq C r^{n-p} 1^{p}$$
for the modification of $v_{1} = \min\{v(x),1\}$. We have achieved that the
bounds in the previous lemma hold for the modified function over the
domain $B_{2r}$ and \emph{a fortiori} for the original $v$, estimated
only over the smaller ball $B_{r}$. Such a \emph{local} estimate is
all that is needed in the proof of the theorem below.
\begin{thm}
\label{Mainell}
Suppose that $v$ is a  $p$-superharmonic
function in $\Omega$. Then 
$$v \in L_{loc}^{q}(\Omega), \text{\, whenever \,} q <
\frac{n(p-1)}{n-p}$$
in the case $1 < p \leq n$ and $v$ is continuous if $p > n$. Moreover,
$\nabla v$ exists in Sobolev's sense\footnote{Strictly speaking, one needs $p >
  2-\frac{1}{n}$ so that $q \geq 1$. This can be circumvented.} and
$$\nabla v \in L_{loc}^{q}(\Omega),\ \mbox{\, whenever \,}\ q <
\frac{n(p-1)}{n-1}$$ in the case $1 < p \leq n$. In the case $p > n$
we have $\nabla v \in L_{loc}^{p}(\Omega)$. Finally,
\begin{equation}
\int_{\Omega}\langle|\nabla v|^{p-2}\nabla v,\nabla \eta \rangle \,
dx \geq 0
\end{equation}
when $\eta \geq 0,\,  \eta \in C_{0}^{\infty}(\Omega)$.
\end{thm}

\emph{Proof:} In view of the local nature of the theorem we may assume
that $v > 0$. According to the previous construction we can further
reduce the proof to the case $v_{k} \in W_{0}^{1,p}(B_{2r})$ for each
truncation at height $k$. The first part of the theorem is included in
Lemma \ref{Kilpmal} when $1 < p < n$. We skip the borderline case $p = n$. The case
$p > n$ is related to the fact that then all functions in the Sobolev space
$W^{1,p}$ are continuos.

 We proceed to the \textsf{estimation of the
  gradient}. First we keep $1 < p < n$ and write
\begin{gather*}
 \int_{B}|\nabla v_{k}|^{q} \, dx = \int_{B}
v_{k}^{\frac{(1+\alpha)q}{p}} \biggl|\frac{\nabla
    v_{k}}{v_{k}^{(1+\alpha)/p}} \biggr|^{q} \, dx\\
 \leq \left\{\int_{B}
v_{k}^\frac{(1+\alpha)q}{p-q} \, dx \right\}^{1- \frac{q}{p}}
\left\{\int_{B}v_{k}^{-1-\alpha}|\nabla v_{k}|^{p} \, dx
  \right\}^{\frac{q}{p}}.
\end{gather*}
Take $q < n(p-1)/(n-1)$ and fix $\alpha$ so that 
$$\frac{(1+\alpha)q}{p-q} < \frac{n(p-1)}{n-p}.$$
Continuing, the Caccioppoli estimate yields the majorant
\begin{equation}
\leq \left\{\int_{B}
v_{k}^\frac{(1+\alpha)q}{p-q} \, dx \right\}^{1- \frac{q}{p}}C
\left\{\int_{2B}v_{k}^{p-1-\alpha} \, dx \right\}^{\frac{q}{p}}.
\end{equation}
We can take $v \geq 1$. Then
let $k \longrightarrow \infty$. Clearly, the resulting majorant is
 finite (Lemma \ref{Kilpmal}). This concludes the case $1 < p <
n$. 

If $p > n$ we obtain that
$$ \int_{B_{r}}|\nabla \log v_{k}|^{p} \, dx \leq C r^{n-p}$$
from (\ref{Gran}), where $C$ is independent of $k$. Hence $\log v_{k}$ is
continuous. 
So is $v$ itself. Now
$$\int_{B_{r}}|\nabla v_{k}|^{p} \, dx
 = \int_{B_{r}}v_{k}^{p}\,|\nabla \log v_{k}|^{p} \, dx \leq C \|
v\|_{\infty}^{p} r^{n-p}$$
implies the desired $p$-summability of the gradient. \quad $\Box$

\medskip

It stands to reason that the lower semicontinuous solutions of (\ref{boundell})
are $p$-superharmonic functions. However, this is not known under the
summability assumption $\nabla v \in L_{loc}^{q}(\Omega)$ accompanying
the differential equation, if $q < p$ and $p < n$.
 In fact, an example of J. Serrin
indicates that even for solutions to linear equations strange
phenomena occur, cf. [S]. False solutions appear, when the \emph{a
priori} summability of the gradient is too poor. About this topic
there  is nowadays a 
theory credited to T. Iwaniec, cf. [L].\footnote{ The ``pathological
 solutions''
of  Serrin are now called``very
weak solutions''.}  

\section{The Evolutionary Equation}

This chapter is rather independent of the previous one. After some
definitions  we first treat bounded supersolutions and then
the unbounded ones. As a mnemonic rule, $v_{t} \geq \Delta_{p} v$ for
smooth supersolutions, $u_{t} \leq \Delta_{p} u$ for smooth
subsolutions. We need the following classes of supersolutions:

\begin{description}
 \item{(1) {\small \textbf{weak supersolutions}}} (test functions under the integral sign)
 \item{(2) {\small \textbf{$p$-supercaloric functions}}} (defined via a comparison
   principle)
 \item{(3) {\small \textbf{viscosity supersolutions}}} (test functions evaluated at points
   of contact)
\end{description}

The weak supersolutions do not form a good closed class under monotone
 convergence.

\subsection{Definitions}

We first define the concept of solutions, p-supercaloric functions and
viscosity supersolutions. 
 The  section ends with an outline of the procedure for the proof.

Suppose that $\Omega$ is a bounded domain in $\mathbb{R^{n}}$ and
consider the space-time cylinder $\Omega_{T} = \Omega \times
(0,T)$. Its \emph{parabolic boundary} consists of the portions $ \Omega
\times \{0\}$ and $\partial \Omega \times [0,T]$.

\begin{definition}
 In the case\footnote{The singular case $1 < p <
2$ requires an extra \emph{a priori} assumption, for example, $u \in
L^{\infty}(0,T;L^{2}(\Omega))$ will do.} $p \geq 2$ we say that $u \in
L^{p}(0,T;W^{1,p}(\Omega))$ is a \emph{weak solution} of the
Evolutionary p-Laplace Equation, if 
\begin{equation}
\integ{(-u\phi_{t} + \langle |\nabla u|^{p-2}\nabla u,\!\nabla \phi
  \rangle)\,}
 = 0
\end{equation}
for all $ \phi \in C^{1}_{0}(\Omega_{T})$. If the integral  is
$\geq 0$ for all test functions $\phi \geq 0$, we say that $u$ is a
\emph{weak supersolution}.
\end{definition}

 In particular, one has the requirement
$$ \integ{(|u|^{p} + |\nabla u|^{p})} <  \infty.$$
Sometimes it is enough to require that $u\in L^p_{loc}(0,T;W^{1,p}_{loc}(\Omega)).$
By the regularity theory one may regard a weak solution $u = u(x,t)$ as
 continuous.\footnote{The weak supersolutions are lower semicontinuous
   according to [K], see Chapter 4.} For
simplicity, we call the continuous weak solutions for \emph{$p$-caloric
  functions}\footnote{One may argue that this is more adequate than
 ``$p$-parabolic functions'', which is in use.}.

The \emph{interior H\"{o}lder estimate}\footnote{This is weaker than
  the estimate in [dB]. See also
  [U] for intrinsic scaling.} takes the following form for solutions 
according to [dB]. In the subdomain $D \times(\delta,T- \delta)$ 
\begin{equation}
\label{Holder}
|u(x_{1},t_{1}) - u(x_{2},t_{2})| \leq \gamma
\|u\|_{L^{\infty}(\Omega_{T})} \left(|x_{1}-x_{2}|^{\alpha}+
  |t_{1}-t_{2}|^{\alpha/p} \right), 
\end{equation}
where the positive exponent $\alpha$ depends only on $n$ and  $p$,
while the constant $\gamma$ depends, in addition, on the distance to
the subdomain. Also an \emph{intrinsic Harnack inequality} is valid, see Lemma \ref{Harnackpar} below.

Recall that the \emph{parabolic boundary} of the domain $\Omega_T = \Omega \times (0,T)$ is
$$\Omega \times \{0\}\,\,\cup\,\, \partial \Omega \times [0,T].$$
The part $\Omega \times \{T\}$ is excluded.

\begin{prop}[Comparison Principle]
\label{Comparison}
 Suppose that $v$ is a weak
  supersolution and $u$ a weak subsolution, $u,v \in
  L^{p}(0,T;W^{1,p}(\Omega))$, satisfying 
$$ \liminf v \geq \limsup u$$
on the parabolic boundary. Then $v \geq u$ almost everywhere in the domain
$\Omega_{T}.$
\end{prop}

\emph{Proof:} This is well-known and we only give a formal proof. For
a non-negative test function $\varphi \in C^{\infty}_{0}(\Omega_{T})$
the equations
\begin{gather*}
\int_{0}^{T}\!\int_{\Omega}(-v\varphi_{t}+\langle|\nabla
v|^{p-2}\nabla v,\nabla \varphi\rangle) \,dx\,dt \geq 0\\
\int_{0}^{T}\!\int_{\Omega}(+u\varphi_{t}-\langle|\nabla
u|^{p-2}\nabla u,\nabla \varphi\rangle) \,dx\,dt \geq 0
\end{gather*}
can be added. Thus
\begin{equation*}
\int_{0}^{T}\!\int_{\Omega}\left((u-v)\varphi_{t}+\langle  |\nabla
v|^{p-2}\nabla v - |\nabla
u|^{p-2}\nabla u,\nabla \varphi\rangle \right) \,dx\,dt \geq 0.
\end{equation*}
These equations remain true if $v$ is replaced by $v+\varepsilon$,
where $\varepsilon$ is any positive constant. To complete the proof we choose
(formally) the test function to be
$$\varphi = (u-v-\varepsilon)_{+}\eta,$$
where $\eta = \eta(t)$ is a cut-off function; the plain choice $\eta(t) = T-t$
will do here. We arrive at
\begin{align*}
&\phantom{=a }\int_{0}^{T}\!\int_{\{u\geq v+\varepsilon\}}\eta(\langle  |\nabla
v|^{p-2}\nabla v - |\nabla
u|^{p-2}\nabla u,\nabla v - \nabla u \rangle) \,dx\,dt \\
&\leq \int_{0}^{T}\!\int_{\Omega}(u-v-\varepsilon)^{2}_{+}\eta'\,dx\,dt
+\frac{1}{2}\int_{0}^{T}\!\int_{\Omega}\eta \frac{\partial}{\partial
  t}(u-v-\varepsilon)^{2}_{+}\,dx\,dt\\
&= \frac{1}{2}\int_{0}^{T}\!\int_{\Omega}(u-v-\varepsilon)^{2}_{+}\eta'\,dx\,dt\\
&= -
\frac{1}{2}\int_{0}^{T}\!\int_{\Omega}(u-v-\varepsilon)^{2}_{+}\,dx\,dt\,
\leq \,0.
\end{align*}
Since the first integral is non-negative by the structural inequality (10), the last integral is, in fact, zero. Hence the integrand $(u-v-\varepsilon)_+ = 0$ almost everywhere. But this means that
$$u\leq v+\varepsilon $$
almost everywhere. Since $\varepsilon >0$ was arbitrary, we have the desired inequality $v \geq u$ \, a.e.. \qquad $\Box$

\bigskip

For continuous functions the Comparison Principle is especially appealing. Then the conclusion is valid at every point. As we will see later, the redefined functions
$$v_* = \essinf v,\qquad u^* = \mathrm{ess\,limsup}\, u$$
are weak super- and subsolutions, and $v_*\geq u^*$ at each point.

\begin{definition}
We say that a function $v: \Omega_{T} \rightarrow (-\infty,\infty]$ is 
\emph{$p$-supercaloric} in $\Omega_{T}$, if 
 \begin{description}
  \item{(i)} \ \  $v$ is finite in a dense subset
  \item{(ii)}\  $v$ is lower semicontinuous
  \item{(iii)} in each cylindrical subdomain $D \times (t_{1},t_{2})
 \subset \subset \Omega_{T}$ $v$ obeys the
    comparison principle:\\ if $h \in C(\overline{D} \times
    [t_{1},t_{2}])$
 is $p$-caloric in
    $D \times (t_{1},t_{2})$, then $v \geq h$ on the parabolic
    boundary of $D \times (t_{1},t_{2})$ implies that $v \geq h$ in
    the whole subdomain.
\end{description}
\end{definition}

 As a matter of fact, every weak supersolution has a semicontinuous representative which is a $p$-supercaloric function. (This is postponed till Section \ref{Semicontinuous})
 The leading example
is the Barenblatt solution, which is a $p$-supercaloric function in the
whole $\mathbb{R^{n+1}}$. Another example is any function of the
form
$$v(x,t) = g(t),$$
where $g(t)$ is an arbitrary monotone increasing lower semicontinuous
function. We also mention
$$v(x,t) + \frac{\varepsilon}{T-t}, \ \ \ \ 0 < t < T,$$
$$v(x,t) = \min \{v_{1}(x,t),\ldots,v_{j}(x,t)\}.$$ The pointwise minimum of (finitely many) $p$-supercaloric functions is employed in Perron's Method. Finally, if $v \geq
0$ is a $p$-supercaloric function, so is the function obtained by
redefining $v(x,t) = 0$ when $t \leq 0$.

\paragraph{A Separable Minorant.} Separation of variables suggests that there are $p$-caloric functions of the type
$$v(x,t) = (t-t_0)^{-\frac{1}{p-2}}u(x).$$
Indeed, if  $\Omega$ is a  domain of finite measure, there exists a $p$-caloric function of the form
\begin{equation}
\label{sepa}
\mathfrak{V}(x,t) =
\frac{\mathfrak{U}(x)}{(t-t_0)^{\frac{1}{p-2}}},\quad \text{when}\quad t > t_0
\end{equation}
where $\mathfrak{U} \in C(\Omega)\cap W^{1,p}_{0}(\Omega)$ is a weak solution to the elliptic equation 
\begin{equation}
\label{apuyh}
 \nabla \!\cdot\!\bigl(|\nabla \mathfrak{U}|^{p-2}\nabla \mathfrak{U}\bigr)\:+\;\tfrac{1}{p-2}\,\mathfrak{U}\:=\;0
\end{equation}
and $\mathfrak{U}>0$ in $\Omega.$  The solution $\mathfrak{U}$ is unique\footnote{Unfortunately, the otherwise reliable paper [\textsc{J. Garci'a Azorero, I. Peral Alonso}: {\it Existence and nonuniqueness for the $p$-Laplacian: Nonlinear eigenvalues}, Communications in Partial Differential Equations \textbf{12}, 1987, pp. 1389--1430], contains a misprint exactly for those parameter values that would yield this function.}. (Actually, $\mathfrak{U} \in C^{1,\alpha}_{loc}(\Omega)$ for some exponent $\alpha = \alpha(n,p) > 0.$) The extended function
\begin{equation}
\label{separable2}
\mathfrak{V}(x,t) = 
\begin{cases}
\dfrac{\mathfrak{U}(x)}{(t-t_0)^{\frac{1}{p-2}}},\quad \text{when}\quad t > t_0\\\quad
0\qquad \text{when} \quad t \leq t_0.
\end{cases}
\end{equation}
is $p$-supercaloric in $\Omega \times \mathbb{R}.$ The existence of $\mathfrak{U}$ follows by the direct method in the Calculus of Variations, when the quotient
$$J(w) = \dfrac{\int_{\Omega}\!|\nabla w|^p\,dx}{\Bigl(\int_{\Omega}\!w^2\,dx\Bigr)^{\frac{p}{2}}}$$
is minimized among all functions $w$ in $W^{1,p}_0(\Omega),\,w\not \equiv 0.$ Replacing $w$ by its absolute value $|w|,$
we may assume that all functions are non-negative. Notice that $J(\lambda w) = J(w)$ for $\lambda = $ constant. Sobolev's and H\"{o}lder's inequalities
yield
$$J(w) \geq c(p,n)|\Omega|^{1-\frac{p}{n}-\frac{p}{2}},\qquad c(p,n) > 0$$
and so $$J_0 = \inf_{w}J(w) > 0.$$ Choose a minimizing sequence of admissible normalized functions
$w_{j}:$ 
$$\lim_{j \to \infty}J(w_{j}) = J_0,\qquad \|w_{j}\|_{L^p(\Omega)} = 1.$$
By compactness, we may extract a subsequence such that
\begin{align*}
\nabla w_{j_k} &\rightharpoonup \nabla w\qquad \text{weakly in} \quad L^p({\Omega})\\
w_{j_k}&\longrightarrow w \qquad \text{strongly in} \quad L^p({\Omega})
\end{align*}
for some function $w$.
The weak lower semicontinuity of the integral implies that
$$J(w) \leq \liminf_{k\to \infty}J(w_{j_k}) = J_0.$$
Since $w\in W^{1,p}_0(\Omega)$ this means that $w$ is a minimizer. We have $w\geq0,\,w\not\equiv 0.$

It follows that $w$ has to be a weak solution of the Euler--Lagrange Equation
$$ \nabla \!\cdot\!\bigl(|\nabla w|^{p-2}\nabla w\bigr) \,+ J_0 \|w\|_{L^p(\Omega)}^{p-2}w\:=\;0$$
where $\|w\|_{L^p(\Omega)} = 1.$ By elliptic regularity theory $w\in C(\Omega),$ see [T1] and [G]. Finally, since\, $ \nabla \!\cdot\!\bigl(|\nabla w|^{p-2}\nabla w\bigr) \leq 0$\, in the weak sense and  $w \geq 0$
we have that $w > 0$ by the Harnack inequality (\ref{avehar}). A normalization remains to be done. The function
$$\mathfrak{U} = Cu,\quad \text{where} \quad J_0C^{p-2} = \tfrac{1}{p-2},$$
will do.

\medskip

The next definition is from the theory of viscosity solutions. One
defines them as being both viscosity super- and subsolutions, since it
is not practical to do it in one stroke.

\begin{definition}
Let $p \geq 2$. A function $v: \Omega_{T} \rightarrow (-\infty,\infty]$ is
called a \emph{viscosity supersolution}, if
 \begin{description}
  \item{(i)}\ \  $v$ is finite in a dense subset
  \item{(ii)}\  $v$ is lower semicontinuous
  \item{(iii)} whenever $(x_{0},t_{0}) \in \Omega_{T}$ and $\phi \in
    C^{2}(\Omega_{T})$ are such that $v(x_{0},t_{0}) =
    \phi(x_{0},t_{0})$ and
$v(x,t) >
    \phi(x,t)$ when $(x,t) \not = (x_{0},t_{0})$, we have
$$\phi_{t}(x_{0},t_{0}) \geq \nabla \cdot \! \left(|\nabla \phi(x_{0},t_{0})|^{p-2}\nabla
  \phi(x_{0},t_{0})
 \right).$$
\end{description}
\end{definition}

\textsf{The p-supercaloric functions are exactly the same as the
  viscosity supersolutions}. For a proof of this fundamental
equivalence we refer to Section 7.2 or   [JLM]. However, the following implication is easyly obtained.

\begin{prop}
\label{suffviscpar} A viscosity supersolution of the Evolutionary $p$-Laplace Equation is a $p$-supercaloric function.
\end{prop}
\emph{Proof:} Similar to the elliptic case in Theorem \ref{uusi}.

\bigskip

Continuous weak supersolutions are $p$-supercaloric functions according to the Comparison Principle (Proposition \ref{Comparison}).

 We aim at proving the
summability results (Theorem \ref{Mainpar}) for a general $p$-supercaloric function $v$. An {\sf outline of the procedure} is the following.
\begin{itemize} {\small
\item{\it Step 1.} Assume first that $v$ is bounded.
\item{\it Step 2.} Approximate $v$ locally with infimal convolutions
  $v_{\varepsilon}$. These are differentiable.
\item{\it Step 3.} The infimal convolutions are $p$-supercaloric functions
  and they are shown to be weak supersolutions of the equation
 (with test functions  under the integral sign).
\item{\it Step 4.} Estimates of the Caccioppoli type  for $v_{\varepsilon}$
 are extracted from the
  equation.
\item{\it Step 5.} The Caccioppoli estimates are passed over from
  $v_{\varepsilon}$ to $v$. This concludes the proof for bounded
  functions.
\item{\it Step 6.} The unbounded case is reached via the bounded $p$-supercaloric
  functions $v_{k} = \min\{v,k\}, k = 1,2,\cdots$, for which the
  results in Step 5 already are available.
\item{\it Step 7.} An iteration with respect to the index $k$ is designed
  so that the final result does not blow up as $k \rightarrow
  \infty$. This works well when the parabolic boundary values (in the
  subdomain studied) are zero.
\item{\it Step 8.} An extra construction is performed to reduce the proof
  to the situation of zero parabolic boundary values (so that the
  iterated result in Step 7 is at our disposal).} \textsf{This is not possible for class} $\mathfrak{M},$ which is singled out.
\end{itemize}

\subsection{Bounded Supersolutions}\label{Bounded Supersolutions}

We aim at proving Theorem \ref{Mainpar}, which was given in the Introduction. The first step is to consider \emph{bounded}
$p$-supercaloric functions. We want to prove that they are weak
supersolutions.  First we approximate them with their infimal
convolutions. Then estimates mainly of the Caccioppoli type are proved
for these approximants. Finally, the so obtained estimates are passed
over to the original functions. Assume therefore that 
$$0 \leq v(x,t) \leq L,\ \ \  (x,t) \in \Omega_{T} = \Omega \times (0,T).$$
The approximants 
$$v_{\varepsilon}(x,t) = \inf_{(y,\tau) \in \Omega_{T}} \left\{v(y,\tau) +
    \frac{|x-y|^{2}+|t-\tau|^{2}}{2\,\varepsilon}\right\}\:,\quad
  \varepsilon > 0,
$$ 
have the properties 
\begin{itemize}
\item $v_{\varepsilon}(x,t) \nearrow v(x,t)$ as $\varepsilon \rightarrow
  0+$
\item $v_{\varepsilon}(x,t) - \frac{|x|^{2}+t^{2}}{2\varepsilon}$ is
  locally
 concave
  in $\Omega_{T}$
\item $v_{\varepsilon}$ is locally Lipschitz continuous in $\Omega_{T}$
\item The Sobolev derivatives $\frac{\partial
    v_{\varepsilon}}{\partial t}$ and $\nabla v_{\varepsilon}$ exist and belong
  to $L^{\infty}_{loc}(\Omega_{T})$
\item The second Alexandrov derivatives of  $v_{\varepsilon}$
  exist\footnote{See Section 7.}
\end{itemize} 
The next to last assertion follows from Rademacher's theorem about Lip\-schitz
functions. Thus these ``infimal convolutions'' are differentiable almost everywhere. The
existence of the time derivative is very useful. A
most interesting property for a bounded viscosity supersolution is
the following:
\begin{prop}\label{infcon} Suppose that $v$ is a viscosity supersolution in $\Omega_T.$ The approximant  $v_{\varepsilon}$ is a viscosity
  supersolution in the open subset of $\Omega_{T}$ where
$$\rm{dist}((x,t),\partial \Omega_{T}) > \sqrt{2L\varepsilon}.$$
Similarly, if $v$ is $p$-supercaloric, so is $v_{\varepsilon}.$
\end{prop}
\emph{Proof}: First, notice that for $(x,t)$ as required above, the infimum
is attained at some point $(y,\tau) = (x^{\star},t^{\star})$ comprised
 in $\Omega_{T}$. The
possibility that $(x^{\star},t^{\star})$ escapes to the boundary
 of $\Omega$ is prohibited
by the inequalities
\begin{gather*}
\frac{|x-x^{\star}|^{2} + |t-t^{\star}|^{2}}{2 \varepsilon}
 \leq  \frac{|x-x^{\star}|^{2} + |t-t^{\star}|^{2}}{2
  \varepsilon} + v(x^{\star},t^{\star}) \\  = v_{\varepsilon}(x,t) \leq
v(x,t)
 \leq L,\\
\sqrt{|x-x^{\star}|^{2} + |t-t^{\star}|^{2}} \leq \sqrt{2L\varepsilon}
 < \rm{dist}((x,t),\partial
\Omega_{T}).
\end{gather*}
Thus the domain shrinks a little. Again there are two proofs.

\textsf{Viscosity proof}: Fix a point $(x_{0},t_{0})$  so that also
 $(x^{\star}_{0},t^{\star}_{0})
\in \Omega_{T}$. Assume that the test function $\varphi$ touches
$v_{\varepsilon}$ from below at  $(x^{\star}_{0},t^{\star}_{0})$. Using 
\begin{eqnarray*}
\varphi(x_{0},t_{0}) & =  v_{\varepsilon}(x_{0},t_{0}) & =
   \frac{|x_{0}-x^{\star}_{0}|^{2}+ |t_{0}-t^{\star}_{0}|^{2}}{2
  \varepsilon} + v(x^{\star}_{0},t^{\star}_{0}) \\
\varphi(x,t) & \leq   v_{\varepsilon}(x,t)  & \leq 
  \frac{|x-y|^{2}+ |t-\tau|^{2}}{2
  \varepsilon} + v(y,\tau)
\end{eqnarray*}
we can verify that the function
$$\psi(x,t) = \varphi(x+x_{0}-x^{\star}_{0},t+t_{0}-t^{\star}_{0})
 -\frac{|x_{0}-x^{\star}_{0}|^{2} +|t_{0}-t^{\star}_{0}|^{2}}{2
  \varepsilon} $$
touches the original function $v$ from below at the point
$(x^{\star}_{0},t^{\star}_{0})$. Since $(x^{\star}_{0},t^{\star}_{0})$
is
 an interior point, the
inequality
$$\divergence\! \left(|\nabla \psi(x^{\star}_{0},t^{\star}_{0})|^{p-2}\nabla
  \psi(x^{\star}_{0},t^{\star}_{0})\right)
 \leq \psi_{t}(x^{\star}_{0},t^{\star}_{0})$$
holds by assumption. Because
$$\psi_{t}(x^{\star}_{0},t^{\star}_{0}) =
\varphi_{t}(x_{0},t_{0}),\ \ 
  \nabla\psi(x^{\star}_{0},t^{\star}_{0}) = \nabla\varphi(x_{0},t_{0}),
\ \ D^{2}\psi(x^{\star}_{0},t^{\star}_{0}) = D^{2}\varphi(x_{0},t_{0})$$
we also have that
$$\divergence\! \left(|\nabla \varphi(x_{0},t_{0})|^{p-2}\nabla
  \varphi(x_{0},t_{0})\right) \leq  \varphi_{t}(x_{0},t_{0}) $$
at the original point $(x_{0},t_{0})$, where $\varphi$ was touching
$v_{\varepsilon}$. Thus $v_{\varepsilon}$ fulfills the requirement in
the definition.

\textsf{Proof by Comparison}: We have to verify the comparison
principle for  $v_{\varepsilon}$ in a subcylinder $D_{t_{1},t_{2}}$
having at least the distance $\sqrt{2L\varepsilon}$ to the boundary of
$\Omega_{T}$. To this end, assume that $h \in
C(\overline{D_{t_{1},t_{2}}})$ is a $p$-caloric function such that
$v_{\varepsilon} \geq h$ on the parabolic boundary. It follows that
the inequality
\begin{equation*}
\frac{|x-y|^{2}+ |t-\tau|^{2}}{2 \varepsilon} + v(y,\tau) \geq h(x,t)
\end{equation*}
is available when $(y,\tau) \in \Omega_{T}$ and $(x,t)$ is on the
parabolic boundary of $D_{t_{1},t_{2}}$. Fix an arbitrary point
$(x_{0},t_{0})$ in $D_{t_{1},t_{2}}$. Then we can take $y =
x+x^{\star}_{0}-x_{0}$ and $\tau = t+t^{\star}_{0} -t_{0}$ in the
inequality above so that
\begin{gather*}
w(x,t) \equiv v(x+x^{\star}_{0}-x_{0},t+ t^{\star}_{0} -t_{0}) +
 \frac{|x_{0}-x^{\star}_{0}|^{2}+ |t_{0}-t^{\star}_{0}|^{2}}{2
  \varepsilon} \geq h(x,t)
\end{gather*}
when $(x,t)$ is on the parabolic boundary. But the translated function
$w$ is $p$-supercaloric in the subcylinder $D_{t_{1},t_{2}}$. By the
comparison principle $w \geq h$  in the whole subcylinder. In
particular,
$$v_{\varepsilon}(x_{0},t_{0}) = w(x_{0},t_{0}) \geq h(x_{0},t_{0}).$$
This proves the comparison principle for $v_{\varepsilon}$. $\Box$

 \medskip

The ``viscosity proof'' did not contain any explicite comparison principle while
the ``proof by comparison'' required the piece of knowledge that the
original $v$ obeys the principle.  This
\emph{parabolic} comparison principle allows comparison in space-time
cylinders. We will encounter domains of a more general shape, but the
following \emph{elliptic} version of the principle turns out to be
enough for our purpose. Instead of the expected parabolic boundary,
the whole boundary (the ``Euclidean'' boundary) appears.

\begin{prop} Given a domain $\Upsilon \subset \subset
  \Omega_{\varepsilon}$
and a $p$-caloric function $h \in C(\overline{\Upsilon})$, then
$v_{\varepsilon} \geq h$ on the whole boundary $\partial \Upsilon$
implies that $v_{\varepsilon} \geq h$ in $\Upsilon$.
\end{prop}

Now $\Upsilon$ does not have to be a space-time cylinder and $\partial
\Upsilon$ is the total boundary in $\mathbb{R}^{n+1}$.

\emph{Proof:} It is enough to realize that the proof is immediate when
$\Upsilon$ is a finite union of space-time cylinders $D_{j} \times
(a_{j},b_{j})$. To verify this, just start with the earliest
cylinder(s) and pay due attention to the passages of $t$ over the $a_{j}$'s
and $b_{j}$'s. Then the general case follows by exhausting  $\Upsilon$
with such unions. Indeed, given $\alpha > 0$ the compact set $\{h(x,t)
\geq v_{\varepsilon}(x,t)\}$ is contained in an open finite union 
$$\bigcup D_{j} \times (a_{j},b_{j})$$
comprised in $\Upsilon$ so that $h < v_{\varepsilon} + \alpha$ on the
(Euclidean) boundary of the union. It follows that  $h \leq
v_{\varepsilon} + \alpha$ in the whole union. Since $\alpha$ was
arbitrary, we conclude that $v_{\varepsilon} \geq h$ in
$\Upsilon$. $\Box$

\medskip
The above \emph{elliptic} comparison principle does not acknowledge
the presence of the parabolic boundary. The reasoning above can easily
be changed so that the latest portion of the boundary is
exempted. For this improvement, suppose that $t< T^{\star}$ for all $(x,t) \in \Upsilon$; in
this case $\partial \Upsilon$ may have a plane portion with
$t=T^{\star}$. It is now sufficient to verify that
$$v_{\varepsilon} \geq h\quad \text{on}\quad \partial \Upsilon \quad
\text{when} \quad t
<  T^{\star}$$
in order to conclude that $v_{\varepsilon} \geq h$ in
$\Upsilon$. To see this, just use 
$$v_{\varepsilon} + \frac{\sigma}{T^{\star} - t}$$
in the place of $v_{\varepsilon}$ and then let $\sigma \rightarrow
0+$. This variant of the comparison principle is convenient for the proof
of the following conclusion.

\begin{lemma}
The approximant $v_{\varepsilon}$ is a weak supersolution in the
shrunken domain, i.e.
\begin{equation}\label{shrunken}
\integ{\Bigl(- v_{\varepsilon}\frac{\partial \phi}{\partial t} +
    \langle |\nabla v_{\varepsilon}|^{p-2}\nabla
    v_{\varepsilon},\nabla \phi \rangle \Bigr)}  \geq 0 \end{equation}
whenever $\phi \in C_{0}^{\infty}(\Omega_{\varepsilon} \times
(\varepsilon, T - \varepsilon )), \ \phi \geq 0.$
\end{lemma}

\emph{Proof:} We show that in a given subdomain $D_{t_{1},t_{2}} = D
\times (t_{1},t_{2})$ of the ``shrunken domain'' our $v_{\varepsilon}$
coincides with the solution of an obstacle problem. The solutions of
the obstacle problem are \emph{per se} weak supersolutions. Hence, so
is $v_{\varepsilon}$. Consider the class  of all functions
\begin{equation*}
\begin{cases}
w \in C(\overline{D_{t_{1},t_{2}}}) \cap
L^{p}(t_{1},t_{2};W^{1,p}(D)),\\
w \geq v_{\varepsilon} \ \text{in} \ D_{t_{1},t_{2}}, \ \text{and}\\
w = v_{\varepsilon} \ \text{on the parabolic boundary of} \
D_{t_{1},t_{2}}.
\end{cases}
\end{equation*}
The function  $v_{\varepsilon}$ itself acts as an obstacle and induces
the boundary values. There exists a (unique) weak supersolution
$w_{\varepsilon}$ in this class satsfying the variational inequality
\begin{gather*}
\int_{t_{1}}^{t_{2}}\!\int_{D}\Bigl((\psi -
  w_{\varepsilon})\frac{\partial \psi}{\partial t} + \langle |\nabla
  w_{\varepsilon}|^{p-2}\nabla w_{\varepsilon},\nabla (\psi -
  w_{\varepsilon})\rangle \Bigr)dx\,dt\\
\geq \frac{1}{2}\int_{D}(\psi(x,t_{2})
-w_{\varepsilon}(x,t_{2}))^{2}\,dx
\end{gather*}
for all smooth $\psi$ in the aforementioned class. Moreover,
$w_{\varepsilon}$ is $p$-caloric in the open set $A_{\varepsilon} =
\{w_{\varepsilon} > v_{\varepsilon} \}$, where the obstacle does not
hinder. We refer to [C].
On the boundary $\partial A_{\varepsilon}$ we know that
$w_{\varepsilon} = v_{\varepsilon}$ except possibly when $t =
t_{2}$. By the elliptic comparison principle we have $v_{\varepsilon}
\geq  w_{\varepsilon}$ in  $A_{\varepsilon}$. On the other hand $w_{\varepsilon}
\geq  v_{\varepsilon}$. Hence  $w_{\varepsilon}
=  v_{\varepsilon}$.

To finish the proof, 
let $\varphi \in C_{0}^{\infty}(D_{t_{1},t_{2}}), \, \varphi \geq 0$,
and choose $\psi = v_{\varepsilon} +  \varphi =  w_{\varepsilon} +
\varphi$ above. ---An easy manipulation yields (\ref{shrunken}). $\Box$

\medskip
Recall that $0 \leq v \leq L$. Then also $0 \leq v_{\varepsilon}  \leq
L$. An estimate for $\nabla v_{\varepsilon}$ is provided in the
well-known lemma below.

\begin{lemma}[Caccioppoli] The inequality
\begin{eqnarray*}
\integ{\zeta^{p}|\nabla v_{\varepsilon}|^{p}} \leq C L^{p}
\integ{|\nabla \zeta|^{p}} + C L^{2} \integ{\Bigl|\frac{\partial
      \zeta^{p}}{\partial t}\Bigr|}
\end{eqnarray*}
 holds whenever  $\zeta \in C_{0}^{\infty}(\Omega_{\varepsilon} \times
(\varepsilon, T - \varepsilon)), \ \zeta \geq 0.$
\end{lemma}

\emph{Proof:} Use the test function 
$$\phi(x,t) = (L -  v_{\varepsilon}(x,t)) \zeta^{p}(x,t).\quad \Box $$

The Caccioppoli estimate above leads to the conclusion that, keeping
$0 \leq v \leq L$, the Sobolev gradient $\nabla v \in L_{loc}^{p}$
exists and
$$ \nabla  v_{\varepsilon} \rightharpoonup \nabla v  \text{\ weakly in
  \ } L_{loc}^{p},$$
at least for a subsequence. For $v$ the Caccioppoli estimate
\begin{eqnarray*}
\integ{\zeta^{p}|\nabla v|^{p}} \leq C L^{p}
\integ{|\nabla \zeta|^{p}} + C L^{2} \integ{\Bigl|\frac{\partial
      \zeta^{p}}{\partial t}\Bigr|},
\end{eqnarray*}
is immediate, because of the lower semicontinuity of the integral
under weak convergence. However the corresponding passage to the limit
under the integral sign of
$$\integ{\Bigl(- v_{\varepsilon}\frac{\partial \phi}{\partial t} +
    \langle |\nabla v_{\varepsilon}|^{p-2}\nabla
    v_{\varepsilon},\nabla \phi \rangle \Bigr)}  \geq 0$$
requires a justification, as $\varepsilon \rightarrow 0.$ The
elementary vector inequality
\begin{equation}
\left| |b|^{p-2}b - |a|^{p-2}a \right|
\leq (p-1)|b -a|\left(|b| + |a| \right)^{p-2},
\end{equation}
$p > 2$, and H\"{o}lder's inequality show that it is sufficient to
establish that
$$\nabla  v_{\varepsilon} \longrightarrow \nabla v \text{\ \ strongly
  in\ \ } L_{loc}^{p-1},$$
to ackomplish the passage. Notice the exponent $p-1$ in place of $p.$
This strong convergence is given in the next theorem, where the
sequence is renamed to $v_{k}.$

\begin{thm}
\label{linmanf}
Suppose that $v_{1},v_{2},v_{3},\dots$ is a sequence of Lipschitz
continuous weak supersolutions, such that
$$0 \leq v_{k} \leq L \text{ \ in \ } \Omega_{T} =   \Omega
 \times (0,t),\ \  v_{k}
\rightarrow v \text{\ in \ } L^{p}(\Omega_{T}).$$
Then
$$\nabla v_{1},\nabla v_{2},\nabla v_{3},\dots$$
is a Cauchy sequence in $ L_{loc}^{p-1}(\Omega_{T}).$
\end{thm}

\emph{Proof:} The central idea is that the measure of the set where
$|v_{j} - v_{k}| > \delta$ is small. Given $\delta > 0$, we have, in
fact, 
\begin{equation}
\label{mitta}
\textrm{mes}\{|v_{j} - v_{k}| > \delta \} 
\leq \delta^{-p}\|v_{j} - v_{k}\|_{p}^{p}
\end{equation}
according to Tshebyshef's inequality. Fix a test function $\theta \in
C_{0}^{\infty}(\Omega_{T}), \ 0 \leq \theta \leq 1.$ From the Caccioppoli
estimate we can extract a bound of the form
$$\int\!\!\int_{\{\theta \not = 0\}} |\nabla v_{k}|^{p}\,dx\,dt \leq A^{p}, \ \
\ k = 1,2,\ldots ,$$since the support is a compact subset.
Fix the indices $k$ and $j$ and use the test function
$$\varphi = (\delta - w_{jk}) \theta$$ where
$$ w_{jk} = \left\{\begin{array}{ll}
\delta, & \mbox{if $v_{j} -v_{k} > \delta$}\\
 v_{j} -v_{k}, & \mbox{if $|v_{j} -v_{k}| \leq \delta$}\\
- \delta, & \mbox{if $v_{j} -v_{k} < - \delta$}
\end{array} \right. $$
in the equation
$$\integ{\Bigl(- v_{j}\frac{\partial \phi}{\partial t} +
    \langle |\nabla v_{j}|^{p-2}\nabla
    v_{j},\nabla \phi \rangle \!\Bigr)}  \geq 0.$$ Since $|w_{jk}| \leq
\delta$, we have $\varphi \geq 0$. In the equation for $v_{k}$ use
$$ \varphi = (\delta +w_{jk}) \theta.$$
Add the two equations and arrange the terms:
\begin{eqnarray*}
\int \! \! \int_{|v_{j} -v_{k}| \leq \delta} \theta \langle |\nabla
v_{j}|^{p-2}\nabla v_{j} -|\nabla v_{k}|^{p-2}\nabla v_{k}, \nabla
v_{j} - \nabla v_{k}\rangle \,dx \, dt \\
\leq \delta \integ{ \langle |\nabla
v_{j}|^{p-2}\nabla v_{j} + |\nabla v_{k}|^{p-2}\nabla v_{k}, \nabla
\theta \rangle } \\
- \integ{w_{jk} \langle |\nabla
v_{j}|^{p-2}\nabla v_{j} - |\nabla v_{k}|^{p-2}\nabla v_{k}, \nabla
\theta \rangle } \\
+ \integ{(v_{j} - v_{k}) \frac{\partial}{\partial t}(\theta w_{jk})}
- \delta \integ{(v_{j} + v_{k}) \frac{\partial \theta}{\partial t}}\\
= I - II + III - IV.
\end{eqnarray*}
The left-hand side is familiar from inequality (\ref{Vect}). As we will see,
the right-hand side is of the magnitude $O(\delta)$. We begin with term
III, which contains
 time derivatives that ought to be
avoided. Integration by parts yields
\begin{eqnarray*}
III = \integ{\theta  \frac{\partial \theta}{\partial
    t}(\frac{w_{jk}^{2}}{2})}
+ \integ{(v_{j} - v_{k})w_{jk} \frac{\partial \theta}{\partial t}} \\
= - \frac{1}{2} \integ{w_{jk}^{2}\frac{\partial \theta}{\partial t}} +
\integ{(v_{j} - v_{k})w_{jk} \frac{\partial \theta}{\partial t}}.
\end{eqnarray*}
We obtain the estimate
$$|III| \leq \frac{1}{2} \delta^{2} \|\theta_{t}\|_{1} + 2 L \delta 
\|\theta_{t}\|_{1} \leq \delta C_{3}.$$
For the last term we immediately have
$$|IV| \leq 2 \delta L  \|\theta_{t}\|_{1} = \delta C_{4}.$$
The two first terms are easy,
$$|I| \leq \delta C_{1}, \ \ \ |II| \leq \delta C_{2}.$$
Summing up,
$$|I| + |II| + |III| + |IV| \leq C \delta.$$
Using the vector inequality (\ref{Vect}) to estimate the left hand side, we
arrive at
\begin{eqnarray*}
\int  \! \int_{|v_{j}-v_{k}| \leq \delta} \theta |\nabla v_{j}
-\nabla v_{k}|^{p} \, dx \, dt &\leq 2^{p-2} \delta C,\\
\int  \! \int_{|v_{j}-v_{k}| \leq \delta} \theta |\nabla v_{j}
-\nabla v_{k}|^{p-1} \, dx \, dt &= O(\delta ^{1- \frac{1}{p}}).
\end{eqnarray*}
We also have in virtue of (\ref{mitta}) 
\begin{gather*}
\int  \! \int_{|v_{j}-v_{k}| > \delta} \theta |\nabla v_{j}
-\nabla v_{k}|^{p-1} \, dx \, dt \\
\leq \delta ^{-1} \|v_{j}-v_{k}\|_{p} \left(\|\nabla v_{j}\|_{p} +
 \|\nabla v_{k}\|_{p} \right)^{p-1}
\leq (2 A)^{p-1} \delta ^{-1} \|v_{j}-v_{k}\|_{p}.
\end{gather*}
Finally, combining the estimates over the sets $|v_{j}-v_{k}| \leq
\delta$ and  $|v_{j}-v_{k}| > \delta$, we have an integral over the
whole domain:
$$\integ{\theta |\nabla v_{j} - \nabla v_{k}|^{p-1}} \leq
 O(\delta^{1- \frac{1}{p}}) + C_{5} \delta^{-1}\|v_{j} -
 v_{k}\|_{p}.$$
Since the left-hand side is independent of $\delta$, we can make it as
small as we please, by first fixing $\delta$ small enough and then
taking the indices large enough. $\Box$

\medskip 

We have arrived at the following result for bounded supersolutions.
\begin{thm}
\label{Thm25}
Let $v$ be a \emph{bounded} $p$-supercaloric function.  Then 
$$\nabla v = \left(\frac{\partial v}{\partial x_{1}}, \ldots ,
\frac{\partial v}{\partial x_{n}} \right)$$
exists in Sobolev's sense, $\nabla v \in L_{loc}^{p}(\Omega_{T}),$ and
$$\integ{\Bigl(- v \frac{\partial \phi}{\partial t} +
    \langle |\nabla v |^{p-2}\nabla
    v,\nabla \phi \rangle \Bigr)}  \geq 0$$
for all non-negative compactly supported test functions $\varphi.$
\end{thm}

Notice that so far the exponent is $p$, as it should for $v$ bounded.

\medskip
\emph{Remark:} It was established that the $p$-supercaloric functions
 are also weak supersolutions. In the opposite
direction, according to [K] every weak supersolution is lower
semicontinuous upon a redefinition in a set of (n+1)-dimensional measure
$0$. Moreover, the representative obtained as
$$\underset{(y,\tau) \rightarrow (x,t)}
{\textrm{ess\,liminf}}\,
v(y,\tau)$$
will do. A proof is given in Chapter 4.

\medskip

 We need a few auxiliary
results.

\begin{lemma}[Sobolev's inequality]
\label{Sobolev}
 If $u \in
  L^{p}(0,T;W^{1,p}_{0}(\Omega))$, then
\begin{equation*}
\integ{|u|^{p(1+ \frac{2}{n})}} \leq \ 
C \integ{|\nabla u|^{p}} \left\{\underset{0<t<T}{\textrm{ess\,sup}}
  \int_{\Omega}{|u(x,t)|^{2}\,dx}
\right\}^{\frac{p}{n}}.
\end{equation*}
\end{lemma}.

\emph{Proof:} See for example [dB, Chapter 1, Proposition 3.1].

\medskip

If the test function $\phi$ is zero on the lateral boundary $\partial
\Omega \times [t_{1},t_{2}]$, then the differential inequality for
the weak supersolution takes the form
\begin{eqnarray*}
\int_{t_{1}}^{t_{2}}\!\int_{\Omega} \Bigl(- v \frac{\partial \phi}{\partial t} +
    \langle |\nabla v |^{p-2}\nabla
    v,\nabla \phi \rangle \Bigr)\,dx\,dt  \\
+ \int_{\Omega}v(x,t_{2})\phi(x,t_{2})\,dx \geq 
\int_{\Omega}v(x,t_{1})\phi(x,t_{1})\,dx.
\end{eqnarray*}
Thus, if $v$ is zero on the lateral boundary, we may take $\phi = v$
above. We obtain
\begin{equation}
\label{future}
\frac{1}{2} \int_{\Omega} v(x,t_{1})^{2} \,dx \leq 
\int_{t_{1}}^{t_{2}}\!\int_{\Omega}|\nabla v|^{p}\,dx\,dt
 + \frac{1}{2} \int_{\Omega} v(x,t_{2})^{2} \,dx,
\end{equation}
which estimates the \emph{past} in terms of the \emph{future} and an
``energy term''.

\subsection{Harnack's Convergence Theorem}

The classical convergence theorem of Harnack states that the limit of an increasing sequence of harmonic functions is either a harmonic function itself or identically $+\infty.$ The convergence is locally uniform. The situation is similar for many elliptic equations. However, the Evolutionary $p$-Laplace Equation exhibits a more delicate behaviour. The limit function can be finite at each point without being even locally bounded!\footnote{This was unfortunately overlooked in [KL] and [KL2]. Corrections appear in [KP] and [KL3].} This is a characteristic feature for the class $\mathfrak{M}$ previously defined.
Consider a sequence of nonnegative $p$-caloric functions 
$$0 \leq h_1 \leq h_2 \leq h_3\leq\,\dots\qquad h = \lim_{k\to\infty}h_k$$
in $\Omega_T.$ There are two different possibilities, depending on whether the limit function $h$ is locally bounded or not. The basic tool is an intrinsic version of Harnack's inequality, which is due to E. DiBenedetto, see [dB2, pp. 157--158] or [dB1].

\bigskip

\begin{lemma}[Harnack]
\label{Harnackpar}
Let $p >2$. There are constants $C$ and $\gamma$, depending only on $n$ and $p$, such that if $u>0$ is a continuous weak solution in
$$B(x_0,4R)\times(t_0-4\theta,t_0+4\theta),\quad \text{where}\quad \theta =\frac{CR^p}{u(x_0,t_0)^{p-2}}$$
then the inequality
\begin{equation}
\label{hartheta}
u(x_0,t_0)\,\leq\,\gamma \inf_{B(x_0,R)}u(x,t_0+\theta)
\end{equation}
is valid.
\end{lemma}

\bigskip

Notice how the waiting time $\theta$ depends on the solution itself. It is very short, if the solution is large.

\bigskip

\begin{prop}
\label{blowup}
Suppose that we have an increasing sequence 
$$0 \leq h_1 \leq h_2 \leq h_3\leq\,\dots\qquad h = \lim_{k\to\infty}h_k$$
of $p$-caloric functions $h_k.$ If for some sequence
$$h_k(x_k,t_k) \to + \infty\quad\text{and}\quad (x_k,t_k) \to(x_0,t_0)$$
where $x_0\in \Omega,\,0<t_0<T,$ then
$$\liminf_{\substack{(y,t)\to(x,t_0)\\t>t_0}}h(y,t)(t-t_0)^{\frac{1}{p-2}}\,>\,0\quad\text{for all}\quad x\in\Omega.$$
Thus, at time $t_0$,
$$\lim_{\substack{(y,t)\to(x,t_0)\\t>t_0}}h(y,t)\,\equiv\,\infty\quad\text{in}\quad \Omega.$$
\end{prop}

\bigskip

\emph{Proof:} Let $B(x_0,4R)\subset\subset \Omega.$ With
$$\theta_k = \frac{CR^p}{h_k(x_k,t_k)^{p-2}} \longrightarrow 0$$
we have by Harnack's inequality (\ref{hartheta})
\begin{equation}
\label{gamma}
h_k(x_k,t_k)\leq \gamma h_k(x,t_k+\theta_k)
\end{equation}
when $x\in B(x_k,R)$ provided that $B(x_k,4R)\times (t_k-4\theta_k,t_k+4\theta_k)\subset\subset \Omega_T.$ The center is moving, but since $x_k \to x_0,$ equation (\ref{gamma}) holds for sufficiently large indices $k.$ Let $\Gamma > 1.$
We want to compare the two solutions
$$\dfrac{\mathfrak{U}^R(x)}{\bigl(t-t_k+(\Gamma-1)\theta_k\bigr)^{\frac{1}{p-2}}}\qquad\text{and}\qquad h_k(x,t)$$
when $t=t_k+\theta_k$ and $x \in B(x_0,R).$ Here $\mathfrak{U}^R$ is the positive solution of the elliptic equation (\ref{apuyh}) with boundary values zero on $\partial B(x_0,R)$. By \ref{gamma} we have
\begin{align*}
&\dfrac{\mathfrak{U}^R(x)}{\bigl(t-t_k+(\Gamma-1)\theta_k\bigr)^{\frac{1}{p-2}}}\Bigg\vert_{t=t_k+\theta_k}\,=\, \,\dfrac{\mathfrak{U}^R(x)}{(\Gamma CR^p)^{\frac{1}{p-2}}}\,h_k(x_k,t_k)\\
&\leq \dfrac{\mathfrak{U}^R(x)}{(\Gamma CR^p)^{\frac{1}{p-2}}}\,\gamma\, h_k(x,t_k+\theta_k)\,\,\leq \,\,\,h_k(x,t_k+\theta_k)
\end{align*}
if we fix $\Gamma$ so large that
$$\gamma \|\mathfrak{U}^R\|_{L^{\infty}(B(x_0,R)}\,\leq\,(\Gamma CR^p)^{\frac{1}{p-2}}.$$
By the Comparison Principle
$$\dfrac{\mathfrak{U}^R(x)}{\bigl(t-t_k+(\Gamma-1)\theta_k\bigr)^{\frac{1}{p-2}}}\leq h_k(x,t)\leq h(x,t)$$
when $t \geq t_k+\theta_k$ and $x\in B(x_0,R).$ Sending $k$ to $\infty$, we arrive at
$$\dfrac{\mathfrak{U}^R(x)}{(t-t_0)^{\frac{1}{p-2}}}\leq h(x,t)\quad\text{when}\quad t_0 < t < T.$$
This yields the desired estimate, though only in a subdomain. Then
 repeat the procedure in suitably chosen balls, thus   extending the estimate  to the entire domain $\Omega.$ \qquad $\Box$

\bigskip

\begin{prop}
\label{increasing}
Suppose that we have an increasing sequence 
$$0 \leq h_1 \leq h_2 \leq h_3\leq\,\dots\qquad h = \lim_{k\to\infty}h_k$$
of $p$-caloric functions $h_k$  in $\Omega_T.$ If the sequence is locally bounded, then the limit function $h$ is $p$-caloric in $\Omega_T.$
\end{prop}

\bigskip

\emph{Remark:} The situation is delicate. It is not enough to assume that $h$ is finite at every point. This is different for elliptic equations! So is it for the Heat Equation.

\bigskip

\emph{Proof:} In a fixed strict subdomain we have H\"{o}lder continuity
$$|h_k(x_1,t_1)-h_k(x_2,t_2)| \leq C\|h_k\|\left(|x_1-x_2|^{\alpha}+|t_1-t_2|^{\frac{\alpha}{p}}\right).$$
Here $\|h_k\| \leq \|h\| <  \infty$ so that the family is locally equicontinuous. Hence the convergence $h_k\to h$ is locally uniform in $\Omega_T$ and, consequently, the limit function $h$ is continuous. 

 From the usual Caccioppoli estimate
\begin{align*}
&\int_{t_1}^{t_2}\!\!\int_{\Omega}\zeta^p|\nabla h_k|^p \,d x d t\\
&\leq C(p)\int_{t_1}^{t_2}\!\!\int_{\Omega}h_k^p|\nabla \zeta|^p\,d x d t + C(p)\int_{\Omega}\zeta(x)^ph_k(x,t)^2\bigg|_{t_1}^{t_2} \,d x\\
&\leq C(p)\int_{t_1}^{t_2}\!\!\int_{\Omega}h^p|\nabla \zeta|^p\,d x d t + C(p)\int_{\Omega}\zeta(x)^ph(x,t)^2\bigg|_{t_1}^{t_2} \,d x
\end{align*}
we can, in a standard way, conclude that $h \in L^p_{loc}(0,T;W^{1,p}_{loc}(\Omega)).$
It is easy to see that $h$ satisfies the Comparison Principle in $\Omega_T,$  since each $h_k$ does it and the convergence is uniform. From Theorem 24 we conclude that the equation
$$\int_0^T\!\!\int_{\Omega}\Bigl(-h\frac{\partial\varphi}{\partial t}+\langle|\nabla h|^{p-2}\nabla h,\nabla \varphi\rangle\Bigr)\, dx  dt\,=\,0$$
is valid.\qquad
$\Box$

\bigskip

\subsection{Unbounded Supersolutions}

We proceed to study an \emph{un}bounded 
$p$-supercaloric function $v$. Let us briefly describe the
method. The starting point is to apply Theorem \ref{Thm25} on the
functions $v_{k} = \min\{v,k\}$ so that estimates depending on $k =1,2,\cdots$
are obtained. To begin with, it is crucial that
$$v_k \in L^p(0,T;W^{1,p}_0(\Omega)).$$
Then an iterative procedure is used to
gradually increase the summability exponent of $v$. First, we achieve
that $v^{\alpha}$ is locally summable for some small exponent $\alpha < p-2$. That result
is iterated, again using the $v_{k}$'s till we come close to the
exponent $\alpha = p-1-0$. Then the passage over $p-1$ requires a
special, although simple, device. At the end we will reach the
desired summability for the function $v$ itself. From this it is not
difficult to obtain the corresponding result also for the gradient
$\nabla v$. Again the $v_{k}$'s are employed. Finally, one has the problem to remove the restriction about zero lateral boundary values. This is not possible for functions of class $\mathfrak{M}$\footnote{Their infinities always hit the lateral boundary.}. For class $\mathfrak{B}$, this is done in Section \ref{Proof}. 

 The considerations are in a bounded
subdomain, which we again call
$\Omega_{T} = \Omega \times (0,T)$, for simplicity. The situation is much
 easier when the function is zero on the whole parabolic boundary:
$$v(x,0) = 0 \, \text{\ when\ } \, x \in \Omega,\ \ \ v = 0 \, \text{\ on\ }\,
\partial \Omega \times [0,T].$$ We  assume that $v \geq 0.$ The functions 
$$v_{k}(x,t) = \min \{v(x,t),k \}, \ \ \ \ k = 0,1,2,\ldots,$$
cut off at the height $k$ are bounded, whence the previous results in Section \ref{Bounded Supersolutions}
apply for them. Fix a
large index $j$.  We may use the test functions
$$\phi_{k} = (v_{k} - v_{k-1}) - (v_{k+1} - v_{k}), \ \ k =
1,2,\ldots,j$$
in the equation 
$$\int_{0}^{\tau}\!\!\int_{\Omega}\langle |\nabla v_{j}|^{p-2}\nabla
v_{j},\nabla \phi_{j}\rangle \,dx\,dt +
\int_{0}^{\tau}\!\!\int_{\Omega}\phi_{k}\frac{\partial v_{j}}{\partial
  t}\,dx\,dt \geq 0,$$
where $0 < \tau \leq T.$ Indeed, $\phi_{k} \geq 0$.
   The ``forbidden'' time derivative can be
avoided through an appropriate regularization. In
principle $v_{j}$ is first replaced by its convolution with a
mollifier
 and later the limit is to be taken. We postpone this complication
 in order to keep the exposition more transparent. The insertion
of the test function yields
\begin{align*}
&\int_{0}^{\tau}\!\!\int_{\Omega}\Bigl(\langle |\nabla v_{j}|^{p-2}\nabla
v_{j},\nabla (v_{k+1} - v_{k})\rangle \ +
 (v_{k+1} - v_{k}) \frac{\partial v_{j}}{\partial  t}\Bigr) dx\,dt \\
\leq&\int_{0}^{\tau}\!\!\int_{\Omega}\Bigl(\langle |\nabla v_{j}|^{p-2}\nabla
v_{j},\nabla (v_{k} - v_{k-1})\rangle \ +
 (v_{k} - v_{k-1}) \frac{\partial v_{j}}{\partial  t}\Bigr) dx\,dt,
\end{align*}
 succinctly written as
$$a_{k+1}(\tau) \leq a_{k}(\tau).$$
It follows that
$$\sum_{k=1}^{j}a_{k}(\tau) \leq j a_{1}(\tau)$$
and, since the sum is ``telescoping'', we have the result below.
\begin{lemma}
\label{eka}
 If each $v_{k} \in L^{p}(0,T;W^{1,p}_{0}(\Omega))$ and
  $v_{k}(x,0) = 0$ when $x \in \Omega$, then
\begin{align*}
&\int_{0}^{\tau}\!\!\int_{\Omega}|\nabla v_{j}|^{p}\,dx\,dt +
  \frac{1}{2} \int_{\Omega}v_{j}^{2}(x,\tau)\,dx\\ \leq
j& \int_{0}^{\tau}\!\!\int_{\Omega}|\nabla v_{1}|^{p}\,dx\,dt +
  j \int_{\Omega}v_{j}(x,\tau)\,dx
\end{align*}
holds for a.e. $\tau$ in the range $0 < \tau \leq T$. Consequently,
\begin{align*}
&\int_{0}^{T}\!\!\int_{\Omega}|\nabla v_{j}|^{p}\,dx\,dt +
  \frac{1}{2}\sup_{0<t<T} \int_{\Omega}v_{j}^{2}(x,t)\,dx\\ &\leq
2j^2 \Bigl(\int_{0}^{T}\!\!\int_{\Omega}|\nabla v_{1}|^{p}\,dx\,dt +
  |\Omega|\Bigr).
\end{align*}
\end{lemma}

\bigskip

Before continuing, we justify the use of the time
derivative in the previous reasoning.

 {\small \textsf{Regularisation
    of the equation.} We use the convolution
\begin{equation}
\label{convo}
(f \star \rho_{\varepsilon})(x,t) =
\int_{-\infty}^{\infty}f(x,t-s)
\rho_{\varepsilon}(s) \,ds, \end{equation}
where $\rho_{\varepsilon}$ is, for instance, Friedrich's
mollifier defined as
$$
\rho_{\varepsilon}(t) = 
\begin{cases}
\frac{C}{\varepsilon}e^{-\varepsilon^2/(\varepsilon^2-t^2)}, &|t| <
  \varepsilon,\\
0, &|t| \geq \varepsilon.
\end{cases}$$
 If the function $v_{j}$ is extended as $0$ when $t \leq 0,\, x
\in \Omega$, the new function is $p$-supercaloric in $\Omega \times
(-\infty,T)$, indeed. To see this, one has only to verify the
comparison principle.

 We have, when $\tau \leq T - \varepsilon$
$$\int_{-\infty}^{\tau}\!\int_{\Omega}\left(\langle(|\nabla
  v_{j}|^{p-2}\nabla v_{j})\star \rho_{\varepsilon},\nabla \varphi
  \rangle + \varphi \frac{\partial}{\partial t}( v_{j} \star 
\rho_{\varepsilon}) \right)\,dx\,dt \geq 0$$
for all test functions $\varphi \geq 0$ vanishing on the lateral
boundary. Replace $v_{k}$ in the previous proof by
$$\tilde{v}_{k} = \min\{v_{j}\star \rho_{\varepsilon},k\}$$
and choose
$$\varphi_{k} = (\tilde{v}_{k} -\tilde{v}_{k-1}) - (\tilde{v}_{k+1}
 -\tilde{v}_{k}).$$
Since the convolution with respect to the time variable does not
affect the zero boundary values on the lateral boundary, we conclude
 that $\tilde{v}_{k}$ vanishes on the parabolic boundary of $\Omega
\times (-\delta/2, T-\delta/2),$ when $\varepsilon < \delta/2$ and
$\delta$ can be taken as small as we wish. (The functions $v_{k}\star
\rho_{\varepsilon}$ instead of the employed $(v\star
\rho_{\varepsilon})_{k}$ do not work well in this proof.) The same
calculations as before yield
$$\tilde{a}_{k+1}(\tau) \leq \tilde{a}_{k}(\tau) \quad \text{and}
\quad \sum_{k=1}^{j}\tilde{a}_{k}(\tau) \leq j \,\tilde{a}_{1}(\tau),$$
where 
\begin{gather*}
\tilde{a}_{k}(\tau) = \int_{-\delta/2}^{\tau}\!\int_{\Omega}
\langle(|\nabla
  v_{j}|^{p-2}\nabla v_{j})\star
  \rho_{\varepsilon},\nabla(\tilde{v}_{k}- \tilde{v}_{k-1})
  \rangle\,dx\,dt\\
 +\int_{-\delta/2}^{\tau}\!\int_{\Omega}(\tilde{v}_{k}-
 \tilde{v}_{k-1})
 \frac{\partial}{\partial t}( v_{j} \star 
\rho_{\varepsilon})\,dx\,dt .
\end{gather*}
Summing up, we obtain
\begin{gather*}
 \sum_{k=1}^{j}\tilde{a}_{k}(\tau)
   = \int_{-\delta/2}^{\tau}\!\int_{\Omega}
\langle(|\nabla
  v_{j}|^{p-2}\nabla v_{j})\star
  \rho_{\varepsilon},\nabla \tilde{v}_{j}
  \rangle\,dx\,dt\\
 +\int_{-\delta/2}^{\tau}\!\int_{\Omega}\tilde{v}_{j}
 \frac{\partial}{\partial t}( v_{j} \star 
\rho_{\varepsilon})\,dx\,dt ,
\end{gather*}
where the last integral can be written as 
$$\frac{1}{2} \int_{\Omega}(v_{j} \star
\rho_{\varepsilon})^{2}(x,\tau)\,dx.$$
Also for $\tilde{a}_{1}(\tau)$ we get an expression free of time
derivatives. Therefore we can safely first let $\varepsilon
\rightarrow 0$ and then $\delta \rightarrow 0$. This leads to
the lemma.}

\medskip
Let us return to the lemma.
Provided that we already have a majorant for the term
$$\int_{0}^{T}\!\!\int_{\Omega}|\nabla v_{1}|^{p}\,dx\,dt,$$ we see
that
$$\int_{0}^{T}\!\!\int_{\Omega}|\nabla v_{j}|^{p}\,dx\,dt =
O(j^{2}).$$
Yet, the right magnitude is $O(j).$
\begin{lemma}\label{Kjj}
 Suppose that $v_{j} \in L^{p}(0,T;W^{1,p}_{0}(\Omega))$
  and
$$\int_{0}^{T}\!\!\int_{\Omega}|\nabla v_{j}|^{p}\,dx\,dt \leq K
j^{2}, \ \ \ \ j = 1,2,3,\ldots$$
Then $v \in L^{q}(\Omega_{T}),$ whenever $q < p-2.$ (Here $p > 2$.)
\end{lemma}

\emph{Proof:} (Recall that the desired bound is $p - 1 + \frac{p}{n}$ and
not the above $p-2$.) The assumption and Sobolev's inequality (Lemma \ref{Sobolev}) will give
us a bound on the measure of the level sets
$$ E_{j} = \{(x,t)|\,j \leq v(x,t) \leq 2j \}$$
so that the integral can be controlled. To this end, denote
$$\kappa = 1 + \frac{2}{n}.$$ We have
\begin{gather*}
j^{\kappa p}|E_{j}| \leq \int \!\! \int_{E_{j}} v^{\kappa
  p}_{2j}\,dx\,dt \leq \integ{ v^{\kappa p}_{2j}}\\  \leq
C \integ{|\nabla v_{2j}|^{p}} \cdot \left(\underset{0<t<T}{\text{ess\,sup}}
\int_{\Omega}v_{2j}^{2}\,dx \right)^{\frac{p}{n}} \leq
C K j^{2}\left(4|\Omega|j^{2}\right)^{\frac{p}{n}}.
\end{gather*}
It follows that
$$ |E_{j}| \leq \text{Constant}\,j^{2-p}.$$

We use this to estimate the $L^{q}$-norm using a dyadic division of
the domain. Thus
\begin{align*}
\integ{v^{q}} \leq \ & T\, |\Omega| + \sum_{j=1}^{\infty}\int \!
\int_{E_{2^{j-1}}} v^{q}\,dx\,dt \\ \leq\ & T\, |\Omega| +
  \sum_{j=1}^{\infty} 2^{jq}|E_{2^{j-1}}|\\
\leq\ & T\, |\Omega| + C \sum_{j=1}^{\infty} 2^{j(q+2-p)},
\end{align*}
which is a convergent majorant when $q < p-2$. $\Box$

\medskip\noindent
\emph{Remark:} If the majorant $K j^{2}$ in the assumption is replaced
by a better $K j^{\gamma}$, then the procedure yields that $|E_{j}|
\approx j^{\gamma - p}$ resulting in $q < p - \gamma.$

\medskip
The previous lemma guarantees that $v^{\varepsilon}$ is summable for some
small positive power $\varepsilon$, since\footnote{It does not work
  for the Heat Equation!} $p > 2$. To improve the exponent,
we start from Lemma \ref{eka} and write the estimate in the form

\begin{equation*}
\int_{0}^{t_{1}}\!\!\int_{\Omega}|\nabla v_{j}|^{p}\,dx\,dt  \leq
j \int_{0}^{T}\!\!\int_{\Omega}|\nabla v_{1}|^{p}\,dx\,dt +
  j \int_{\Omega}v_{j}(x,\tau)\,dx,
\end{equation*}
where $0 < t_{1} \leq \tau \leq T$. Integrate with respect to $\tau$
over the interval $[t_{1},T]$:
\begin{align*}
(T-t_{1}) \int_{0}^{t_{1}}\!\!\int_{\Omega}|\nabla v_{j}|^{p}\,dx\,dt  \leq
j(T-t_{1})K +
  j \int_{t_{1}}^{T}\!\int_{\Omega}v_{j}(x,t)\,dx \,dt \\ \leq
j(T-t_{1})K +
  j^{2-\varepsilon} \int_{t_{1}}^{T}\!\int_{\Omega}v^{\varepsilon}(x,t)\,dx \,dt.
\end{align*}
Thus we have reached the estimate
\begin{equation}
 \int_{0}^{t_{1}}\!\!\int_{\Omega}|\nabla v_{j}|^{p}\,dx\,dt
  \leq j^{2-\varepsilon}K_{1}, \ \ \ \ \ (t_{1} < T).
\end{equation}
This is an improvement from $j^{2}$ to $j^{2-\varepsilon}$, but we
have to obey the restriction that $\varepsilon \leq 1$, because the
term $j(T-t_{1})K$ was absorbed. Estimating
again the measures $|E_{j}|$, but this time starting with the bound
$j^{\gamma}K,\  \gamma = 2 -\varepsilon$ in place of $j^{2}K$, yields
$$|E_{j}| \approx j^{2 - p -\varepsilon},\ \ \ \ \ \ 0 < \varepsilon
\leq 1.$$
The result is that 
$$\int_{0}^{t_{1}}\!\!\int_{\Omega} v^{q}\,dx\,dt < \infty \quad \text{\
  when \ } \quad 0 < q < p - \gamma = p - 2 + \varepsilon.$$
Iterating, we have the scheme
\begin{align*}
& q_{0} =  \varepsilon \ \ \ \ \ \ & T \\
& q_{1} = p - 2 + \varepsilon \ \       & t_{1} \\
& q_{2} = 2(p - 2) +  \varepsilon \     & t_{2} 
\end{align*}
We can continue till we reach
$$q_{k} = k(p - 2) +  \varepsilon  > p-1.$$
We have to stop, because the previous exponent $(k-1)(p-2) +
\varepsilon$ has to obey the rule not to become larger than $1$.
 This way we can
 reach that
 \begin{equation}
  \label{L1T} 
  v \in L^{1}(\Omega_{T'}),
  \end{equation}
with a $T' < T$, which will do to proceed.
 In fact, adjusting we can reach any
exponent strictly below $p-1$, but {\textsf the passage over the
exponent $p-1$} requires a special device. Since only a finite
number of steps were involved, we can take $T'$ as close to $T$ as we wish.

 We use inequality (\ref{future}) in the
form
\begin{equation}
\label{future1}
\frac{1}{2} \int_{\Omega} v_{j}(x,t)^{2} \,dx \leq 
\int_{0}^{\tau}\!\int_{\Omega}|\nabla v_{j}|^{p}\,dx\,dt
 + \frac{1}{2} \int_{\Omega} v_{j}(x,\tau)^{2} \,dx,
\end{equation}
where $t < \tau$. For $t_{1} <\tau < T$ it follows that
\begin{align*}
\underset{0<t<t_{1}}{\textrm{ess\,sup}} \int_{\Omega} v_{j}(x,t)^{2} \,dx \leq 
2 &\int_{0}^{\tau}\!\int_{\Omega}|\nabla v_{j}|^{p}\,dx\,dt
 + \int_{\Omega} v_{j}(x,\tau)^{2} \,dx\\ \leq
2 & j \int_{0}^{\tau}\!\int_{\Omega}|\nabla v_{1}|^{p}\,dx\,dt +
2 j \int_{\Omega} v_{j}(x,\tau) \,dx \\ \leq
2 & j \int_{0}^{T}\!\int_{\Omega}|\nabla v_{1}|^{p}\,dx\,dt +
2 j \int_{\Omega} v_{j}(x,\tau) \,dx,
\end{align*}
where the middle step is from Lemma \ref{eka}. We integrate the resulting
inequality
 with
respect to $\tau$ over the interval $[t_{1},T]$, which affects only
the last integral. Upon division by $T-t_{1}$, the last term is replaced by
$$\frac{2 j}{T - t_{1}} \int_{t_{1}}^{T}\!\!\int_{\Omega}
v_{j}\,dx\,dt.$$
We can combine this and the earlier estimate
\begin{equation*}
\int_{0}^{t_{1}}\!\!\int_{\Omega}|\nabla v_{j}|^{p}\,dx\,dt  \leq
j \int_{0}^{T}\!\!\int_{\Omega}|\nabla v_{1}|^{p}\,dx\,dt +
 \frac{j^{2-\varepsilon}}{T-t_{1}} \int_{t_{1}}^{T}\!\!\int_{\Omega}
v^{\varepsilon}\,dx\,dt,
\end{equation*}
taking $\varepsilon = 1$, so that we finally arrive at
\begin{align}
  \label{orderj}
\int_{0}^{t_{1}}\!\!\int_{\Omega}|\nabla v_{j}|^{p}\,dx\,dt +
\underset{0<t<t_{1}}{\textrm{ess\,sup}} \int_{\Omega} v_{j}(x,t)^{2}
\,dx \\
\leq
3 j \int_{0}^{T}\!\!\int_{\Omega}|\nabla v_{1}|^{p}\,dx\,dt +
\frac{3j}{T - t_{1}} 
\int_{t_{1}}^{T}\!\!\int_{\Omega} v\,dx\,dt. \nonumber
\end{align}
The majorant is now $O(j)$, which is of the right order,
as the following lemma shows with  its sharp exponents.

\begin{lemma}\label{jK}
If $v_j \in L^p(0,T;W^{1,p}_0(\Omega))$ and
\begin{equation*}
\int_{0}^{T}\!\!\int_{\Omega}|\nabla v_{j}|^{p}\,dx\,dt +
\underset{0<t<T}{\textrm{ess\,sup}} \int_{\Omega} v_{j}(x,t)^{2}
\,dx  \leq jK
\end{equation*}
when $j = 1,2,\ldots$ then $\nabla v$ exists and
\begin{align*}
v &\in L^{q}(\Omega_{T})\ \ \  \mbox{\ whenever\ } \ \ \ \ \ 0 < q <
  p-1+\frac{p}{n},\\
\nabla v &\in L^{q}(\Omega_{T})\ \ \  \mbox{\ whenever\ }\ \ \ 0 < q <
p-1+\frac{1}{n+1}.
\end{align*}
\end{lemma}

\emph{Proof:} The first part is a repetition of the proof of Lemma
\ref{Kjj}. Denote again
$$ E_{j} = \{(x,t)|\,j \leq v(x,t) \leq 2j \}, \ \ \ \
\kappa = 1 + \frac{2}{n}.$$ We have as before
\begin{align*}
& \ j^{\kappa p}|E_{j}|  \leq \int \!\! \int_{E_{j}} v^{\kappa
  p}_{2j}\,dx\,dt \leq \integ{ v^{\kappa p}_{2j}}\\  & \leq
C \integ{|\nabla v_{2j}|^{p}} \cdot \left(\underset{0<t<T}{\text{ess\,sup}}
\int_{\Omega}v_{2j}^{2}\,dx \right)^{\frac{p}{n}} \leq
C K ^{1+\frac{p}{n}} (2j)^{1+\frac{p}{n}}.
\end{align*}
It follows that
$$ |E_{j}| \leq \text{Constant}\times j^{1-p-\frac{p}{n}}.$$

We estimate the $L^{q}$-norm using the subdivision of
the domain. Thus
\begin{align*}
\integ{v^{q}} \leq \ & T\, |\Omega| + \sum_{j=1}^{\infty}\int \!
\int_{E_{2^{j-1}}} v^{q}\,dx\,dt \\ \leq\ & T\, |\Omega| +
  \sum_{j=1}^{\infty} 2^{jq}|E_{2^{j-1}}|\\
\leq\ & T\, |\Omega| + C \sum_{j=1}^{\infty} 2^{j(q+1-p-\frac{p}{n})},
\end{align*}
which  converges in the desired range for $q$. Thus the first part is
proved.

For \textsf{the summability of the gradien}t, we use the bound on the measure
of the level sets $E_{j}$ and also the growth assumed for the energy of the
truncated functions. Fix a large index $k$ and
write, using that $|\nabla v_{k}| \leq |\nabla v_{2^{j}}|$ on $E_{2^{j-1}}$:
\begin{align*}
\integ{|\nabla v_{k}|^{q}} \lesssim \ & \sum_{j=1}^{\infty}\int \!
\int_{E_{2^{j-1}}}| \nabla v_{k}|^{q}\,dx\,dt \\ \leq\ & 
  \sum_{j=1}^{\infty}\left(  \int \!
\int_{E_{2^{j-1}}}| \nabla v_{k}|^{p}\,dx\,dt
\right)^{\frac{q}{p}}|E_{2^{j-1}}|^{1- \frac{q}{p}} \\
\leq\ & \sum_{j=1}^{\infty} \left(  \integ{
| \nabla v_{2^{j}}|^{p}}
\right)^{\frac{q}{p}}2^{(j-1)(1-\frac{q}{p})(1-p- \frac{p}{n})}\\
\leq &\sum_{j=1}^{\infty}2^{(j-1)(1-
  \frac{q}{p})(1-p- \frac{p}{n})}(2^{j}K)^{\frac{q}{p}},
\end{align*}
where the geometric series converges provided that $q < p-1+1/(n+1).$
Strictly speaking, the ``first term'' $K^{q/p}(T|\Omega|)^{1-q/p}$
ought to be added to the sum, since the integral over the set
$\{0<v<1\}$ was missing.
Now we may let $k$ go to infinity. $\Box$

\bigskip

A combination of the results in this section (formula \ref{L1T}, equation (\ref{orderj}), and Lemma \ref{jK}) yields the following

\bigskip

\begin{lemma}\label{auxiliaryl} Suppose that $v \geq 0$ is a $p$-supercaloric function in $\Omega_T$ with initial values $v(x,0) = 0$ in $\Omega$. If every $ v_k \in L^p(0,T;W_0^{1,p}(\Omega))$, then
  \begin{align*}
v &\in L^{q}(\Omega_{T_1})\ \ \  \mbox{\ whenever\ } \ \ \ \ \ 0 < q <
  p-1+\frac{p}{n},\\
\nabla v &\in L^{q}(\Omega_{T_1})\ \ \  \mbox{\ whenever\ }\ \ \ 0 < q <
p-1+\frac{1}{n+1}
  \end{align*}
  when $T_1 < T.$ In particular, $v$ is of class $\mathfrak{B}.$
\end{lemma}

\bigskip

In principle, this lemma is Theorem \ref{Mainpar} in the special case with zero lateral boundary values.
 
\section{Proof of the Theorem}
\label{Proof}

For the proof of Theorem \ref{Mainpar} we start with a non-negative $p$-supercaloric function $v$ defined in $\Omega_T$.
A simple device  is used for the initial values:  fix a small $\delta > 0$ and redefine $v$ so that  $v(x,t) \equiv 0$ when $t\leq \delta.$ This function is $p$-supercaloric, since it obviously satisfies the comparison principle. This does not affect the statement of the theorem, since we can take $\delta$ as small as we please. The initial condition $v(x,0) = 0$ required in Lemma \ref{auxiliaryl} is now in order. 

Let $Q_{2l} \subset \subset \Omega$ be a cube with side length $4l$ and consider the concentric cube
$$Q_{l} = \{x\big\vert\,|x_i-x^0_i| < l,\, i = 1,2,\dots n\}$$
of side length $2l.$ The center is at $x^0.$ The main difficulty is that $v$ is not zero on the lateral boundary, neither does $v_j$ obey Lemma \ref{auxiliaryl}.  We aim at correcting $v$ outside $Q_{l} \times (0,T)$ so that also the new function  is $p$-supercaloric and, in addition,  satisfies the requirements of zero boundary values in Lemma \ref{auxiliaryl}. Thus we study the function
\begin{equation}
\label{noll}
w = 
\begin{cases}
v \qquad \text{in}\quad Q_l \times (0,T)\\
h\qquad\text{in}\quad (Q_{2l} \setminus Q_l) \times (0,T)
\end{cases}
\end{equation}
where the function $h$ is, in the outer region,   the weak solution to the boundary value problem
\begin{equation}
\label{ringh}
\begin{cases} h = 0 \qquad \text{on}\qquad \partial Q_{2l} \times (0,T)\\
h = v \qquad \text{on}\qquad \partial Q_{l} \times (0,T)\\
h = 0 \qquad \text{on}\qquad ( Q_{2l} \setminus  Q_{l}) \times \{0\}
\end{cases}
\end{equation}
An essential observation is that the solution $h$ does not always exist!  
This counts for the dichotomy described in the main Theorem \ref{Mainpar}. If it exists, the truncations $w_j$ satisfy the assumptions in Lemma \ref{auxiliaryl}, as we shall see.

For the construction we use the infimal convolutions
$$v^{\varepsilon}(x,t) = \inf_{(y,\tau)\in\Omega_T}\Bigl\{v(y)+\frac{1}{2\varepsilon}(|x-y|^2+|t-\tau|^2)\Bigr\}$$
described in Section \ref{Bounded Supersolutions}.
They  are Lipschitz continuous in $\overline{Q_{2l}}\times [0,T].$ They are  weak supersolutions when $\varepsilon$ is small enough according to Proposition \ref{infcon} and Theorem \ref{Thm25}. Then we define the solution $h^{\varepsilon}$ as in formula (\ref{ringh}) above, but with $v^{\varepsilon}$ in place of $v$. Then we construct 
\begin{equation*}
w^{\varepsilon} =
\begin{cases}
v^{\varepsilon} \qquad \text{in}\quad Q_l \times (0,T)\\
h^{\varepsilon} \qquad\text{in}\quad (Q_{2l} \setminus Q_l) \times (0,T)
\end{cases}
\end{equation*}
and $w^{\varepsilon}(x,0)= 0$ in $\Omega.$ Now $h^{\varepsilon} \leq v^{\varepsilon},$ and when $t \leq \delta$ we have $ 0 \leq h^{\varepsilon} \leq v^{\varepsilon} = 0$ so that  $h^{\varepsilon}(x,t) = 0 $ when $t \leq \delta.$ The function $w^{\varepsilon}$ satisfies the comparison principle and is therefore a $p$-supercaloric function. Here it is essential that $h^{\varepsilon} \leq v^{\varepsilon}$\,! The function  $w^{\varepsilon}$ is also (locally) bounded; thus 
we have arrived at the conclusion that  $w^{\varepsilon}$ is a weak supersolution in $Q_{2l} \times (0,T).$
Se Theorem \ref{Thm25}.

There are two possibilities, depending on whether the sequence $\{h^{\varepsilon}\}$ is bounded or not, when $\varepsilon \searrow 0$ through a sequence of values.

\textsf{Bounded case.} Assume that there does not exist any sequence of points such that
$$ \lim_{\varepsilon \to 0}h^{\varepsilon}(x_{\varepsilon},t_{\varepsilon})\, =\, \infty,\quad (x_{\varepsilon},t_{\varepsilon}) \to (x_0,t_0)$$
where $x_0 \in Q_{2l} \setminus \overline{Q_{l}}$ and $0  < t_0 < T$ (that is an \emph{interior} limit point). By Proposition \ref{increasing}
$$h =  \lim_{\varepsilon \to 0}h^{\varepsilon}$$
is a $p$-caloric function in its domain. The function\, $w = \lim w^{\varepsilon}$\, itself is $p$-supercaloric and agrees with formula (\ref{noll}). 

By Theorem \ref{Thm25} the truncated functions
$$w_j = w_j(x,t) = \min\{w(x,t),j\},\qquad j = 1,2,3,\dots,$$
are weak supersolutions in $Q_{2l} \times (0,T).$ We claim that
$w_j \in L^p(0,T';W^{1,p}_0(Q_{2l}))$ when $T' < T.$ This requires an estimation where we use
$$L = \sup \{h\}\quad\text{over}\quad (Q_{2l}\setminus Q_{5l/4})\times (0,T').$$
Let $\zeta = \zeta(x)$ be a smooth function such that $$0 \leq \zeta \leq 1,\quad  \zeta = 1\quad\text{in}\quad Q_{2l} \setminus Q_{3l/2},\quad \zeta = 0\quad\text{in}\quad Q_{5l/4}.$$ Using the test function 
$\zeta^ph$ when deriving the Caccioppoli estimate we get
\begin{align*}
&\int_{0}^{T'}\!\!\int_{Q_{2l}\setminus Q_{3l/2}}|\nabla w_j|^p \dif x \dif t\\ &\leq \,
\int_{0}^{T'}\!\!\int_{Q_{2l}\setminus Q_{3l/2}}|\nabla h|^p \dif x \dif t\, \leq \,
\int_{0}^{T'}\!\!\int_{Q_{2l}\setminus Q_{5l/4}}\zeta^p|\nabla h|^p \dif x \dif t \\ &\leq
C(p)\Bigl\{ \int_{0}^{T'}\!\!\int_{Q_{2l}\setminus Q_{l}}h^p|\nabla \zeta|^p \dif x \dif t +
\int_{Q_{2l}\setminus Q_{5l/4}}h(x,T')^2\dif x \Bigr\}\\ &\leq C(n,p)\bigl(L^pl^{n-p}T + L^2l^n\Bigr),
\end{align*}
where we first used 
$$|\nabla w_j| = |\nabla \min\{h,j\}| \leq |\nabla h|$$
in the outer region.
Thus we have an estimate over the outer region $Q_{2l}\setminus Q_{3l/2}.$ Concerning the inner region $Q_{3l/2},$ we first choose a test function $\eta = \eta(x,t)$
$$0 \leq \eta \leq 1,\quad \eta \equiv 1\quad
\text{in}\quad Q_{3l/2},\quad \eta =0\quad\text{in}\quad Q_{2l} \setminus Q_{9l/4}.$$
 Then the Caccioppoli estimate for the truncated functions
$$w_j = \min\{w,j\}, \qquad j = 1,2,3, \dots,$$
takes the form
\begin{gather*}
\int_{0}^{T'}\!\!\int_{Q_{3l/2}}|\nabla w_j|^p\dif x \dif t \leq  \int_{0}^{T'}\!\!\int_{Q_{2l}}\eta^p|\nabla w_j|^p\dif x \dif t \\
\leq Cj^p\int_{0}^{T'}\!\!\int_{Q_{2l}}|\nabla \eta|^p\dif x \dif t + Cj^p\int_{0}^{T'}\!\!\int_{Q_{2l}}|\eta_t|^p\dif x \dif t.
\end{gather*}

 Therefore we have obtained an estimate over the whole domain $Q_{2l} \times (0,T')$:
$$\int_{0}^{T'}\!\!\int_{Q_{2l}}|\nabla w_j|^p\dif x \dif t \,\,\leq \,\,Cj^p$$
and it follows that $w_j \in L(0,T';W^{1,p}_0(Q_{2l})).$ In particular, the crucial estimate
$$\int_{0}^{T'}\!\!\int_{Q_{2l}}|\nabla w_1|^p\dif x \dif t < \infty,$$
which was ``assumed'' in [KL], is now established.\footnote{The class $\mathfrak{M}$ passed unnoticed in [KL2].}

From  Lemma  \ref{auxiliaryl} we conclude that $v \in L^q(Q_{l}\times (0,T'))$ and $\nabla v \in L^{q'}(Q_{l}\times (0,T'))$ with the correct summability exponents. Either we can proceed like this for all interor cubes, or then the following case happens.

{\sf Unbounded case}\footnote{This case is described in [KP] and [KL3].}. If 
$$ \lim_{\varepsilon \to 0}h^{\varepsilon}(x_{\varepsilon},t_{\varepsilon}) = \infty,\quad (x_{\varepsilon},t_{\varepsilon}) \rightarrow (x_0,t_0)$$
for some $x_0 \in Q_{2l}\setminus \overline{Q_{l}},\,\, 0 < t_0 < T,$\, then
$$v(x,t) \geq h(x,t) \geq (t-t_0)^{-\frac{1}{p-2}}\mathfrak{U}(x),\quad \text{when}\quad t > t_0,$$
according to Proposition \ref{blowup}. Therefore
$$v(x,t_o+)\,=\,\infty$$
in $Q_{2l}\setminus \overline{Q_{l}}.$ But in this construction we can replace the outer cube $Q_{2l}$ with $\Omega,$ that is, a new $h$ is defined in $\Omega \setminus \overline{Q_l}.$ The proof is the same as above.  Then by comparison
$$v\, \geq \,  h^{\Omega}\, \geq \,  h^{Q_{2l}}$$
and so $v(x,t_0+) = \infty$ in the whole boundary zone $\Omega \setminus  \overline{Q_l}.$

It remains to include the inside, the cube $Q_{l}.$ This is easy. Reflect $h=h^{Q_{2l}}$
in the plane  $x_1 = x_1^0 + l,$
which contains one side of the small cube:
$$h^*(x_1,x_2.\dots,x_n) = h(2x^0_1+2l-x_1,x_2,\dots,x_n),$$
so that
$$\dfrac{x_1 + \bigl(2(x^0_1+l)-x_1\bigr)}{2} = x_1^0 + l$$
as it should. Recall that $x^0$ was the center of the cube. (The same can be done earlier for all the $h^{\varepsilon}.$) The reflected function $h^*$ is $p$-caloric. Clearly, $v\geq h^*$ by comparison. This forces $v(x,t_0+) = 0$ when
$x \in Q_{l},\, x_1 > x_1^0.$ A similar reflexion in the plane $x_1 = x_1^0 - l$ includes the other half  $x_1 < x_1^0.$
We have achieved that $v(x,t_0+) = \infty$ also in the inner cube $Q_l.$
This proves that 
$$v(x,t_o+)\,\equiv\,\infty\quad \text{in the whole}\quad \Omega.$$\qquad $\Box$

\bigskip 

From the proof we can extract that 
\begin{equation}
\label{extracted}
v(x,t) \geq \dfrac{\mathfrak{U}(x)}{(t-t_0)^{\frac{1}{p-2}}}\quad \text{in}\quad \Omega\times(t_0,T),
\end{equation}
where $\mathfrak{U}$ is from equation (\ref{apuyh}).


\section{Weak Supersolutions Are Semicontinuous}
\label{Semicontinuous}

Are the weak supersolutions $p$-supercaloric functions
 (= viscosity supersolutions)? The question is seemingly trivial, but there is a requirement.
To qualify they have to obey the comparison principle and to be
semicontinuos. The comparison principle is rather immediate. The
semicontinuity is a  delicate issue.
For a weak supersolution defined in the classical way with test
functions under the integral sign (Definition 16) the Sobolev
derivative is assumed to exist, but {\sf the semicontinuity}, which now is
 not assumed, {\sf has to be
established}. The proof requires parts of the classical regularity
theory\footnote{The preface of \textsc{Giuseppe Mingione'}s work [M] 
is worth
  reading as an enlightment.}. We will use a variant of the Moser iteration, for practical
reasons worked out for weak \emph{sub}solutions bounded from below. Our
proof of the theorem below is essentially the same as in [K], but
we avoid the use of infinitely stretched infinitesimal space-time
cylinders.

\begin{thm}\label{semico}
 Suppose that $v = v(x,t)$ is a weak supersolution of the
  Evolutionary $p$-Laplace equation. Then it is locally bounded from
  below and at almost every point $(x_{0},t_{0})$ it holds that
$$v(x_{0},t_{0}) =  \mathrm{ess}\!
\liminf_{(x,t) \rightarrow (x_{0},t_{0})}v(x,t).$$
In particular, $v$ is lower semicontinuous after a redefinition in a
set of measure zero.
\end{thm}

Functions like \, $ \mathrm{ess}\!
\liminf v(x,t)$ \, are lower semicontinuous, if they are bounded from
below. Thus the problem is the formula. The hardest part of the
proof is to establish  that the supremum norm of a non-negative
weak subsolution is  $1^{0}$) bounded (Lemma \ref{subsol34}) and $2^{0}$) bounded in terms
 of quantities
that can carry information from the Lebesgue points (Theorem \ref{cylinder}).
  With such
estimates the proof follows easily (at the end of this
section). Before entering into the semicontinuity proof we address
 the comparison principle.

\begin{prop}[Comparison Principle] Let $\Omega$ be bounded. Suppose that $v$ is a weak
  supersolution and $u$ a weak subsolution, $u,v \in
  L^{p}(0,T;W^{1,p}(\Omega))$, satisfying 
$$ \liminf v \geq \limsup u$$
on the parabolic boundary. Then $v \geq u$ almost everywhere in the  domain
$\Omega_{T}.$
\end{prop}

\emph{Proof:} This is well-known and we only give a formal proof. For
a non-negative test function $\varphi \in C^{\infty}_{0}(\Omega_{T})$
the equations
\begin{gather*}
\int_{0}^{T}\!\int_{\Omega}(-v\varphi_{t}+\langle|\nabla
v|^{p-2}\nabla v,\nabla \varphi\rangle) \,dx\,dt \geq 0\\
\int_{0}^{T}\!\int_{\Omega}(+u\varphi_{t}-\langle|\nabla
u|^{p-2}\nabla u,\nabla \varphi\rangle) \,dx\,dt \geq 0
\end{gather*}
can be added. Thus
\begin{equation*}
\int_{0}^{T}\!\int_{\Omega}\left((u-v)\varphi_{t}+\langle  |\nabla
v|^{p-2}\nabla v - |\nabla
u|^{p-2}\nabla u,\nabla \varphi\rangle \right) \,dx\,dt \geq 0.
\end{equation*}
These equations remain true if $v$ is replaced by $v+\varepsilon$,
where $\varepsilon$ is any constant. To complete the proof we choose
(formally) the test function to be
$$\varphi = (u-v-\varepsilon)_{+}\eta,$$
where $\eta = \eta(t)$ is a cut-off function; even $\eta(t) = T-t$
will do here. We arrive at
\begin{align*}
&\phantom{=a }\int_{0}^{T}\!\int_{u\geq v+\varepsilon}\eta(\langle  |\nabla
v|^{p-2}\nabla v - |\nabla
u|^{p-2}\nabla u,\nabla v - \nabla u \rangle) \,dx\,dt \\
&\leq \int_{0}^{T}\!\int_{\Omega}(u-v-\varepsilon)^{2}_{+}\eta'\,dx\,dt
+\frac{1}{2}\int_{0}^{T}\!\int_{\Omega}\eta \frac{\partial}{\partial
  t}(u-v-\varepsilon)^{2}_{+}\,dx\,dt\\
&= \frac{1}{2}\int_{0}^{T}\!\int_{\Omega}(u-v-\varepsilon)^{2}_{+}\eta'\,dx\,dt\\
&= -
\frac{1}{2}\int_{0}^{T}\!\int_{\Omega}(u-v-\varepsilon)^{2}_{+}\,dx\,dt
\leq 0.
\end{align*}
Since the first integral is non-negative by the vector
inequality (\ref{Vect}), the last integral is, in fact, zero. Hence the integrand 
$(u-v-\varepsilon)^{2}_{+} = 0$ almost everywhere. But this means that 
$$u \leq v + \varepsilon$$
almost everywhere. Since $\varepsilon > 0$ we have the desired
inequality $v\geq u$ a.e..  $\Box$

\medskip

We need some estimates for the semicontinuity proof and
 begin with the well-known Caccioppoli estimates, which are
extracted directly from the differential equation.

\begin{lemma}[Caccioppoli estimates] For a non-negative weak subsolution
  $u$ in $\Omega \times (t_{1},t_{2})$ we have the estimates
\begin{gather*}
\underset{t_{1}<t<t_{2}}{\esssup}\int_{\Omega}\zeta^{p}u^{\beta+1}\,dx
\leq
\int_{t_{1}}^{t_{2}}\!\int_{\Omega}u^{\beta+1}\Bigl|\frac{\partial}{\partial
    t}\zeta^{p}\Bigr|\,dx\,dt\\
+2p^{p-1}\beta^{2-p}\int_{t_{1}}^{t_{2}}\!\int_{\Omega}u^{p-1+\beta}
\left|\nabla\zeta\right|^{p}\,dx\,dt\\
\intertext{and}
\int_{t_{1}}^{t_{2}}\!\int_{\Omega}\left|\nabla(\zeta
  u^{\frac{p-1+\beta}{p}})\right|^{p}\,dx\,dt \leq C \beta^{p-2}
\int_{t_{1}}^{t_{2}}\!\int_{\Omega}u^{\beta+1}\Bigl|\frac{\partial}{\partial
    t}\zeta^{p}\Bigr|\,dx\,dt\\
+
C\int_{t_{1}}^{t_{2}}\!\int_{\Omega}u^{p-1+\beta}
\left|\nabla\zeta\right|^{p}\,dx\,dt,\\
\end{gather*}
where the exponent $\beta \geq 1$, $C = C(p)$, and $
\zeta \in C^{\infty}(\Omega \times [t_{1},t_{2})),$\\ $ \zeta(x,t_{1})
= 0,\,\zeta \geq 0.$
\end{lemma}

\medskip
\emph{Proof:} Use the test function $\varphi = u^{\beta}\zeta^{p}$
in the equation
\begin{gather*}
\int_{t_{1}}^{\tau}\!\int_{\Omega}\left(-u\varphi_{t} +
  \langle|\nabla u|^{p-2}\nabla u, \nabla \varphi\rangle\right)dx\,dt\\
+ \int_{\Omega}u(x,\tau)\varphi(x,\tau)\,dx \leq
\int_{\Omega}u(x,t_{1})\varphi(x,t_{1})\,dx = 0,
\end{gather*}
where $t_{1} < \tau \leq t_{2}$. (The intermediate $\tau$ is needed 
to match the supremum in the first estimate.) Strictly
speaking, the ``forbidden'' time derivative $u_{t}$ is required at the
intermediate steps. This can be handled through a regularization,
which we omit. Proceeding, integration by parts leads to
\begin{gather*}
\int_{t_{1}}^{\tau}\!\int_{\Omega}-u\varphi_{t}\,dx\,dt +
 \int_{\Omega}u(x,\tau)\varphi(x,\tau)\,dx\\
=
\frac{1}{\beta+1}\int_{\Omega}\zeta(x,\tau)^{p}u(x,\tau)^{\beta+1}\,dx
-\frac{1}{\beta+1}\int_{t_{1}}^{\tau}\!\int_{\Omega}u^{\beta+1}\Bigl|\frac{\partial}{\partial
    t}\zeta^{p}\Bigr|\,dx\,dt
\end{gather*}
valid for a.e. $\tau$. To treat the ``elliptic term'', we use 
$$\nabla \varphi = \beta \zeta^{p}u^{\beta-1}\nabla u + p
\zeta^{p-1}u^{\beta}\nabla \zeta$$
and obtain
\begin{align*}
\frac{1}{\beta+1}\int_{\Omega}\zeta(x,\tau)^{p}u(x,\tau)^{\beta+1}\,dx
+
\beta\int_{t_{1}}^{\tau}\!\int_{\Omega}\zeta^{p}u^{\beta-1}|\nabla
u|^{p}\,dx\,dt\\
\leq \frac{1}{\beta+1}\int_{t_{1}}^{\tau}\!\int_{\Omega}u^{\beta+1}\Bigl|\frac{\partial}{\partial
    t}\zeta^{p}\Bigr|\,dx\,dt + p\int_{t_{1}}^{\tau}\!\int_{\Omega}\zeta^{p-1}u^{\beta}|\nabla
u|^{p-1}|\nabla \zeta|\,dx\,dt.
\end{align*}
As much as possible of the last integral must be absorbed by the
double integral in the left-hand member. It is convenient to employ
Young's inequality
$$ab \leq \frac{a^{q}}{q} + \frac{b^{p}}{p}$$
to achieve the splitting
\begin{align*}
&\zeta^{p-1}u^{\beta}|\nabla
u|^{p-1}|\nabla \zeta|\\
&=
\overbrace{\left(\frac{\beta}{p}\right)^{\frac{p-1}{p}}\zeta^{p-1}u^{(\beta-1)\frac{p-1}{p}}|\nabla
u|^{p-1}}^{a} \times \overbrace{
\left(\frac{p}{\beta}\right)^{\frac{p-1}{p}}u^{\frac{p-1+\beta}{p}}|\nabla
\zeta|}^{b}\\
&\leq
\frac{p-1}{p}\left(\frac{\beta}{p}\right)\zeta^{p}u^{\beta-1}|\nabla
u|^{p}
+ \frac{1}{p}\left(\frac{p}{\beta}\right)^{p-1}u^{p-1+\beta}|\nabla
\zeta|^{p},
\end{align*}
which has to be multiplied by $p$ and integrated. Absorbing one
integral into the left-hand member, we arrive at the fundamental
estimate
\begin{gather*}
\frac{1}{\beta+1}\int_{\Omega}\zeta(x,\tau)^{p}u(x,\tau)^{\beta+1}\,dx
+\frac{\beta}{p}\int_{t_{1}}^{\tau}\!\int_{\Omega}\zeta^{p}u^{\beta-1}|\nabla
u|^{p}\,dx\,dt\\
\leq \frac{1}{\beta+1}\int_{t_{1}}^{\tau}\!\int_{\Omega}u^{\beta+1}|\frac{\partial}{\partial
    t}\zeta^{p}|\,dx\,dt + \left(\frac{p}{\beta}\right)^{p-1}\int_{t_{1}}^{\tau}\!\int_{\Omega}u^{p-1+\beta}|\nabla
\zeta|^{p}\,dx\,dt.
\end{gather*}

Since the integrands are positive it follows that
\begin{gather*}
\frac{1}{\beta+1}\int_{\Omega}\zeta(x,\tau)^{p}u(x,\tau)^{\beta+1}\,dx\\
\leq \frac{1}{\beta+1}\int_{t_{1}}^{t_{2}}\!\int_{\Omega}u^{\beta+1}\Bigl|\frac{\partial}{\partial
    t}\zeta^{p}\Bigr|\,dx\,dt + \left(\frac{p}{\beta}\right)^{p-1}\int_{t_{1}}^{t_{2}}\!\int_{\Omega}u^{p-1+\beta}|\nabla
\zeta|^{p}\,dx\,dt,
\end{gather*}
where the majorant now is free from $\tau$.
Taking the supremum over $\tau$ we obtain the first Caccioppoli
inequality.

 To derive the second Caccioppoli inequality, we start from
\begin{gather*}
\frac{\beta}{p}\int_{t_{1}}^{t_{2}}\!\int_{\Omega}\zeta^{p}u^{\beta-1}|\nabla
u|^{p}\,dx\,dt\\
\leq \frac{1}{\beta+1}\int_{t_{1}}^{t_{2}}\!\int_{\Omega}u^{\beta+1}\Bigl|\frac{\partial}{\partial
    t}\zeta^{p}\Bigr|\,dx\,dt + \left(\frac{p}{\beta}\right)^{p-1}\int_{t_{1}}^{t_{2}}\!\int_{\Omega}u^{p-1+\beta}|\nabla
\zeta|^{p}\,dx\,dt
\end{gather*}
and notice that
$$\zeta^{p}u^{\beta-1}|\nabla u|^{p} =
 \biggl(\frac{p}{p-1+\beta}\biggr)^{p}|\zeta \nabla u^{\frac{p-1+\beta}{p}}|^{p}
.$$
Then the triangle inequality
$$ |\nabla(\zeta u^{\frac{p-1+\beta}{p}})| \leq |\zeta \nabla
  u^{\frac{p-1+\beta}{p}}| + |u^{\frac{p-1+\beta}{p}} \nabla
      \zeta|$$
and a simple calculation yield  the desired result.  $\Box$

\medskip

In the following version of Sobolev's inequality  the exponents are
adjusted to our need. For a proof, see [dB, Chapter 1].

\begin{prop} [Sobolev] For $\zeta \in C^{\infty}(\Omega_{T})$
  vanishing on the lateral boundary $\partial \Omega \times [0,T]$ we have
\begin{gather*}
\int_{0}^{T}\!\int_{\Omega}
\zeta^{p\gamma}|u|^{p-2+(\beta+1)\gamma}\,dx\,dt\\
\leq S
\int_{0}^{T}\!\int_{\Omega}|\nabla(\zeta|u|^{\frac{p-1+\beta}{p}})|^{p}\,dx\,dt
\left\{\underset{0<t<T}{\esssup}\int_{\Omega}\zeta^{p}|u|^{\beta+1}\,dx
\right\}^{\frac{p}{n}},
\end{gather*}
where $\gamma = 1 + \frac{p}{n}.$
\end{prop}

Now we can control the right-hand member in the Sobolev inequality by
the quantities in the Caccioppoli estimates for the weak subsolution:
\begin{gather*}
\left(\int_{t_{1}}^{t_{2}}\!\int_{\Omega}
\zeta^{p\gamma}u^{p-2+(\beta+1)\gamma}\,dx\,dt\right)^{\frac{1}{\gamma}}\\
\leq C \beta^{\frac{(2-p)p}{n+p}}\left(\beta^{p-2}
\int_{t_{1}}^{t_{2}}\!\int_{\Omega}u^{\beta+1}\Bigl|\frac{\partial}{\partial
    t}\zeta^{p}\Bigr|\,dx\,dt
+
\int_{t_{1}}^{t_{2}}\!\int_{\Omega}u^{p-1+\beta}
\left|\nabla\zeta\right|^{p}\,dx\,dt\right).
\end{gather*}
We select the test function $\zeta$ so that it is equal to $1$ in the
cylinder
$B_{R-\Delta R} \times (T+\Delta T,t_{2})$, $\zeta(x,T) = 0$, and so that
$\zeta(x,t) = 0$ when x is outside $B_{R}.$ Then we can write

\begin{gather*}
\left(\int_{T+\Delta T}^{t_{2}}\!\int_{B_{R-\Delta R}}
u^{p-2+(\beta+1)\gamma}\,dx\,dt\right)^{\frac{1}{\gamma}}\\
\leq C \beta^{\frac{(2-p)p}{n+p}}\left(\frac{\beta^{p-2}}{\Delta T}
\int_{T}^{t_{2}}\!\int_{B_{R}}u^{\beta+1}\,dx\,dt
+
 \frac{1}{\left(\Delta R \right)^{p}}\int_{T}^{t_{2}}\!\int_{B_{R}}
  u^{p-1+\beta}
\,dx\,dt\right),
\end{gather*}
where C is a new constant.  Recall that $\gamma > 1.$ 
This is a \emph{reverse H\"{o}lder
  inequality}, which is most transparent for $p = 2$. It will be important to
keep $\Delta T = (\Delta R)^{p}$. This is the basic inequality for the
celebrated \emph{Moser iteration}, which we will employ. The power of
$u$ increases to $p-2 +(\beta+1)\gamma$, but the integral is taken
over a smaller cylinder. In order to iterate over a chain of
shrinking cylinders $U_{k} = B(x_{0},R_{k}) \times (T_{k},t_{2})$ , starting with 
$$ U_{0} = B(x_{0},2R) \times (\frac{T}{2},t_{2})$$ and
ending up with an estimate over the cylinder 
$$ U_{\infty} = B(x_{0},R) \times (T,t_{2}),$$ 
we introduce the quantities
\begin{align*}
R_{k} &= R + \frac{R}{2^{k}},  &R_{k} - R_{k+1} &= \frac{R}{2^{k+1}}\\
T_{k} &= T - \frac{T}{2^{kp+1}}, &T_{k+1} - T_{k} &=
  \frac{T}{2^{(k+1)p}}s,\\
\omega &= \frac{R^{p}}{Ts} =  \frac{(\Delta R_{k})^{p}}{\Delta T_{k}} &s &=
\frac{2^{p-1}-1}{2}.
\end{align*}
We remark that \emph{$\omega$ is independent of the index $k$}. Further, we
write
$\alpha = \beta +1$, so that $\alpha \geq 2$. Thus
\begin{gather}
\label{maja}
\Biggl(\underset{U_{k+1}}{\int\!\int}
u^{p-2+\alpha\gamma}\,dx\,dt\Biggr)^{\frac{1}{\gamma}}\\
\leq C\,\frac{2^{(k+1)p}\beta^{\frac{(2-p)p}{n+p}}}{R^{p}}\Biggl(\beta^{p-2}\omega
\underset{U_{k}}{\int\!\int}u^{\alpha}\,dx\,dt
+
 \underset{U_{k}}{\int\!\int}
  u^{p-2+\alpha}
\,dx\,dt\Biggr).\nonumber
\end{gather}

It is inconvenient to deal with \emph{two}  different integrals in the
majorant. For simplicity we will perform two iteration procedures,
 depending on which
integral is dominating. For the \textsf{first procedure} we assume that
$$\omega \leq u^{p-2}.$$
Then we have the simpler expression
\begin{gather*}
\Biggl(\underset{U_{k+1}}{\int\!\int}
u^{p-2+\alpha\gamma}\,dx\,dt\Biggr)^{\frac{1}{\gamma}}
\leq C_{1}\,\frac{2^{kp}\alpha^{\frac{(p-2)}{\gamma}}}{R^{p}}
 \underset{U_{k}}{\int\!\int}
  u^{p-2+\alpha}
\,dx\,dt.
\end{gather*}
We start the iteration with $\alpha = 2$ and $k = 0$. Thus
\begin{gather*}
\Biggl(\underset{U_{1}}{\int\!\int}
u^{p-2+2\gamma}\,dx\,dt\Biggr)^{\frac{1}{\gamma}}
\leq C_{1}\,\frac{2^{0p}2^{\frac{(p-2)}{\gamma}}}{R^{p}}
 \underset{U_{0}}{\int\!\int}
  u^{p}
\,dx\,dt.
\end{gather*}
Then take $\alpha = 2\gamma$ and $k=1$ so that
\begin{gather*}
\Biggl(\underset{U_{2}}{\int\!\int}
u^{p-2+2\gamma^{2}}\,dx\,dt\Biggr)^{\frac{1}{\gamma^{2}}}
\leq
\Biggl(C_{1}\,\frac{2^{1p}(2\gamma)^{\frac{(p-2)}{\gamma}}}{R^{p}} 
 \underset{U_{1}}{\int\!\int}
  u^{p-2+2\gamma}
\,dx\,dt  \Biggr)^{\frac{1}{\gamma}} \\ 
\leq
\Biggl(C_{1}\,\frac{2^{1p}(2\gamma)^{\frac{(p-2)}{\gamma}}}{R^{p}}\Biggr)
^{\frac{1}{\gamma}}\times C_{1}\,\frac{2^{0p}2^{\frac{(p-2)}{\gamma}}}{R^{p}}
 \underset{U_{0}}{\int\!\int}
  u^{p}
\,dx\,dt.
\end{gather*}
The result of the next step is
\begin{gather*}
\Biggl(\underset{U_{3}}{\int\!\int}
u^{p-2+2\gamma^{3}}\,dx\,dt\Biggr)^{\frac{1}{\gamma^{3}}}
\\
\leq
\Biggl(\frac{C_{1}2^{\frac{p-2}{\gamma}}}{R^{p}}\Biggr)^{1+\frac{1}{\gamma}
+\frac{1}{\gamma^{2}}}2^{p(\frac{1}{\gamma}+\frac{2}{\gamma^{2}})}\gamma^{(p-2)
(\frac{1}{\gamma^{2}}+\frac{2}{\gamma^{3}})}
 \underset{U_{0}}{\int\!\int}
  u^{p}
\,dx\,dt.
\end{gather*}
Continuing the chain and noticing that the geometric series
$$1+\frac{1}{\gamma}
+\frac{1}{\gamma^{2}}+\frac{1}{\gamma^{3}}+ \cdots = 1 + \frac{n}{p}$$
and the series $\sum k\gamma^{-k}$ appearing in the exponents converge,
 since $\gamma > 1$, we arrive at
\begin{gather*}
\Biggl(\underset{U_{k+1}}{\int\!\int}
u^{p-2+2\gamma^{k+1}}\,dx\,dt\Biggr)^{\frac{1}{\gamma^{k+1}}}
\leq K R^{-p(1+\frac{1}{\gamma}
+\frac{1}{\gamma^{2}}+\cdots+\frac{1}{\gamma^{k}})} \underset{U_{0}}{\int\!\int}
  u^{p}
\,dx\,dt.
\end{gather*}
Here $K$ is a numerical constant. As $k\rightarrow \infty$, we obtain
the final estimate
\begin{gather*}
\underset{B_{R}\times(T,t_{2})}{\esssup} (u^{2}) \leq  \frac{K}{R^{n+p}}
 \int_{\frac{T}{2}}^{t_{2}}\!\int_{B_{2R}} u^{p}\,dx\,dt
= \frac{K}{\omega s\,TR^{n}} \int_{\frac{T}{2}}^{t_{2}}\! \int_{B_{2R}} u^{p}\,dx\,dt,
\end{gather*}
where the square came from the factor $2$ in $2\gamma^{k+1}$. The sum
of the geometric series determined the power of $R$.

Finally, if the assumption $\omega \leq u^{p-2}$ is relaxed to $u \geq
0$, we can apply the previous estimate to the function
$$u(x,t) + \omega ^{\frac{1}{p-2}}= u(x,t) +
\Bigl(\frac{R^{p}}{Ts}\Bigr)^{\frac{1}{p-2}}.$$
 A simple calculation
gives us the bound in the next lemma.

\begin{lemma}
\label{subsol34}
 Suppose that $u \geq 0$ is a weak subsolution in the
  cylinder\\ $B_{2R}\times (\frac{T}{2},t_{2})$. Then  
\begin{gather*}
\underset{B_{R}\times(T,t_{2})}{\esssup}\{u^{2}\} \leq C
\biggl\{\Bigl(\frac{R^{p}}{T}\Bigr)^{\frac{2}{p-2}} +
  \frac{T}{R^{p}}\Bigl(\frac{1}{TR^{n}}\int_{\frac{T}{2}}^{t_{2}}\!
 \int_{B_{2R}} u^{p}\,dx\,dt\Bigr)\biggr\},
\end{gather*}
where $C = C(n,p).$ 
\end{lemma}

\medskip \noindent
We can extract the following piece of information.

\begin{cor} A weak supersolution that is bounded from above, is
  locally bounded from below.
\end{cor}

\emph{Proof:} Use $u(x,t) = L - v(x,t).$\, $\Box$

\medskip

The estimate in the lemma suffers from the defect that it is not sharp
when $ u \approx 0$ because of the presence of the constant term. Our remedy
is a \textsf{second iteration procedure}, this time under the assumption that
$$ 0 \leq u \leq j,$$
where we take $j$ so large that also
$$j^{p-2} \geq \omega.$$  Read $j^{p-2}$ as $\max\{\omega,j^{p-2}\}$. The
previous lemma shows that $j$ is finite, but the point now is that $u$
is not bounded away from zero. Then
the
 first integral in the majorant of
(\ref{maja}) is dominating and we can begin with the bound
\begin{gather*}
\Biggl(\underset{U_{k+1}}{\int\!\int}
u^{p-2+\alpha\gamma}\,dx\,dt\Biggr)^{\frac{1}{\gamma}}
\leq C\,j^{p-2}\frac{2^{kp}\alpha^{\frac{(p-2)}{\gamma}}}{R^{p}}
 \underset{U_{k}}{\int\!\int}
  u^{\alpha}
\,dx\,dt.
\end{gather*}
We start the iteration with $\alpha = p$ and $k = 0$. Thus
\begin{gather*}
\Biggl(\underset{U_{1}}{\int\!\int}
u^{p-2+p\gamma}\,dx\,dt\Biggr)^{\frac{1}{\gamma}}
\leq Cj^{p-2}\,\frac{2^{0p}p^{\frac{(p-2)}{\gamma}}}{R^{p}}
 \underset{U_{0}}{\int\!\int}
  u^{p}
\,dx\,dt.
\end{gather*}
Then take $\alpha = p-2 +p\gamma$, which is $<n\gamma^{2}$, and $k=1$ so that
\begin{gather*}
\Biggl(\underset{U_{2}}{\int\!\int}
u^{(p-2)(1+\gamma)+p\gamma^{2}}\,dx\,dt\Biggr)^{\frac{1}{\gamma^{2}}}
\leq
\Biggl(Cj^{p-2}\,\frac{2^{1p}(n\gamma^{2})^{\frac{(p-2)}{\gamma}}}{R^{p}} 
 \underset{U_{1}}{\int\!\int}
  u^{p-2+p \gamma}
\,dx\,dt  \Biggr)^{\frac{1}{\gamma}} \\ 
\leq
\biggl(Cj^{p-2}\,\frac{2^{1p}( n\gamma^{2} )^{\frac{(p-2)}{\gamma}}}{R^{p}}\biggr)
^{\frac{1}{\gamma}}\times Cj^{p-2}\,
\frac{2^{0p}(n\gamma)^{\frac{(p-2)}{\gamma}}}{R^{p}}
 \underset{U_{0}}{\int\!\int}
  u^{p}
\,dx\,dt.
\end{gather*}
At the next step $\alpha = (p-2)(1+\gamma)+ p\gamma^{2} <
n\gamma^{3}$ and $k=2$. The result is
\begin{gather*}
\Biggl(\underset{U_{3}}{\int\!\int}
u^{(p-2)(1+\gamma+\gamma^{2})+p\gamma^{3}}\,dx\,dt\Biggr)^{\frac{1}{\gamma^{3}}}
\\
\leq
\biggl(\frac{Cj^{p-2}n^{\frac{p-2}{\gamma}}}{R^{p}}\biggr)^{1+\frac{1}{\gamma}
+\frac{1}{\gamma^{2}}}2^{p(\frac{1}{\gamma}+\frac{2}{\gamma^{2}})}\gamma^{(p-2)
( \frac{1}{\gamma}+ \frac{2}{\gamma^{2}}+\frac{3}{\gamma^{3}})}
 \underset{U_{0}}{\int\!\int}
  u^{p}
\,dx\,dt.
\end{gather*}
Continuing like this we end up with an estimate integrated over $U_{k+1}$
with the power $\alpha_{k+1} = p-2 +\alpha_{k}\gamma$, where
\begin{align*}
\alpha_{k} &= (p-2)(1+\gamma
+\gamma^{2}+\cdots+\gamma^{k-1})+p \gamma^{k}\\ & = 
\frac{n(p-2)}{p}(\gamma^{k}-1)+p \gamma^{k} \approx
\Bigl(n+p-\frac{2n}{p}\Bigr)\gamma^{k} 
\end{align*}
and $\alpha_{k} < n \gamma^{k+1}$. As $k \rightarrow \infty$ we find that
\begin{gather*}
\underset{B_{R}\times(T,t_{2})}{\esssup}\{u^{n+p-\frac{2n}{p}} \}
 \leq C
\frac{j^{\frac{(p-2)(n+p)}{p}}}{R^{n+p}}\int_{\frac{T}{2}}^{t_{2}}\!
 \int_{B_{2R}} u^{p}\,dx\,dt.
\end{gather*}
We can summarize the result.

\begin{thm}\label{cylinder} A weak subsolution $u$ that is non-negative in the
  cylinder $U = B(x_{0},2R)\times (t_{0}-3T/2,t_{0}+T)$ has the bound
\begin{equation}
\label{koko}
\underset{B_{R}\times(t_{0}-T,t_{0}+T )}{\esssup}\{u^{n+p-\frac{2n}{p}} \}
 \leq K
\frac{\left(\frac{R^{p}}{T} + \|u\|_{\infty}^{p-2}\right)
  ^{1+\frac{n}{p}}}{TR^{n}} \int_{ t_{0}-\frac{3T}{2}   }^{ t_{0}+T   }\!
 \int_{B_{2R}} u^{p}\,dx\,dt,
\end{equation}
where  $0 \leq u \leq \|u\|_{\infty}$\, in $U$.
\end{thm}

We need the fact that \textsf{the positive part $(u)_{+}$ of a weak
subsolution is again a weak subsolution}. Here the proof has to avoid the
comparison principle, which is not yet available. It reduces to the
following lemma.

\begin{lemma}\label{minwe}
If $v$ is a weak supersolution, so is $v_{L} = \min\{v,L\}.$
\end{lemma}

\medskip
\emph{Proof:} Formally, the test function\footnote{This is from Lemma
  2.109 on page 122 of J. Maly \& W. Ziemer:''Fine Regularity of
  Solutions of Elliptic Partial Differential Equations'',
  Math. Surveys Monogr. 51, AMS, Providence 1998.}
$$\varphi = \min\{k(L-v)_{+},1\}\zeta = \chi_{k} \zeta$$
inserted into
$$\int_{0}^{T}\!\int_{\Omega}\left(-v\varphi_{t}+ \langle|\nabla
  v|^{p-2}\nabla v,\nabla \varphi \rangle \right)dx\,dt \geq 0$$
implies the desired inequality
$$\int_{0}^{T}\!\int_{\Omega}\left(-v_{L}\zeta_{t}+ \langle|\nabla
  v_{L}|^{p-2}\nabla
  v_{L},\nabla \zeta \rangle \right)dx\,dt \geq 0$$
at the limit $k=\infty$. As usual, $\zeta \in
C^{\infty}_{0}(\Omega_{T}), \quad \zeta \geq 0$. The explanation is that $\lim
\chi_{k} = $ the characteristic function of the set $\{v<L\}$. Under
the assumption that the ``forbidden'' time derivative $u_{t}$ is
available at the intermediate steps we have
\begin{align*}
&\int_{0}^{T}\!\int_{\Omega}\chi_{k}\!\left(-v\zeta_{t}+ \langle|\nabla
  v|^{p-2}\nabla v,\nabla \zeta \rangle \right)dx\,dt \\ &\geq
k \underset{L-\frac{1}{k}<v<L}{\int\!\int} \zeta |\nabla v|^{p}\,dx\,dt + 
\int_{0}^{T}\!\int_{\Omega} v \zeta \frac{\partial}{\partial
  t}\chi_{k}\,dx\,dt\\
&\geq \int_{0}^{T}\!\int_{\Omega} v \zeta \frac{\partial}{\partial
  t}\chi_{k}\,dx\,dt = -\frac{1}{2k}\int_{0}^{T}\!\int_{\Omega}\zeta \frac{\partial}{\partial
  t}\left(\chi_{k}\right)^{2}\,dx\,dt\\
&= +\frac{1}{2k}\int_{0}^{T}\!\int_{\Omega}
\left(\chi_{k}\right)^{2}\zeta_{t}\,dx\,dt \,\longrightarrow\, 0.
\end{align*}
The formula $\partial \chi_{k}/\partial t = -vk$ or $=0$ was used. 
The result follows.

Finally, to handle the problem with the time derivative, one has first
to regularize the equation and then to use the test function
$$\varphi^{\varepsilon} = \min\{k(L-v^{\varepsilon})_{+},1\}\zeta =
\chi_{k} \zeta,$$
where $v^{\varepsilon}$ is the convolution in (\ref{convo}).
The term
$$ \int_{0}^{T}\!\int_{\Omega}-v^{\varepsilon}\frac{\partial
  \varphi^{\varepsilon}}{\partial t}\!dx\,dt$$
can be written so that the derivative $\partial
v^{\varepsilon}/\partial t$ disappears. Then one may safely let
$\varepsilon \rightarrow 0$. The result follows as before. $\Box$

\medskip

\emph{Proof of Theorem \ref{semico}:} Let $(x_{0},t_{0})$ be a Lebesgue point
for the weak supersolution $v$. Then
$$\lim_{TR^{n}\rightarrow
  0}\,\frac{1}{TR^{n}}\int_{t_{0}-2T}^{t_{0}+T}\!\int_{B_{2R}}|v(x_{0},t_{0})-v(x,t)|^{p}\,dx\,dt \, = \, 0.$$
\emph{A fortiori}
\begin{equation}
\lim_{TR^{n}\rightarrow
  0}\,\frac{1}{TR^{n}}\int_{t_{0}-2T}^{t_{0}+T}\!\int_{B_{2R}}(v(x_{0},t_{0})-v(x,t))_{+}^{p}\,dx\,dt \, = \, 0.
\end{equation}\label{averi}
We claim that
\begin{equation}\label{ess4}
v(x_{0},t_{0}) \leq \mathrm{ess}\!
\liminf_{(x,t) \rightarrow (x_{0},t_{0})}v(x,t).
\end{equation}
It is sufficient to establish that 
$$\mathrm{ess}\!
\limsup_{(x,t) \rightarrow (x_{0},t_{0})}( v(x_{0},t_{0})- v(x,t))_{+}
= 0,$$
since those points where $v(x,t) \geq v(x_{0},t_{0})$ can do no harm to
inequality (\ref{ess4}).

To this end, notice that the function $v(x_{0},t_{0})- v(x,t)$ is a
weak subsolution and so is its positive part, the function
$$u(x,t) = (v(x_{0},t_{0})- v(x,t))_{+}$$
by Lemma \ref{minwe}.
It is locally bounded according to Lemma \ref{subsol34}. Thus the essliminf is $>
- \infty$ in (\ref{ess4}). Use Theorem \ref{cylinder} and let
$TR^{n} \rightarrow 0$, keeping $R^{p}/T \leq$ Constant. In virtue of
(\ref{koko}) it follows that
$$ \mathrm{ess}\!
\limsup_{(x,t) \rightarrow
  (x_{0},t_{0})}\left\{u(x,t)^{n+p-\frac{2n}{p}}\right\} = 0 $$
and the exponent can be erased.
This proves the claim (\ref{ess4}) at the given Lebesgue point.

Furthermore, the Lebesgue points have the property that
\begin{align*}
&\phantom{ab}v(x_{0},t_{0})  \leq \mathrm{ess}\!
\liminf_{(x,t) \rightarrow (x_{0},t_{0})}v(x,t)\\&\leq \lim_{TR^{n}\rightarrow
  0}\,\frac{1}{3TR^{n}|B_{2}|}\int_{t_{0}-2T}^{t_{0}+T}\!\int_{B_{2R}}v(x,t)\,dx\,dt \, =
 v(x_{0},t_{0}).
\end{align*}
Since almost every point is a Lebesgue point, we have established that
$$v(x_{0},t_{0}) = \mathrm{ess}\!
\liminf_{(x,t) \rightarrow (x_{0},t_{0})}v(x,t)$$
almost everywhere. The right-hand member is a semicontinuous
function.\, $\Box$

\section{The Equation With Measure Data}
\label{Measuredata}

There is a close connexion between supersolutions and equations where
the right-hand side is a Radon measure. The Barenblatt solution has
the Dirac measure (multiplied by a suitable constant) as the
right-hand side, and hence it is, indeed, a \emph{solution} to an equation.
The equation 
 $$\frac{\partial v}{\partial t} - \nabla \cdot (|\nabla v|^{p-2}\nabla
v) = \mu$$
with a Radon measure $\mu$ has been much studied. For example, in [BD]
a summability result is given for the spatial gradient $\nabla v$ of
the solution. There the starting point was the given measure and the
above equation. However, we can do the opposite and \textsf{produce the
measure}. Indeed, every  $p$-supercaloric
function belonging to $L^{p-2}_{loc}(\Omega_T)$ induces a Radon measure $\mu \geq 0$. This follows from our
summability theorem, combined with the Riesz Representation Theorem for
linear functionals. However, if it so happens that $v$ belongs to class $\mathfrak{M}$, then
for some time $t_0,$
$$v(x,t)\geq (t-t_0)^{-\frac{1}{p-2}}\mathfrak{U}(x,t)$$
and it cannot induce any sigma finite measure, let alone a Radon measure.

\begin{thm}\label{Radon} Let $v$ be a $p$-supercaloric function in $\Omega \times
  (0,T)$. If $v$ is of class $\mathfrak{B}$  there exists a non-negative Radon measure $\mu$ such
  that 
$$\int_{0}^{T}\!\int_{\Omega}\Bigl(-v\frac{\partial \varphi}{\partial
    t}+ \langle |\nabla v|^{p-2}\nabla v,\nabla \varphi \rangle
\Bigr) dx\,dt = \int_{\Omega \times (0,T)} \varphi\, d \mu$$
for all $\varphi \in C_{0}^{\infty}(\Omega \times (0,T)).$
\end{thm}

\emph{Proof:} We already know that $v, \nabla v \in
 L_{loc}^{p-1}(\Omega \times (0,T)).$ In order to use Riesz's
 Representation Theorem we define the linear functional
\begin{gather*}
\Lambda_{v}:\,C_{0}^{\infty}(\Omega \times (0,T)) \longrightarrow
\mathbb{R},\\
\Lambda_{v}(\varphi) = \int_{0}^{T}\!\int_{\Omega}\Bigl(-v\frac{\partial
    \varphi}
{\partial
    t}+ \langle |\nabla v|^{p-2}\nabla v,\nabla \varphi \rangle
\Bigr) dx\,dt.
\end{gather*}
Now $\Lambda_{v}(\varphi) \geq 0$ for $\varphi \geq 0$ according to
Theorem \ref{Mainpar}. Thus the functional is positive and the existence of the
Radon measure follows from Riesz's theorem, cf. [EG, Section 1.8].\qquad  $\Box$

\medskip
Some further results can be found in [KLP]. The elliptic case has been thoroughly treated in [KL]. See also [Kuusi-Mingione].

\section{Pointwise Behaviour}

The viscosity supersolutions are defined at each point, not only
almost everywhere. Actually, the results in this section imply that two
\textsf{viscosity supersolutions that coincide almost everywhere do so at
each point}. 

\subsection{The Stationary Equation}
We begin with the stationary case. At each point a $p$-superharmonic
function $v$
satisfies
$$ v(x) \leq \liminf_{y \rightarrow x} v(y) \leq \mathrm{ess}
\liminf_{y \rightarrow x} v(y)$$
by lower semicontinuity. \emph{Essential limes inferior} means that sets of
Lebesgue measure zero be neglected in the calculation of the lower
limit. The reverse inequalities also hold. To see this, we start by a
lemma, which requires a pedantic formulation.

\begin{lemma}
Suppose that $v$ is $p$-superharmonic in the domain $\Omega$. If $v(x)
\leq \lambda$ at each point $x$ in $\Omega$ and if $v(x) = \lambda$ at almost
every point  $x$ in $\Omega$, then $v(x) = \lambda$ at each point  $x$ in
$\Omega$.
\end{lemma}

\emph{Proof:} The proof is trivial for continuous functions and the
idea is that $v$ is everywhere equal to a $p$-harmonic function,
which, of course, must coincide with the constant $\lambda$. We
approximate $v$ by the infimal convolutions $v_{\varepsilon}$. We can
assume that the  function $v$ is bounded also from below in  a given ball
$B_{2r}$, strictly
interior in $\Omega$. We may even take $0 \leq v \leq \lambda$ by
adding a constant. We
approximate $v$ by the infimal convolutions $v_{\varepsilon}$.
 Replace $v_{\varepsilon}$ in  $B_{r}$ by
the $p$-harmonic function $h_{\varepsilon}$ having boundary values 
$v_{\varepsilon}$. Thus we have the function
$$ w _{\varepsilon} = \left\{ \begin{array}{ll}
h_{\varepsilon} \ \mbox{in $B_{r}$}\\
v_{\varepsilon} \ \mbox{in $B_{2r} \backslash B_{r}$}
\end{array}
\right.$$
As we have seen before, also $w_{\varepsilon}$ is
$p$-superharmonic. By comparison
 $$ w _{\varepsilon} \leq
v_{\varepsilon} \leq v$$ pointwise in $B_{2r}$. As $\varepsilon$
approaches zero via a
decreasing sequence, say $1,1/2,1/3,\cdots$, the $h_{\varepsilon}$'s
converge to a $p$-harmonic function $h$, which is automatically
continuous because the family is uniformly equicontinuous so that
Ascoli's theorem applies. The equicontinuity is included in the
H\"{o}lder estimate (\ref{ellholder}), because $0 \leq h_{\varepsilon} \leq
\lambda$. Thus
$$ h \leq v \leq \lambda$$
at \emph{each} point in $B_{r}$. Since $\lambda - v_{\varepsilon} \geq
\lambda - v \geq 0$, the Caccioppoli estimate 
\begin{gather*}
\int_{B_{r}}\!|\nabla h_{\varepsilon}|^{p}\,dx \leq  
  \int_{B_{r}}\!|\nabla v_{\varepsilon}|^{p}\,dx   \\
\leq p^{p}\int_{B_{2r}}\!(\lambda -  v_{\varepsilon})^{p}|\nabla
\zeta|^{p}\,dx
 \leq Cr^{-p}\int_{B_{2r}}\!(\lambda -  v_{\varepsilon})^{p}\,dx 
\end{gather*}
is valid. The weak lower semicontinuity of the integral implies that
$$\int_{B_{r}}\!|\nabla h|^{p}\,dx \leq \lim_{\varepsilon \rightarrow
  0}
 \int_{B_{r}}\!|\nabla
h_{\varepsilon}|^{p}\,dx \leq Cr^{-p}\int_{B_{2r}}\!(\lambda -
v)^{p}\,dx = 0.$$
  The
conclusion is that $h$ is constant almost everywhere, and hence
everywhere by continuity. The constant must be $\lambda$, because it
has boundary values $\lambda$ in Sobolev's sense. We have proved that
also $v(x) = \lambda$ at \emph{each} point in the ball $B_{r}$.
 The result follows. $\Box$

\begin{lemma} If $v$ is $p$-superharmonic in $\Omega$ and if $v(x) >
  \lambda$ for a.e. $x$ in  $\Omega$, then $v(x) \geq
  \lambda$ for every $x$ in $\Omega$.
\end{lemma}

\emph{Proof:} If $\lambda = - \infty$, there is nothing to
prove. Applying the previous lemma to the $p$-superharmonic function
defined by
$$ \min\{v(x),\lambda\}$$
we obtain the result in the case $\lambda > - \infty$. $\Box$

\begin{thm} 
At each point a $p$-superharmonic
function $v$
satisfies
$$ v(x) =  \mathrm{ess}\liminf_{y \rightarrow x} v(y).$$
\end{thm}

\emph{Proof:} Fix an arbitrary point $x \in \Omega$. We must show only
that
$$\lambda =  \mathrm{ess}\liminf_{y \rightarrow x} v(y) \leq v(x),$$
since the opposite inequality was clear. Given any $\varepsilon > 0$,
there is a $\delta$ such that $v(y) > \lambda - \varepsilon$ for a.e.
$y \in B(x,\delta).$ By the lemma $v(y) \geq \lambda - \varepsilon$
for \emph{each} such $y$. In particular, $v(x) \geq \lambda -
\varepsilon$.
Because $\varepsilon$ was arbitrary, we have established that $v(x)
\geq \lambda$.  $\Box$

\subsection{The Evolutionary Equation}

We turn to the pointwise behaviour for the Evolutionary $p$-Laplacian
Equation. At each point in its domain a lower semicontinuous function satisfies

$$ v(x,t) \leq \liminf_{(y,\tau) \rightarrow (x,t)} v(y,\tau) \leq \mathrm{ess}\!
\liminf_{(y,\tau) \rightarrow (x,t)} v(y,\tau) \leq \mathrm{ess}\!
\liminf_{\underset{\tau < t}{(y,\tau) \rightarrow (x,t)}} v(y,\tau).$$
We show that for a $p$-supercaloric function also the reverse
inequalities hold, thus establishing Theorem \ref{Ess} in the Introduction.
In principle, the proof is similar to the stationary case, but now a
delicate issue of regularization arises. We first consider a
non-positive  $p$-supercaloric function
 $v=v(x,t)$ which is equal to zero at
\emph{almost} each point and, again, we show that locally it coincides with
the $p$-caloric function having the same boundary values, now in a space-time
cylinder. Then one has to conclude that $v$ was identically zero.

We seize the opportunity to describe a useful procedure of
\textsf{regularizing} by taking the
convolution\footnote{The origin of this function is unknown to me. In
  connexion with the Laplace transform it would be the convolution of
  $u$ and $\sigma^{-1}e^{-t/\sigma}$.}
$$u^{\star}(x,t) =
\frac{1}{\sigma}\int_{0}^{t}e^{(s-t)/\sigma}u(x,s)\,ds, \quad \sigma >
0.
$$
The notation hides the dependence on the parameter $\sigma$. For
continuous and for bounded semicontinuous functions $u$ the averaged
function $u^{\star}$ is defined at each point. We will stay within
this framework. Observe that
$$\sigma \frac{\partial u^{\star}}{\partial t} +u^{\star} = u.$$ Some
of its properties are listed in the next lemma.

\begin{lemma}\label{luutio}
\begin{description}
\item{(i)} If $u \in L^{p}(D_{T})$, then $$\|u^{\star}\|_{L^{p}(D_{T})}
 \leq \|u\|_{L^{p}(D_{T})}$$ and
$$ \frac{\partial u^{\star}}{\partial t} = \frac{u-u^{\star}}{\sigma}
 \in L^{p}(D_{T}).$$ Moreover, $u^{\star} \rightarrow u$ in
$L^{p}(D_{T})$ as $\sigma \rightarrow 0$.
\item{(ii)} If, in addition, $\nabla u \in L^{p}(D_{T})$, then $\nabla
  (u^{\star}) = (\nabla u)^{\star}$ componentwise,
$$\|\nabla u^{\star}\|_{L^{p}(D_{T})}
 \leq \|\nabla u\|_{L^{p}(D_{T})},$$
and $\nabla  u^{\star} \rightarrow \nabla u$ in $L^{p}(D_{T})$ as
$\sigma \rightarrow 0$.
\item{(iii)} Furthermore, if $u_{k} \rightarrow u$ in $ L^{p}(D_{T})$ 
  then also 
$$ u^{\star}_{k} \rightarrow  u^{\star}  \quad \text{and} \quad
 \frac{\partial u_{k}^{\star}}{\partial t} \rightarrow 
 \frac{\partial u^{\star}}{\partial t}$$
in $L^{p}(D_{T})$.
\item{(iv)} If $\nabla u_{k} \rightarrow \nabla u$ \ in \  $ L^{p}(D_{T})$,
  then $\nabla u^{\star}_{k} \rightarrow \nabla u^{\star}$\ in \ $
  L^{p}(D_{T})$.
\item{(v)} Finally, if $\varphi \in C(\overline{D_{T}})$, then 
$$ \varphi^{\star}(x,t) + e^{-t/\sigma}\varphi(x,0) \rightarrow
\varphi(x,t)$$
uniformly in $D_{T}$ as $\sigma \rightarrow 0$.
\end{description}
\end{lemma}

\emph{Proof:} We leave this as an exercise. (Some details are worked
out on page 7 of [KL1].)  $\Box$

\medskip

The \emph{averaged equation} for a weak supersolution $u$ in $D_{T}$
reads as follows:
\begin{gather*}
\int_{0}^{T}\!\!\int_{D}\Bigl(\langle(|\nabla u|^{p-2}\nabla u)^{\star},\nabla
  \varphi \rangle - u^{\star}\frac{\partial \varphi}{\partial
    t}\Bigr)dx\,dt + \int_{D}u^{\star}(x,T)\varphi(x,T)\,dx\\
\geq \int_{D} u(x,0)
\left(\frac{1}{\sigma}\int_{0}^{T}\varphi(x,s)e^{-s/\sigma}\,ds
\right)dx
\end{gather*}
valid for all test functions $\varphi \geq 0$ vanishing on the
lateral boundary $\partial D\times [0,T]$ of the space-time cylinder. For
solutions one has equality. Notice the typical difficulty with
obtaining $(|\nabla u|^{p-2}\nabla u)^{\star}$ and not $|\nabla
u^{\star}|^{p-2} \nabla u^{\star}$, except in the linear case. The
averaged equation follows from the equation for the retarded
supersolution $u(x,t-s)$, where $0 \leq s \leq T$:
\begin{gather*}
\int_{s}^{T}\!\!\int_{D}\Bigl(\langle|\nabla u(x,t-s)|^{p-2}\nabla
  u(x,t-s),\nabla \varphi(x,t)\rangle - u(x,t-s)\frac{\partial
    \varphi}{\partial t}(x,t)\Bigr)dx\,dt\\
+ \int_{D}u(x,T-s)\varphi(x,T)\,dx \geq \int_{D}u(x,0)\varphi(x,s)\,dx.
\end{gather*}
Notice that $(x,t-s) \in \overline{D_{T}}$ when $0 \leq s \leq t\leq
T$. Multiply by $\sigma ^{-1}e^{-s/\sigma}$, integrate over $[0,T]$
with respect to $s$, and, finally, interchange the order of
integration between $s$ and $t$. This yields the averaged equation
above.

The advantage of this procedure over more conventional convolutions is
that no values outside the original space-time cylinder are evoked.

We begin with a simple situation.

\begin{lemma} Suppose that $v$ is a $p$-supercaloric function in a
  domain containing the closure of $B_{T} = B \times (0,T)$. If
\begin{description}
\item{(i)\phantom{i}} $v \leq 0$ at each point in $B_{T}$ and
\item{(ii)} $v = 0$ at almost every point in  $B_{T}$,
\end{description}
then $v = 0$ at each point in  $B \times (0,T]$.
\end{lemma}

\medskip
\emph{Proof:} We may assume that $v$ is bounded. Construct the infimal
convolution $v_{\varepsilon}$ with respect to a larger domain than
$B_{T}$. Fix a small time $t' > 0$ and let $h^{\varepsilon}$ be the
$p$-caloric function with boundary values induced by
$v_{\varepsilon}$ on the parabolic boundary of the cylinder $B \times
(t',T)$ and define the function
\begin{equation*}
w_{\varepsilon} = 
\begin{cases}
h^{\varepsilon},\quad \text{in} \quad B \times
(t',T]\\
v_{\varepsilon} \quad \text{otherwise}.
\end{cases}
\end{equation*}
To be on the safe side concerning the validity at the terminal time
$T$ we may solve the boundary value problem in a slightly larger
domain with terminal time $T' > T$. Also $w_{\varepsilon}$ is a
$p$-supercaloric function. By comparison
$$w_{\varepsilon} \leq v_{\varepsilon} \leq 0 \quad \text{pointwise in}
\quad B_{T}.$$
We let $\varepsilon$ go to zero through a monotone sequence, say
$1,\frac{1}{2},\frac{1}{3}, \cdots .$ Then the limit 
$$h = \lim_{\varepsilon \rightarrow 0}h^{\varepsilon}$$
exists pointwise and it follows from the uniform H\"{o}lder estimates
(\ref{Holder}) that this $h$ is continuous without any correction made in a set
of measure zero. It is important to preserve the information at
\emph{each} point. Thus $h$ is a $p$-caloric function.  The so
obtained function
\begin{equation*}
w = 
\begin{cases}
h,\quad \text{in} \quad B \times
(t',T')\\
v \quad \text{otherwise}
\end{cases}
\end{equation*}
is a $p$-supercaloric function. For the verification of the
semicontinuity and the comparison
principle, which proves this, the fact that $h \leq v$ is essential.

We know that $w \leq v \leq 0$ everywhere in  a domain containing 
$B \times (0,T)$. In particular,
 $$h \leq v \leq 0 \quad \text{everywhere in}\quad
 B \times (0,T).$$
 We claim that $h = 0$ at each point. The claim
 immediately implies that $v = 0$ at each point in  $B \times (0,T)$.
Concerning the statement at the terminal time $T$, we notice that $v
\geq h$ and
$$v(x,T) \geq h(x,T)
= \lim_{\phantom{a}t \rightarrow T-}h(x,t) = 0,$$
since $h$ is continuous. On the other hand $v(x,T) \leq 0$ by the
lower semicontinuity. Thus also $v(x,T) = 0$.

Therefore it is sufficient to prove the claim. To conclude that $h$ is
identically zero we use the averaged
equation for $w^{\star}$ and write
\begin{gather*}
\int_{0}^{T}\!\!\int_{B}\Bigl(\langle(|\nabla w|^{p-2}\nabla w)^{\star},\nabla
  \varphi \rangle + \varphi \frac{\partial w^{\star}}{\partial
    t}\Bigr)dx\,dt \\
\geq \int_{B} w(x,0)
\left(\frac{1}{\sigma}\int_{0}^{T}\varphi(x,s)e^{-s/\sigma}\,ds
\right)dx,
\end{gather*}
where the test function vanishes on the parabolic boundary (an
integration by parts has been made with respect to time.)
Select the test function $\varphi = (v_{\varepsilon} -
w_{\varepsilon})^{\star}$ and let $\varepsilon$ approach zero. Taking
into account that $\varphi = 0$ when $t<t'$, we
arrive at
\begin{gather*}
\int_{t'}^{T}\!\!\int_{B}\Bigl(\langle(|\nabla h|^{p-2}\nabla
  h)^{\star},\nabla v^{\star}- \nabla
  h^{\star} \rangle + (v^{\star}-  h^{\star})\frac{\partial h^{\star} }{\partial
    t}\Bigr)dx\,dt\\
\geq \int_{B} v(x,0)
\left(\frac{1}{\sigma}\int_{t'}^{T}(v^{\star}(x,s)-h^{\star}(x,s))e^{-s/\sigma}\,ds
\right)dx.
\end{gather*}

The last integral (which could be negative) approaches zero as the
 regularization parameter
$\sigma$ goes  to zero, because $t' > 0$, so that the exponential decays.
Integrating
$$(v^{\star}-  h^{\star})\frac{\partial h^{\star} }{\partial
    t} = - (v^{\star}-  h^{\star})\frac{\partial (v^{\star}-
    h^{\star}) }{\partial
    t} + (v^{\star}-  h^{\star})\frac{\partial v^{\star} }{\partial
    t}$$
we obtain
\begin{gather*}
\int_{t'}^{T}\!\!\int_{B}(v^{\star}-  h^{\star})
\frac{\partial h^{\star} }{\partial
    t}\,dx\,dt\\ = -\frac{1}{2}\int_{B}(v^{\star}(x,T)- 
 h^{\star}(x,T))^{2}\,dx + \int_{t'}^{T}\!\!\int_{D}(v^{\star}-  h^{\star})
\frac{\partial v^{\star} }{\partial
    t}\,dx\,dt.
\end{gather*}
because $v^{\star}(x,t')-h^{\star}(x,t') = 0.$
The last integral is zero because $v^{\star}$ and $\frac{\partial
  v^{\star}}{\partial t}$ are zero almost everywhere according to
property (i) in Lemma \ref{luutio}. Erasing this integral and letting the
regularization parameter $\sigma$ go to zero (so that the $\star$'s
disappear) we finally obtain
$$\int_{t'}^{T}\!\!\int_{B}|\nabla h|^{p}\,dx\,dt
+\frac{1}{2}\int_{B}h^{2}(x,T)\,dx \,\leq \,0\quad \text{i.e.}\quad = 0.$$
In fact\footnote{It is the validity of$$
\lim_{\sigma \rightarrow 0}\frac{1}{2}\int_{B}(v^{\star}(x,T)- 
 h^{\star}(x,T))^{2}\,dx =\frac{1}{2}\int_{B}h^{2}(x,T)\,dx$$ that
 requires some caution. We know that $v^{\star}$ is zero almost
 everywhere but with respect to the $(n+1)$-dimensional measure.}
, the proof guarantees this only for almost all values of $T$
in the range $t' < T < T'$.
From this it is not difficult to conclude that $h$ is identically
zero. Thus our claim has been proved. $\Box$

\begin{lemma}
Suppose that $v$ is a $p$-supercaloric function in a domain containing
$B_{T} = B \times (0,T)$. If $v(x,t)> \lambda$ for almost every $(x,t)
\in B_{T}$, then $v(x,t)\geq \lambda$ for every $(x,t)
\in B \times (0,T]$.
\end{lemma}

\medskip
\emph{Proof:} The auxiliary function
$$u(x,t) = \min\{v(x,t),\lambda\} - \lambda$$
in place of $v$ satisfies the assumptions in the previous lemma. Hence
$u= 0$ everywhere in  $B \times (0,T]$. This is equivalent to the
assertion. $\Box$

\medskip
\emph{Proof of Theorem \ref{Ess}:} Denote
$$\lambda = \mathrm{ess}\!
\liminf_{\underset{t < t_{0}}{(x,t) \rightarrow (x_{0},t_{0})}}
v(x,t).$$
According to the discussion in the beginning of this subsection, it is
sufficient to prove that $\lambda \leq  v(x_{0},t_{0})$. Thus we can
assume that $\lambda > - \infty$.

 First, we consider the case $\lambda <
\infty$. Given $\varepsilon > 0$, we can find a $\delta > 0$ and a
ball $B$ with centre $x_{0}$ such that the closure of $B \times
(t_{0}- \delta,t_{0})$ is comprised in the domain and
$$v(x,t) > \lambda - \varepsilon$$
for almost every $(x,t) \in B \times
(t_{0}- \delta,t_{0})$. According to the previous lemma 
$$v(x,t) \geq  \lambda - \varepsilon$$
for \emph{every}  $(x,t) \in B \times
(t_{0}- \delta,t_{0}]$. In particular, we can take $(x,t) =
(x_{0},t_{0})$.
Hence $v(x_{0},t_{0}) \geq \lambda - \varepsilon$. Since $\varepsilon$
was arbitrary, we have proved that $\lambda \leq  v(x_{0},t_{0})$, as
desired.

Second, the case $\lambda = \infty$ is easily reached via the
truncated functions $v_{k} = \min\{v(x,t),k\},\,\, k = 1,2,\ldots.$
Indeed,
$$v(x_{0},t_{0}) \geq v_{k}(x_{0},t_{0}) \geq \min\{\infty,k\} = k,$$
in view of the previous case. This concludes the proof of Theorem \ref{Ess}.\qquad$\Box$

\section{Viscosity Supersolutions Are Weak\\
  Supersolutions}

In this chapter\footnote{The previous chapters, do in fact, not
    require familiarity with the viscosity theory of second order
    equations, but now it is desirable that the reader knows the
    basics of this theory. Some chapters of Koike's book [Ko] are
    enough. A more advanced source is [CIL].} we give a simple proof, due to Julin and Juutinen, of
the fact that the viscosity supersolutions are the same as those
obtained in potential theory, cf. [JJ]. The proof in [JLM], which is more
complicated, will be bypassed. (Thus we can avoid the uniqueness
machinery for second order equations, the doubling of variables, and
Jensen's auxiliary equations.) The proof is based on the fact that the
\emph{infimal convolutions have second derivatives in the sense of
Alexandrov}, which can be used in the testing with so-called superjets.
These occur in a reformulation of the definition of viscosity
supersolutions. The idea is that a \emph{one}-sided estimate  makes it possible to use
Fatou's lemma and finally
pass to the limit in an integral. 

In establishing the equivalence between the two concepts of
supersolutions, the easy part is to show that $p$-superharmonic or
$p$-supercaloric functions are viscosity supersolutions. The proof
comes from the fact that an antithesis produces a touching test
function which is $p$-\emph{sub}harmonic or $p$-\emph{sub}caloric in a
neighbourhood, in which situation the comparison principle leads to a
contradiction in the indirect proof. This was accomplished in the
proof of Proposition \ref{uusi} for the stationary case. The evolutionary case
is similar, and we omit it here. ---We now turn to the preparations
for the difficult part of the equivalence proof.

\begin{thm}[Alexandrov]
Let $f: \mathbf{R}^n\rightarrow \mathbf{R}$ be a convex function. Then $f$ has second
derivatives in the sense of Alexandrov: for a. \!e. point $x$
there is a symmetric $n\times n$-matrix $\mathbb{A} = \mathbb{A}(x)$
such that the expansion
$$f(y) = f(x) + \langle\nabla f(x),y-x\rangle + \frac{1}{2}\langle
y-x, \mathbb{A}(x)(y-x)\rangle + o(|y-x|^2)$$
is valid as $y \rightarrow x$.
\end{thm}

For a proof\footnote{Some details in [GZ, Lemma 7.11, p. 199] are helpful to
  understand the singular part of the Lebesgue decomposition, which is
  used in
  the proof in [EG].} we refer to [EG, Section 6.4, pp. 242--245]. The
problem is not the
first derivatives, since by Rademacher's Theorem they  are Sobolev derivatives and $\nabla f \in
L^{\infty}_{\mathrm{loc}}$. The question is about the second ones. We will use
the \emph{notation} $D^2f = \mathbb{A}$, although the Alexandrov derivatives
are not always second Sobolev derivatives, because a singular Radon
measure may be present. The proof in [EG] establishes that pointwise
we have a.\!\! e. that
\begin{equation}
\label{evans}
\mathbb{A} = \lim_{\varepsilon \rightarrow 0}(D^{2}(f\star
\varrho_{\varepsilon}))
\end{equation}
where $\varrho_{\varepsilon}$ is Friedrich's mollifier.

Alexandrov's theorem is applicable to the concave 
functions 
$$v_{\varepsilon}(x) - \frac{|x|^{2}}{2 \varepsilon}, \qquad
v_{\varepsilon}(x,t) - \frac{|x|^{2}+ t^{2}}{2 \varepsilon}$$
encountered in Section 2 and Section 3.2. They are in
fact defined in the whole space (although the infima are taken over
bounded sets). Then the theorem is also applicable to the
$v_{\varepsilon}$, since the subtracted smooth functions have second
derivatives.

\subsection{The Stationary Case}

The concept of viscosity solutions can be reformulated in terms of
so-called jets.  \emph{Super}solutions require \emph{sub}jets (and
\emph{sub}solutions \emph{super}jets). We say that the pair $(\xi,\mathbb{X})$, where $\xi$
is a vector in $\mathbf{R}^n$ and $\mathbb{X}$ is a symmetric $n\times
n$-matrix, belongs to the \emph{subjet} $J^{2,-}u(x)$ if
$$u(y) \geq u(x)   + \langle\xi,y-x\rangle + \frac{1}{2}\langle
y-x, \mathbb{X}(x)(y-x)\rangle + o(|y-x|^2)$$
as $y \rightarrow x$. See [Ko, Section 2.2, p. 17]. Notice the similarity with 
a Taylor polynomial. If it so happens that $u$ has continuous
second derivatives at the point $x$, then we must have $\xi = \nabla
u(x),\,  \mathbb{X} = D^2 u(x) = $ the Hessian matrix. In other words,
$$J^{2,-}u(x) = \{(\nabla u(x), D^2 u(x))\}.$$ The essential feature is
that the Alexandrov derivatives always do as members of the jets.

For a wide class of second order equations the subjets can be used to give an equivalent characterization of the
viscosity supersolutions. We need only the following
necessary\footnote{Testing with subjets is also a sufficient condition
  when their ''closures'' are employed.} 
condition.

\begin{prop}
Let $p \geq 2$. Suppose that $\Delta_pv \leq 0$ in the viscosity
sense. If $(\xi, \mathbb{X}) \in J^{2,-}v(x)$, then
\begin{equation}
\label{eq:jet}
|\xi|^{p-2}\mathsf{trace}(\mathbb{X}) + (p-2)|\xi|^{p-4}\langle \xi,\mathbb{X}\,
\xi \rangle \, \leq \, 0.
\end{equation}
\end{prop}

\emph{Proof:} A simple proof is given in [Ko, Proposition 2.6, pp. 18--19].
$\Box$

\medskip
After these preparations we are in the position of proving that a
bounded viscosity supersolution of the Stationary $p$-Laplace
Equation, $p \geq 2$, is also a weak supersolution. This is the analogue of Theorem \ref{boundell}
in Section 2, but for viscosity supersolutions. It was based on Corollary \ref{mant}. We will
now prove Corollary \ref{mant} for viscosity supersolutions without evoking the reference [JLM]. To this
end, assume that $0 \leq v(x) \leq L$ and that $\Delta_p v(x) \leq 0$
in the viscosity sense in $\Omega$. The infimal convolution
$v_{\varepsilon}$ defined by formula (\ref{ellinfimal}) is, according to Proposition \ref{sitten}, also a viscosity supersolution
in the shrunken domain $\Omega_{\varepsilon}.$ Given a non-negative
test function $\psi$ in $C_{0}^{\infty}(\Omega)$, we have to prove the
following 
$$\mathsf{Claim:}\qquad \qquad  
\int_{\Omega_{\varepsilon}}\big\langle |\nabla
v_{\varepsilon}|^{p-2}\nabla v_{\varepsilon},\nabla \psi\big\rangle\,dx
\geq 0$$
when $\varepsilon$ is so small that $\Omega_{\varepsilon}$ contains
the support of $\psi$.

  As we saw above, the
second Alexandrov derivatives $D^2v_{\varepsilon}(x)$ exist a.e. in
$\Rn$ and therefore $\bigl(\nabla v_{\varepsilon}(x), D^2
v_{\varepsilon}(x)\bigr) \in J^{2,-}v_{\varepsilon}(x)$ at almost every point
$x$. Hence, by the Proposition, the inequality
\begin{gather}
\label{eq:A}
\Delta_{p}v_{\varepsilon}(x) \nonumber\\
=\,\,  |\nabla
v_{\varepsilon}(x)|^{p-4}\Bigl\{|\nabla v_{\varepsilon}(x)|^2\Delta
v_{\varepsilon}(x)+(p-2)\big\langle \nabla  v_{\varepsilon}(x),
D^2v_{\varepsilon}(x)\nabla  v_{\varepsilon}(x)\big\rangle \Bigr\} \nonumber\\
\leq\,\,  0
\end{gather}
is valid almost everywhere in $\Omega_{\varepsilon}.$   Here $\Delta
v_{\varepsilon} = \mathsf{trace}(D^2v_{\varepsilon}) .$

We need one further mollification. For
$$f_{\varepsilon}(x) = v_{\varepsilon}(x) -
\frac{|x|^2}{2\varepsilon}$$
we define the convolution
$$f_{\varepsilon,j} = f_{\varepsilon}\star
\varrho_{\varepsilon_{j}} \quad
\text{where}\quad
\varrho_{\varepsilon_{j}} =
\begin{cases}\frac{C}{\varepsilon^n_{j}}\exp\bigl(-\frac{\varepsilon^2_{j}}{\varepsilon_{j}^2-|x|^2}\bigr),\quad \text{when}\,\,
  |x|<\varepsilon_{j}\\
0, \quad\text{otherwise}.
\end{cases}$$
The smooth functions $v_{\varepsilon,j} = v_{\varepsilon}\star
\varrho_{\varepsilon_{j}}$ satisfy the identity
$$\int_{\Omega_{\varepsilon}}\big\langle |\nabla  v_{\varepsilon,j}
|^{p-2}v_{\varepsilon,j},\nabla \psi \big\rangle \, dx =
\int_{\Omega_{\varepsilon}}\psi\bigl(-\Delta_p
v_{\varepsilon,j}\bigr)\,dx,$$
which identity we  will extend to the function $v_{\varepsilon}$ by
passing to the limit. However, they are not viscosity supersolutions
themselves! 
By the linearity of the convolution,  we can from (\ref{evans})  conclude that
$$\lim_{j\rightarrow \infty}D^{2}v_{\varepsilon,j} =
D^{2}v_{\varepsilon}$$
almost everywhere. Therefore we have 
$$\lim_{j\rightarrow \infty}\Delta_p v_{\varepsilon,j}(x) = \Delta_p
v_{\varepsilon}(x)$$ at a.\!\! e. point $x$ in the support of $\psi.$
Obviously, the convolution has preserved the concavity, and hence
$D^2f_{\varepsilon,j} \leq 0.$ It follows that
$$D^2v_{\varepsilon,j} \leq \frac{\mathbb{I}_{n}}{\varepsilon}, \qquad \Delta
v_{\varepsilon,j} \leq \frac{n}{\varepsilon}$$
a.\!\! e.. Here  $\mathbb{I}_{n}$ is the unit matrix. It is immediate that
$$|\nabla v_{\varepsilon,j}| \leq \|\nabla  v_{\varepsilon}\|_{\infty}
= C_{\varepsilon}.$$ These estimates yield
 the  bound
\begin{equation}
\label{eq:fatou}
- \Delta_p v_{\varepsilon,j} \geq -
C_{\varepsilon}^{p-2}\frac{n+p-2}{\varepsilon}
\end{equation}
valid almost everywhere in the support of $\psi$. This lower bound
justifies the use of Fatous lemma below:
\begin{align*}
\int_{\Omega_{\varepsilon}}\big\langle |\nabla
v_{\varepsilon}|^{p-2}\nabla v_{\varepsilon},\nabla \psi\big \rangle\,dx 
&= \lim_{j \rightarrow \infty}\int_{\Omega_{\varepsilon}}\big\langle |\nabla
v_{\varepsilon,j}|^{p-2}\nabla v_{\varepsilon,j},\nabla \psi \big\rangle\,dx\\
&=  \lim_{j \rightarrow \infty}  \int_{\Omega_{\varepsilon}}\psi
(-\Delta_{p}v_{\varepsilon,j})\, dx\\
&\geq  \int_{\Omega_{\varepsilon}} \liminf_{j \rightarrow
  \infty}  \Bigl(\psi(-\Delta_{p}v_{\varepsilon,j})\Bigr)\,dx\\
&= \int_{\Omega_{\varepsilon}}\psi (-\Delta_{p}v_{\varepsilon})\,dx\quad
\geq \quad 0.
\end{align*}
 In the very last step we used the inequality
 $-\Delta_{p}v_{\varepsilon} \geq 0$, which, as we recall, needed Alexandrov's
 theorem in its proof.
This proves our claim. $\Box$

\subsection{The Evolutionary Equation}

Since the parabolic proof is very similar to the elliptic one, we only
sketch the proof of the 
$$\mathsf{Claim:} \qquad  
\int_{0}^{T}\!\!\int_{\Omega_{\varepsilon}}\Bigl(-v_{\varepsilon} \frac{\partial \psi}{\partial t} +\big\langle |\nabla
v_{\varepsilon}|^{p-2}\nabla v_{\varepsilon},\nabla \psi \big\rangle\!
\Bigr)dx\,dt \geq 0$$
for all non-negative test functions $\psi \in C_0^{\infty}(\Omega).$
As in Section 3.2 the infimal convolution of the given bounded viscosity
supersolution $v$, $0 \leq v(x,t) \leq L,$ is defined as
$$v_{\varepsilon}(x,t) = \inf_{(y,\tau)\in \Omega_{T}} \Bigl\{v(y,\tau) + \frac{|y-x|^{2} +
  (\tau -t)^{2}}{2
\varepsilon}\Bigr\}$$
and the function
$$f_{\varepsilon}(x,t) =v_{\varepsilon}(x,t) - \frac{|x|^{2} +
  t^{2}}{2\varepsilon}$$
is  concave  in $n+1$ variables. Therefore it has second
derivatives in the sense of Alexandrov. So has $v_{\varepsilon}$
itself, since the quadratic term has no influence on this matter. Thus
\begin{align*}
v_{\varepsilon}(&y,\tau)\\ & = v_{\varepsilon}(x,t) + \langle \nabla
v_{\varepsilon}(x,t),y-x \rangle  + \frac{1}{2} \langle
y-x,D^2_{x}v_{\varepsilon}(x,t)( y-x) \rangle \\ & + \frac{\partial
  v_{\varepsilon}(x,t)}{\partial t}(t-\tau)
+ \Bigl\langle \nabla  \frac{\partial
  v_{\varepsilon}(x,t)}{\partial t},y-x \Bigr\rangle (\tau-t) +  \frac{1}{2} \frac{\partial^2
  v_{\varepsilon}(x,t)}{\partial t ^{2}}(\tau-t)^2\\
&+  o\bigl(|y-x|^2 +
|\tau-t|^{2}\bigr)
\end{align*}
as $(y,\tau) \rightarrow (x,t).$ Here $D^2_xv_{\varepsilon}$ is not
the complete Hessian but the $n\times n$-matrix consisting of the second space
derivatives; the time derivatives are separately written. This implies
that
\begin{align*}
v_{\varepsilon}(y,\tau)& = v_{\varepsilon}(x,t)  + \frac{\partial
  v_{\varepsilon}(x,t)}{\partial t}(t-\tau) +  \langle \nabla
v_{\varepsilon}(x,t),y-x \rangle \\
  & + \frac{1}{2} \langle
y-x,D^2_{x}v_{\varepsilon}(x,t)( y-x) \rangle 
+  o\bigl(|y-x|^2 +
|\tau-t|\bigr),
\end{align*}
where the error term is no longer quadratic in $\tau -t.$

The \emph{parabolic subjet} $\mathcal{P}^{2,-}u(x,t)$ consists of all
triples $(a,\xi,\mathbb{X}),$ where $a = a(x,t)$ is a real number, $\xi
= \xi(x,t)$ is a vector in $\Rn$ and $\mathbb{X} = \mathbb{X}(x,t),$
such that 
\begin{align*}
u(y,\tau) \geq & \,u(x,t) + a\,(\tau -t) + \langle \xi,y-x \rangle +
\frac{1}{2} \langle y-x,\mathbb{X}(y-x) \rangle\\ + &\, o\bigl(|y-x|^2 +
|\tau-t|\bigr)
\end{align*}
as $(y,\tau) \rightarrow (x,t).$ See [CIL, equation(8.1), p.48]. The Alexandrov (and
the Rademacher) derivatives will do in the parabolic subjet and the
characterization of viscosity supersolutions in terms of jets yields
now
the pointwise inequality
\begin{equation}
\label{eq:pos}
 -\Delta_{p}v_{\varepsilon} +  \frac{\partial
  v_{\varepsilon}}{\partial t} \geq 0
\end{equation}
valid almost everywhere in the support of $\psi.$ 

Again we have to use a convolution. Because the second time
derivatives will not be needed, we take the convolution
$f_{\varepsilon,j} =f_{\varepsilon} \star \varrho_{\varepsilon_{j}}$
only with respect to the space variables:  $ \varrho_{\varepsilon_{j}}
= \varrho_{\varepsilon_{j}}(x).$ (This does not matter.) We have
$$\int \!\!\int\Bigl( -v_{\varepsilon,j} \frac{\partial
  \psi}{\partial t} + \big\langle |\nabla
v_{\varepsilon,j}|^{p-2}v_{\varepsilon,j},\nabla \psi \big\rangle \!\Bigr) \,
dx\,dt = \int \!\!\int \psi \Bigl(\frac{\partial v_{\varepsilon,j}
  }{\partial t} - \Delta_p v_{\varepsilon,j}\Bigr)\,dx\,dt,$$
where the integrals are taken over the support of $\psi$, and $\varepsilon$
is small. One can clearly pass to the limit under the integral signs
above, as $j \rightarrow \infty$, except that the integral of $
 - \Delta_p v_{\varepsilon,j}$ requires a
  justification. Actually, the estimate (\ref{eq:fatou}) is valid also
  in the parabolic case, whence Fatou's lemma can be used. We obtain
\begin{align*}
\int \!\!\int\Bigl( -v_{\varepsilon} \frac{\partial
  \psi}{\partial t} +& \big\langle |\nabla
v_{\varepsilon}|^{p-2}\nabla v_{\varepsilon},\nabla \psi \big\rangle \!\Bigr) \,
dx\,dt
\\ = &\, \lim_{j \rightarrow \infty}\int \!\!\int\Bigl( -v_{\varepsilon,j} \frac{\partial
  \psi}{\partial t} + \big\langle |\nabla
v_{\varepsilon,j}|^{p-2}v_{\varepsilon,j},\nabla \psi \big\rangle \!\Bigr) \,
dx\,dt \\= &\, \lim_{j \rightarrow \infty} \int \!\!\int \psi \Bigl(\frac{\partial v_{\varepsilon,j}
  }{\partial t} - \Delta_p v_{\varepsilon,j}\Bigr)\,dx\,dt
\\\geq &\, \int \!\!\int \liminf_{j \rightarrow \infty} \psi \Bigl (\frac{\partial v_{\varepsilon,j}
  }{\partial t} - \Delta_p v_{\varepsilon,j}\Bigr)\,dx\,dt
\\= &\, \int \!\!\int \psi \Bigl (\frac{\partial v_{\varepsilon}
  }{\partial t} - \Delta_p v_{\varepsilon}\Bigr)\,dx\,dt \quad
\geq \quad 0,
\end{align*}
where we used (\ref{eq:pos}) in the last step. This proves our claim. $\Box$
$$\Bigg\langle \text{THE END} \Bigg\rangle$$

\end{document}